\documentclass[english]{elsarticle}
\usepackage[T1]{fontenc}
\usepackage[latin9]{inputenc}
\usepackage{amsmath}
\usepackage{amsthm}
\usepackage{amssymb}
\usepackage{graphicx}

\makeatletter
\renewcommand{\vec}[1]{\boldsymbol #1}

\renewcommand\[{\begin{equation}}
\renewcommand\]{\end{equation}}

\@ifundefined{showcaptionsetup}{}{%
 \PassOptionsToPackage{caption=false}{subfig}}
\usepackage{subfig}
\makeatother

\usepackage{babel}
\begin{document}

\begin{frontmatter}{}

\title{WLS-ENO Remap: Superconvergent and Non-Oscillatory \\
Weighted Least Squares Data Transfer on Surfaces}

\author{Yipeng Li}

\author{Qiao Chen}

\author{Xuebin Wang}

\author{Xiangmin Jiao\corref{cor1}}

\ead{xiangmin.jiao@stonybrook.edu}

\cortext[cor1]{Corresponding author. }

\address{Department of Applied Mathematics \& Statistics and Institute for
Advanced Computational \\
Science, Stony Brook University, Stony Brook, NY 11794, USA}
\begin{abstract}
Data remap between non-matching meshes is a critical step in multiphysics
coupling using a partitioned approach. The data fields being transferred
often have jumps in function values or derivatives. It is important
but very challenging to avoid spurious oscillations (a.k.a. the \emph{Gibbs
Phenomenon}) near discontinuities and at the same time to achieve
high-order accuracy away from discontinuities. In this work, we introduce
a new approach, called \emph{WLS-ENOR}, or \emph{Weighted-Least-Squares-based
Essentially Non-Oscillatory Remap}, to address this challenge. Based
on the WLS-ENO reconstruction technique proposed by Liu and Jiao (\emph{J.
Comput. Phys.} vol 314, pp 749-{}-773, 2016),\emph{ WLS-ENOR} differs
from WLS-ENO and other WENO schemes in that it resolves not only the
$\mathcal{O}(1)$ oscillations due to $C^{0}$ discontinuities, but
also the accumulated effect of $\mathcal{O}(h)$ oscillations due
to $C^{1}$ discontinuities. To this end, WLS-ENOR introduces a robust
detector of discontinuities and a new weighting scheme for WLS-ENO
near discontinuities. We also optimize the weights at smooth regions
to achieve superconvergence. As a result, WLS-ENOR is more than fifth-order
accurate and highly conservative in smooth regions, while being non-oscillatory
and minimally diffusive near discontinuities. We also compare WLS-ENOR
with some commonly used methods based on $L^{2}$ projection, moving
least squares, and radial basis functions.
\end{abstract}
\begin{keyword}
weighted least squares; data transfer; superconvergence; discontinuities;
Gibbs phenomenon; essentially non-oscillatory scheme
\end{keyword}

\end{frontmatter}{}

\section{Introduction}

In multiphysics applications with a partitioned approach, data often
must be mapped between different meshes \citep{keyes2013multiphysics}.
This problem is often referred to as \emph{data transfer} or \emph{remap},
and the two meshes are often referred to as the \emph{source} (or
\emph{donor}) and \emph{target} (or \emph{donee}) mesh, respectively.
In many of these applications, the data approximate a piecewise smooth
function, where discontinuities (or jumps) may be present in the solutions,
for example, due to shocks in high-speed fluids or multi-material
solids. In addition, nowadays the physics models are often discretized
with third or higher-order methods, such as $hp$ finite elements
\citep{karniadakis2013spectral}, spectral elements \citep{lauritzen2018ncar},
discontinuous Galerkin \citep{cockburn2000development}, etc. The
combination of the discontinuities and high-order methods lead to
a fundamental challenge for remap: high-order interpolations (such
as splines \citep{de1990quasiinterpolants}) or least squares (such
as moving least squares (MLS) \citep{wendland2001local}) are prone
to spurious oscillations near discontinuities \citep{fornberg2007gibbs,jerri2013gibbs}.
In particular, $C^{0}$ discontinuities, a.k.a. (\emph{simple}) \emph{jumps},
often lead to $\mathcal{O}(1)$ oscillations (overshoots or undershoots),
analogous to the \emph{Gibbs-Wilbraham phenomenon} in Fourier transformation
\citep{gibbs1899,wilbraham1848certain,jerri2013gibbs}. In addition,
$C^{1}$ discontinuities, a.k.a. \emph{derivative} \emph{jumps}, can
lead to $\mathcal{O}(h)$ oscillations, which may accumulate in repeated
transfer or in time-dependent settings. Even some second and lower
order methods, such as $L^{2}$ projection \citep{jiao2004common}
and radial basis function (RBF) interpolation \citep{beckert2001multivariate,flyer2007transport,rendall2008unified}
can suffer from such oscillations \citep{fornberg2007gibbs,foster1991gibbs}.
Among the existing methods, nearest-point or linear interpolation
have no oscillation, but they do not deliver sufficient accuracy and
are excessively diffusive in repeated transfer. The objective of this
work is to propose a new data-transfer method that is at least third-order
accurate for smooth functions, is non-oscillatory at discontinuities,
and is as non-diffusive as possible.

Resolution of the Gibbs phenomenon or similar oscillations is an important
and challenging subject in numerical analysis and computational physics.
There have been various methods developed for it in various contexts;
see e.g. \citep{gottlieb1997gibbs}. Among these, the most closely
related ones are probably the so-called \emph{reconstruction} methods,
such as flux reconstructions based on WENO \citep{liu1994weighted,shu2009high}
and their variants in finite volume methods. The connection is even
closer in the context of hyperbolic conservation laws, in that they
are time dependent. However, in the context of hyperbolic PDEs, it
typically suffices for the reconstruction to resolve only simple jumps,
because small oscillations due to $C^{1}$ discontinuities can be
effectively controlled by the numerical diffusion in flux limiters
\citep{leveque1992numerical} and TVD time-integration schemes \citep{gottlieb1998total}.
In contrast, a general remap technique must be stable without assuming
numerical diffusion in the physics models. Hence, it is important
to control the oscillations due to both simple and derivative jumps
in remap. This property makes remap significantly more challenging.

In this work, we introduce a new method, called \emph{WLS-ENOR}, or
\emph{Weighted-Least-Squares-based Essentially Non-Oscillatory Remap}.
WLS-ENOR is based on the WLS-ENO scheme proposed by Liu and Jiao \citep{liu2016wls}
for the finite volume methods. The original WLS-ENO controls oscillations
by adapting the weights in a weighted-least-squares (WLS) approximation,
where the weights are approximately based on the inverse of a \emph{non-smoothness}
\emph{indicator}. However, the indicator in WLS-ENO only considers
simple jumps. WLS-ENOR introduces a novel technique to detect both
simple and derivative jumps by introducing two novel indicators and
a new dual-thresholding strategy. We then apply WLS-ENO near the detected
singularities by using a new weighting scheme to take into account
$C^{1}$ discontinuities. We also optimize WLS-ENOR for smooth functions
to achieve superconvergence. As a result, WLS-ENOR is more than fifth-order
accurate for smooth functions, non-oscillatory and minimally diffusive
near $C^{0}$ and $C^{1}$ discontinuities, and is highly conservative
in repeated transfers. We compare WLS-ENOR with some commonly used
methods, including linear interpolation, RBF interpolation \citep{beckert2001multivariate},
$L^{2}$ projection \citep{jiao2004common}, and modified moving least
squares (MMLS) \citep{joldes2015modified,slattery2016mesh}, and show
that WLS-ENOR is more accurate, stable, and less diffusive in almost
all cases. The presentation and results in this work focus on spherical
geometries in earth and climate modeling. However, the methodology
also applies to other smooth surfaces.

The remainder of this paper is organized as follows. In Section~\ref{sec:preliminaries},
we review some preliminaries on remap and reconstruction methods and
the Gibbs phenomenon at discontinuities. In Section~\ref{sec:Superconvergent-weighted-least},
we describe WLS fittings for smooth functions that exhibit superconvergence.
In Section~\ref{sec:Detection-and-resolution}, we describe the detection
of discontinuities in data remapping and a new weighting scheme for
WLS-ENO for discontinuities. In Section~\ref{sec:Numerical-results},
we compare WLS-ENOR with other methods for transferring both smooth
and discontinuous functions. Finally, Section~\ref{sec:Conclusions}
concludes the paper with a discussion on future research directions.

\section{\label{sec:preliminaries}Preliminaries and related work}

In this section, we review some methods for remap and reconstruction,
especially the high-order methods based on weighted-least-squares
(WLS). In addition, we review some concepts and methods related to
discontinuities, including Gibbs phenomena and WLS-ENO.

\subsection{Remap and reconstruction methods}

The remap problem has been studied extensively in the past few decades.
Several software packages have been developed over the past two decades.
These include those developed for climate and earth modeling (such
as those in Earth System Modeling Framework (ESMF) \citep{hill2004architecture},
Community Earth System Model (CESM) \citep{hurrell2013community}
and the next-generation Energy Exascale Earth System Model (E3SM),
including SCRIP \citep{jones1999first}, Model Coupling Toolkit (MCT)
\citep{larson2005model}, OASIS \citep{craig2017development}, Tempest
Remap \citep{ullrich2015arbitrary,ullrich2016arbitrary}, etc.), fluid-structure
interaction or heat transfer (such as MpCCI \citep{joppich2006mpcci}
and preCICE \citep{bungartz2016precice}), and general-purpose remap
software (such as PANG \citep{gander2013algorithm} and Data Transfer
Kit (DTK) \citep{slattery2013data,slattery2016mesh}). We first review
the methodology behind some of these packages, with a focus on node-to-node
transfer, and then review some closely related methods for high-order
reconstruction.

\subsubsection{Low-order nonconservative and conservative remap}

In 1990s and 2000s, data-transfer methods primarily on only first-
or second-order accuracy, because the prevailing numerical discretization
methods in engineering (including finite volume methods and linear
finite element methods) only had low-order accuracy. One of the most
primitive data-transfer methods is piecewise interpolation. This approach
is particularly convenient if the source mesh has an associated finite-element
function space, in that one could use the same function space for
interpolation. In this case, we refer to the approach as \emph{consistent
interpolation}. If the basis functions are linear (with simplicial
elements) or bilinear (with quadrilateral elements), the interpolation
is second order. In a finite-volume method and scatter-data interpolation,
the nearest-point interpolation is sometimes used, which is only first-order
accurate. If the basis functions are not known, some interpolation
(and more precisely, \emph{quasi-interpolation} \citep{de1990quasiinterpolants})
methods can use some alternative basis functions, such as thin-plate
splines \citep{harder1972interpolation} or other RBF \citep{beckert2001multivariate,bungartz2016precice,flyer2007transport,rendall2008unified}.

The major disadvantage of low-order interpolation is its high diffusiveness
in repeated transfer; see e.g. \citep{de2008comparison}. To overcome
this, one approach is to enforce \emph{conservation}, so that the
integrals over the source and target meshes are equal. Examples of
lower-order conservative methods include the first- and second-order
conservative remap (such as those in \citep{chesshire1994scheme,grandy1999conservative,gander2013algorithm,jones1999first}),
common-refinement-based $L^{2}$ projection \citep{jiao2004common},
etc. A commonality of these methods is that they all require computing
some integrals of the function defined on the source mesh over some
control volumes of each target node (or cell), and the numerical integration
is computed over the intersections of the elements (or cells) of the
source and target meshes. The collection of these intersections forms
the common refinement \citep{jiao2004common} or supermesh \citep{farrell2009conservative,farrell2011conservative},
whose computations require some sophisticated computational-geometry
algorithm for efficiency and for robustness in the presence of truncation
and rounding errors \citep{jiao2004overlayingI,jiao2004overlayingII}.
We also note that some low-order constrained interpolation and conservative
remap methods have been developed in the context of Arbitrary Lagrangian-Eulerian
(ALE) methods, such as \citep{bochev2005constrained,margolin2003second}.

\subsubsection{High-order remap methods}

Due to the increasing use of high-order discretization methods, more
recent development of data-transfer methods had focused on third and
higher order methods. If the source mesh uses a high-order finite
element or spectral element method, one can apply consistent interpolation
to achieve high-order accuracy. Such an approach may be implemented
with a discretization library, such as MOAB \citep{tautges2004moab,mahadevan2015sigma}.
More generally, a remap method can be independent of the discretization
methods. In this context, the method may be nonconservative or conservative.
Examples of the former include modified moving least square (MMLS)
\citep{joldes2015modified,slattery2016mesh}, etc. An example of the
latter is the conservative remap in Tempest Remap \citep{ullrich2015arbitrary,ullrich2016arbitrary}
and conservative interpolation \citep{adam2016higher}. Note that
high-order nonconservative remap methods tend to be significantly
less diffusive in repeated transfer than low-order nonconservative
methods, even though conservation is not enforced explicitly; see
e.g. \citep{slattery2016mesh}. However, because high-order methods
are more prone to oscillations, special attention is needed to preserve
positivity and convexity \citep{ullrich2015arbitrary,ullrich2016arbitrary}.

\subsubsection{High-order reconstructions}

The remap methods, especially the nonconservative ones, are closely
related to the \emph{reconstruction} problem, i.e., to reconstruct
a piecewise smooth function or its approximations at some discrete
points on a mesh given some known qualities at discrete points or
cells on the same mesh. The primary difference between remap and reconstruction
is that the former involves two meshes while the latter typically
involves only a single mesh. Hence, remap is in general more complicated
from a combinatorial point of view. From a numerical point of view,
conservation in remap is also more complicated than in reconstruction.
However, nonconservative high-order remap and reconstruction share
many similarities. In particular, the techniques used in high-order
remap, such as spline interpolation \citep{de1990quasiinterpolants},
moving least square (MLS) \citep{lancaster1981surfaces}, and variants
of MLS (known as MMLS) with regularization \citep{joldes2015modified,slattery2016mesh},
etc., originated from high-order reconstruction.

The remap method in this work utilizes two high-order reconstruction
techniques, known as \emph{CMF} (\emph{Continuous Moving Frames})
and \emph{WALF} (\emph{Weighted Average of Least-squares Fittings})
\citep{Jiao2011RHO}, both of which are based on \emph{weighted least
squares} (\emph{WLS}). These techniques were first proposed for reconstructing
surfaces and later adapted to reconstruct functions on surfaces \citep{RayWanJia12}.
CMF shares some similarities with some variant of MMLS (such as that
in \citep{slattery2016mesh}), in that CMF constructs a least-squares
fitting at each reconstruction point and it achieves accuracy and
stability through some local adaptivity (instead of the global construction
of the original MLS for smoothness \citep{Lancaster1981}). However,
CMF differs from MLS in terms of the choice of the coordinate systems,
stencils, and weighting schemes. WALF constructs a least-squares fitting
at each node of a mesh using CMF and then blend the fittings taking
a weighted average using piecewise linear or bilinear basis functions
within each element (such as a triangle or quadrilateral). In general,
CMF is more accurate than WALF. For meshing applications, WALF tends
to be more efficient when many points need to be evaluated within
a single element. For remap in multiphysics applications, however,
a transfer operator can be built in a preprocessing step and be reused
in repeated transfers. Hence, In this work, we use CMF as the basis
for remapping smooth regions. From the mathematical point of view,
the major challenge is to overcome the potential oscillations (a.k.a.
the Gibbs phenomenon) associated with CMF, for which we will leverage
WALF to detect singularities. 

\subsection{Gibbs phenomena at $C^{0}$ and $C^{1}$ discontinuities}

It is well known that discontinuities require special attentions in
numerical approximations. The most studied type is the $C^{0}$ discontinuities,
also known as \emph{edges} \citep{tadmor2007filters}, (\emph{simple})\emph{
jump discontinuities }\citep{cates2007detecting,hewitt1979gibbs},
\emph{faults} (or \emph{vertical faults}) \citep{bozzini2013detection,de2002vertical},
etc. These discontinuities are the most prominent because they tend
to lead to $\mathcal{O}(1)$ oscillations (or ``ringing''), which
do not vanish under mesh refinement. These oscillations were first
analyzed by Wilbraham in 1848 \citep{wilbraham1848certain}, about
50 years before the famous analysis of the overshoots in Fourier series
for saw-tooth functions by Gibbs in 1899 \citep{gibbs1899} in response
to an observation by Michelson in 1898 \citep{michelson1898letter}.
For this reason, such oscillations are properly referred to as the
\emph{Gibbs-Wilbraham }phenomenon \citep{hewitt1979gibbs}, or more
commonly as the \emph{Gibbs phenomenon}. Besides occurring in Fourier
series, the $\mathcal{O}(1)$ oscillations also occur in interpolation
or approximation using cubic splines \citep{richards1991gibbs,zhang1997convergence},
orthogonal polynomials \citep{gottlieb1997gibbs}, RBF \citep{fornberg2007gibbs,jung2007note},
least squares \citep{foster1991gibbs}, wavelets \citep{kelly1996gibbs,shim1996gibbs},
$L^{2}$ projection \citep{jiao2004common}, etc. Some tight constant
factors can be established in 1D for some of these problems. However,
for least-squares-based approximations in 2D and higher dimensions,
it is impractical to establish such a precise bound of the oscillations.
Instead, we focus on the asymptotic effect of discontinuities, which
is easy to understand from the Taylor series expansion and is sufficient
for developing robust techniques to detect and resolve them in this
work.

Consider a piecewise smooth function $f(\vec{x}):\mathbb{R}^{d}\rightarrow\mathbb{R}$.
The $d$-dimensional Taylor series of $f$ about a point $\vec{x}_{0}$
is given by
\begin{align}
f(\vec{x}_{0}+\vec{h}) & =\sum_{k=0}^{p}\frac{1}{k!}\partial_{\hat{\vec{h}}}^{k}f(\vec{x}_{0})\left\Vert \vec{h}\right\Vert _{2}^{k}+\frac{\partial_{\hat{\vec{h}}}^{p+1}f(\vec{x}_{0}+\vec{\epsilon})}{(p+1)!}\Vert\vec{h}\Vert_{2}^{p+1},\label{eq:Taylor-polynomial-residual}
\end{align}
where $\partial_{\hat{\vec{h}}}^{k}$ denotes the $k$th directional
derivative along the direction $\hat{\vec{h}}=\vec{h}/\left\Vert \vec{h}\right\Vert _{2}$,
and $\Vert\vec{\epsilon}\Vert\leq\Vert\vec{h}\Vert$; see, e.g., \citep{humpherys2017foundations}
for a proof. The summation term in (\ref{eq:Taylor-polynomial-residual})
defines the degree-$p$ \emph{Taylor polynomial} $T_{p}(\vec{x})$
at $\vec{x}_{0}$. The last term in (\ref{eq:Taylor-polynomial-residual})
is the \emph{remainder} $r(\vec{x})$, which corresponds to the residual
of the approximation $f(\vec{x})$ by $T_{p}(\vec{x})$. If $f$ is
continuously differentiable to at least $p$th order, the residual
$r(\vec{x})$ is $\mathcal{O}(\Vert\vec{h}\Vert_{2}^{p+1})$, i.e.,
\[
\left|r(\vec{x})\right|=\left|f(\vec{x})-T_{p}(\vec{x})\right|=\mathcal{O}(\Vert\vec{h}\Vert_{2}^{p+1}).
\]
However, if $f$ has a $C^{q}$ discontinuity in the direction $\hat{\vec{h}}$
at $\vec{x}_{0}+\vec{\epsilon}$, i.e., $\partial_{\hat{\vec{h}}}^{q}f$
has a jump at $\vec{x}_{0}+\vec{\epsilon}$, then $\partial_{\hat{\vec{h}}}^{p+1}f$
is unbounded. As a result, the approximation of $f(\vec{x})$ by $T_{p}(\vec{x})$
in general are expected to lead to overshoots or undershoots (or oscillations)
of $\mathcal{O}(\Vert\vec{h}\Vert_{2}^{p})$ or larger. Note that
even if $f$ is approximated by non-polynomials, this argument based
on Taylor series also applies with a simple application of triangle
inequality. In particular, let $g(\vec{x})$ be a smooth approximation
of $f$ and let $\tilde{T}_{p}$ denote the degree-$p$ Taylor polynomial
of $g(\vec{x})$ at $\vec{x}_{0}$. Then,
\[
\left|f(\vec{x})-g(\vec{x})\right|\leq\left|r(\vec{x})\right|+\left|g(\vec{x})-\tilde{T}_{p}(\vec{x})\right|,
\]
where the second term is $\mathcal{O}(\Vert\vec{h}\Vert_{2}^{p+1})$
under the smoothness assumption of $g$, but the first term is $\mathcal{O}(\Vert\vec{h}\Vert_{2}^{p})$
or larger and it dominates the overall error. Hence, for $C^{0}$
discontinuities, we can expect to observe $\mathcal{O}(1)$ oscillations,
even when using linear approximations, as in $L^{2}$ projection using
linear polynomials.

From the Taylor series expansions, it is easy to see that overshoots
can occur not only at $C^{0}$ discontinuities but also at $C^{1}$
and even $C^{2}$ discontinuities. In general, the oscillations due
to $C^{1}$ discontinuities are expected to be $\mathcal{O}(h)$.
Theoretically, these oscillations should vanish as $h$ tends to $0$,
so they are often ignored. However, its treatment has drawn some attention
recently, and it is often referred to as \emph{derivative jump discontinuities
}\citep{cates2007detecting} or \emph{oblique }(or \emph{gradient})\emph{
faults} \citep{bozzini2013detection,de2002vertical}. In the context
of data transfer in multiphysics coupling, it is important to resolve
the $\mathcal{O}(h)$ oscillations at $C^{1}$ discontinuities, because
they may accumulate over time and lead to $\mathcal{O}(1)$ oscillations
and even instabilities. Although these $\mathcal{O}(h)$ oscillations
may be controlled by numerical diffusion in some numerical schemes,
such as the flux limiters and TVD time integration schemes in FVM
for hyperbolic conservation laws, for generality, it is desirable
for the data-transfer method to be stable without relying on numerical
diffusion in physics models. In addition, controlling $\mathcal{O}(h)$
oscillations also helps safeguarding the accumulation of the $\mathcal{O}(h^{2})$
oscillations near $C^{2}$ discontinuities.

\subsection{Resolution of Gibbs oscillations}

\subsubsection{Resolution of oscillations in 1D}

Because the mathematical analysis of the Gibbs phenomenon in the literature
primarily focused on 1D reconstructions, the overwhelming majority
of the resolution techniques were also for 1D problems. We briefly
mention two classes of classical techniques in 1D. The first class
is \emph{filtering} and \emph{mollification}, which recovers the piecewise
solutions in Fourier spaces and in physical spaces, respectively.
Examples of the former include Fejer averaging, Lanczos-like filtering,
Vandeven and Daubechies filters, etc., which are specific to Fourier
expansions; see e.g. \citep{zygmund1959trigonometric} and \citep[Chapter 2]{jerri2013gibbs}
for some surveys. Examples of the latter include spectral mollifiers
based on Gegenbauer polynomials in \citep{gottlieb1997gibbs,jung2004generalization,shizgal2003towards}
and adaptive mollifiers in \citep{tadmor2007filters}.  The fundamental
idea of adaptive mollifiers is to detect the edges \citep{gelb1999detection,gelb2000detection}
and then apply spectral (or other) mollifiers locally. Similar ideas
have been applied for approximating piecewise smooth functions in
1D using splines and moving least squares in \citep{lipman2009approximating}.

The second class is the \emph{Godunov} and \emph{ENO-type methods}.
Examples include the original Godunov's method \citep{godunov1959difference},
Van Leer's MUSCL scheme \citep{van1979towards}, Colella and Woodward's
piecewise parabolic method (PPM) \citep{colella1984piecewise}, Harten's
subcell-resolution method \citep{harten1989eno}, and some WENO reconstructions
for FVM \citep{liu1994weighted,jiang1996efficient,shu1999high}. These
methods were initially developed for the finite volume methods (FVM)
for hyperbolic problems, and they typically involve reconstructing
piecewise smooth functions from cell averages. Among these techniques,
the ENO and WENO techniques are the most general, and their key ideas
lie at the approximation level, where a nonlinear adaptive procedure
is used to automatically choose the locally smoothest stencil, hence
avoiding crossing discontinuities in the interpolation process as
much as possible. This basic idea can be generalized to reconstructing
nodal values, such as in the so-called \emph{WENO interpolation} in
finite difference methods \citep{shu2009high}. Similar to mollifiers,
the WENO schemes may be applied globally without explicitly identifying
the discontinuities, or be applied adaptively by coupling with some
edge-detection techniques \citep{qiu2005comparison}.

\subsubsection{Resolution of oscillations in higher dimensions}

The two broad classes of 1D techniques may be applied to structured
meshes in a dimension-by-dimension fashion. For example, in climate
modeling, Lauritzen and Nair \citep{lauritzen2008monotone} developed
monotonicity-preserving filters on latitude--longitude and cubed-sphere
grids by adapting those developed by Zerroukat et al. using monotonic
parabola \citep{zerroukat2005monotonic} or parabolic spline method
(PSM) \citep{zerroukat2006parabolic}. However, when generalizing
the techniques to unstructured meshes or scattered data, most of the
techniques in 1D or for structured meshes encounter significant technical
difficulties. However, some of the basic ideas can still be preserved
on higher dimensions. Hence, we also broadly categorize the methods
into two classes.

The \emph{filtering} techniques in 2D and higher dimensions typically
reply on post-processing to resolve oscillations. This requires detecting
the discontinuities and then adapt the kernels (i.e., basis functions)
accordingly. Various techniques have been developed in the literature
to detect discontinuities. The detection of $C^{0}$ discontinuities
is often referred to as \emph{edge detection} \citep{archibald2008determining,canny1987computational,gelb2006adaptive,petersen2012hypothesis,romani2019edge,viswanathan2008detection,ziou1998edge}
or \emph{troubled-cell detection} \citep{qiu2005comparison,yang2009parameter},
and it has many applications in signal and image processing, data
compression, shock capturing in CFD, etc. These techniques typically
use some indicators for the jump and then compare the indicator against
some thresholds. Similarly, detecting $C^{1}$ discontinuities requires
an indicator of the derivative jump followed by some thresholding.
The indicators are often motivated by finite differences (such as
in \citep{bozzini2013detection,cates2007detecting,de2002vertical})
or Taylor series expansions (such as in the \emph{polynomial annihilation
edge detection} in \citep{archibald2005polynomial,archibald2008determining,saxena2009high}).
In \citep{jung2011iterative}, Jung et al. detected discontinuities
by taking advantage of the instabilities of the expansion coefficients
of multiquadric RBF at local jumps. Romani et al. \citep{romani2019edge}
developed some variants of the approach using variably scaled kernels
(VSK) \citep{bozzini2014interpolation} and Wendland's $C^{2}$ RBF
\citep{wendland1995piecewise}. In terms of thresholding, they are
often based on some statistical analysis (i.e., outlier detection)
\citep{fleishman2005robust} or based on some double (or hysteresis)
thresholding \citep{canny1987computational}. After the discontinuities
have been detected, various reconstruction strategies have been used
in practice. For reconstructions over structured grids, such as in
image processing, it is relatively straightforward to adapt some limiters
used in hyperbolic conservation laws to resolve discontinuities in
a dimension-by-dimension fashion, such as in \citep{archibald2008determining,gelb2006adaptive}.
However, such an approach does not generalize to unstructured meshes
easily. For unstructured meshes or scattered data points, Rossini
\citep{rossini2018interpolating} proposed to use VSK to resolve discontinuities,
but it required the exact location of discontinuities in order to
shift the radial functions.

In terms of Godunov and ENO-type methods, the generalization of the
low-order methods tends to be straightforward within FVM. However,
for the high-order ENO and WENO schemes, the generalization to unstructured
meshes is challenging. In \citep{hu1999weighted}, Hu and Shu proposed
a WENO scheme on triangular meshes, and it was further extended to
tetrahedral meshes by Zhang and Shu in \citep{zhang2009third}. However,
some point distributions may lead to negative weights and in turn
lead to potential instability. Shi, Hu and Shu \citep{shi2002technique}
proposed a technique to mitigate the issue, but the method had limited
success due to large condition numbers of the local linear systems.
In \citep{xu2011point}, Xu et al. proposed a hierarchical reconstruction
technique for discontinuous Galerkin and finite volume methods on
triangular meshes, and Luo et al. \citep{luo2013reconstructed} extended
the ideas to tetrahedral meshes and developed the so-called hierarchical
WENO (or HWENO). In \citep{liuzhang2013robust}, Liu and Zhang proposed
a hybrid of two different reconstruction strategies to improve stability
of WENO reconstruction. In \citep{liu2016wls}, Liu and Jiao proposed
WLS-ENO, which adapts the weights in a weighted least squares formulation
instead of taking a convex combination of interpolations. WLS-ENO
overcomes the stability issues associated negative weights, and it
was shown to be robust in FVM and is easy to implement.

The resolution technique proposed in this work may be considered as
a hybrid of adaptive mollification and ENO-type method, in that we
detect $C^{0}$ and $C^{1}$ discontinuities in high-order CMF and
then apply quadratic WLS-ENO as mollifier (or filters in the physical
space) near discontinuities. Our detection technique is different
from others in that it introduces two novel discontinuity indicators
and a novel dual thresholding strategy to detect discontinuities.
The indicators use the observation that WALF tends to be more prone
to oscillations than CMF near discontinuities, so that we can detect
and resolve discontinuities before oscillations manifest themselves
in CMF. In addition, we will introduce a new weighting scheme for
WLS-ENO based on one of the new discontinuity indicators.

\section{\label{sec:Superconvergent-weighted-least}Superconvergent weighted
least squares (WLS)}

In this section, we describe WLS remap for smooth functions on surface.

\subsection{WLS with Cartesian coordinates}

We first derive WLS in Cartesian coordinates using Taylor series.
Consider a function $f(\vec{u}):\mathbb{R}^{2}\rightarrow\mathbb{R}$
at a given point $\boldsymbol{u}_{0}=\left[0,0\right]^{T}$, and assume
its value is known at a sample of $m$ points $\vec{u}_{i}$ near
$\vec{u}_{0}$, where $1\leq i\leq m$. We refer to these points as
the \emph{stencil} for WLS. Suppose $f$ is continuously differentiable
to $(p+1)$st order for some $p>1$. We approximate $f(\boldsymbol{u})$
to (\emph{$p+1$})st\emph{ }order accuracy about $\boldsymbol{u}_{0}$
as 
\begin{equation}
f(\boldsymbol{u})=\sum_{q=0}^{p}\sum_{j,k\ge0}^{j+k=q}c_{jk}u^{j}v^{k}+\mathcal{O}\left(\Vert\boldsymbol{u}\Vert^{p+1}\right),\label{eq:FTaylorseries_disc}
\end{equation}
which is an alternative form of that in (\ref{eq:Taylor-polynomial-residual}),
where $c_{jk}=\dfrac{1}{j!k!}\dfrac{\partial^{j+k}}{\partial u^{j}\partial v^{k}}f(\boldsymbol{0})$.
Supper there are $n$ coefficients, i.e., $n=(p+1)(p+2)/2$, and assume
$m\ge n$. Let $f_{i}$ denote $f(\vec{u}_{i})$. We then obtain a
system of $m$ equations 
\begin{equation}
\sum_{q=0}^{p}\sum_{j,k\ge0}^{j+k=q}c_{jk}u_{i}^{j}v_{i}^{k}\approx f_{i}\label{eq:Vandermone-system}
\end{equation}
for $1\leq i\leq m$. The equation can be written in matrix form as
$\vec{A}\vec{x}\approx\vec{b}$, where $\vec{A}\in\mathbb{R}^{m\times n}$
is an \emph{generalized Vandermonde matrix}, $\vec{x}\in\mathbb{R}^{n}$
is composed of $c_{jk}$, and $\vec{b}\in\mathbb{R}^{m}$ is composed
of $f_{i}$.

The generalized Vandermonde system from (\ref{eq:Vandermone-system})
is rectangular, and we solve it by minimizing a weighted norm of the
residual vector $\vec{r}=\vec{b}-\vec{A}\vec{x}$, i.e.,
\begin{equation}
\min_{\vec{x}}\Vert\vec{r}\Vert_{\vec{W}}\equiv\min_{\vec{x}}\Vert\vec{W}(\vec{A}\vec{x}-\vec{b})\Vert_{2},\label{eq:weighted_norm}
\end{equation}
where $\vec{W}=\mbox{diag}\{\omega_{1},\omega_{2},\dots,\omega_{m}\}$
is diagonal. The weighting matrix $\vec{W}$ plays a fundamental role
in $\vec{W}$, in that each diagonal entry $w_{i}$ in $\vec{W}$
assigns a weight to each row in the generalized Vandermonde system.
If both $\vec{A}$ and $\vec{W}$ are both nonsingular, then $\vec{W}$
has no effect on the solution. However, if $m\neq n$ or $\vec{A}$
is singular, then a different $\vec{W}$ leads to a different solution;
we will address the choice of the weights shortly. Given $\vec{W}$,
we further scale the columns of $\vec{W}\vec{A}$ with a diagonal
matrix $\vec{T}$, so that the condition number of the rescaled Vandermonde
matrix $\tilde{\vec{A}}=\vec{W}\vec{A}\vec{T}$ is improved. This
process is known as \emph{column equilibration} \citep[p. 265]{Golub13MC}.
Then, WLS reduces to the least squares problem 
\begin{equation}
\vec{W}\vec{A}\vec{T}\vec{y}\approx\vec{W}\vec{b}.\label{eq:least-squares}
\end{equation}
The coefficients in $c_{ij}$ are then given by $\vec{x}=\vec{T}\vec{y}$.

We give two important algorithmic details. First, we solve (\ref{eq:least-squares})
using truncated QR factorization with column pivoting (QRCP) \citep{Golub13MC}.
The QRCP is
\[
\tilde{\vec{A}}\vec{P}=\vec{Q}\vec{R},
\]
where $\vec{Q}$ is $m\times n$ with orthonormal column vectors,
$\vec{R}$ is an $n\times n$ upper-triangular matrix, $\vec{P}$
is a permutation matrix, and the diagonal entries in $\vec{R}$ are
in descending order. Second, we control the condition number of $\tilde{\vec{A}}$
to avoid numerical instabilities. To this end, we estimate the 1-norm
condition number of $\vec{R}$ \citep{Higham87SCN}. If the condition
number is large, then we first try to enlarge the stencils and then
incrementally drop the right-most highest-degree terms in $\tilde{\vec{A}}\vec{P}$
if the stencils cannot be enlarged further. These strategies add more
rows to, and remove columns from, the Vandermonde systems, respectively,
both of which reduce the condition number of $\tilde{\vec{A}}$ and
in turn ensure the stability of WLS.

\subsection{WLS on discrete surfaces}

When applying WLS on surfaces, we use the same construction as in
\emph{Continuous Moving Frames} (\emph{CMF}) for high-order surface
reconstruction \citep{Jiao2011RHO}. We adapt CMF to remap between
discrete surfaces composed of triangles and quadrilaterals. In the
following, we will use WLS and CMF interchangeable, unless otherwise
noted.

Given a node $\vec{x}_{0}$ on a target mesh, we first construct a
local $uv$ coordinate frame as follows. Let $\vec{m}_{0}$ denote
an approximate normal at $\vec{x}_{0}$, which can be first-order
estimation by averaging face normals or be computed from the analytical
surface (such as a sphere in earth modeling). Let $\vec{t}_{1}$ and
$\vec{t}_{2}$ be an orthonormal basis of approximate tangent plan
orthogonal to $\vec{m}_{0}$, which we obtain using Gram-Schmidt orthogonalization.
The local $uv$ coordinate frame is then centered at $\vec{x}_{0}$
with axes $\vec{t}_{1}$ and $\vec{t}_{2}$.

The construction of stencils requires some special attention in remap.
Unlike reconstruction, the stencil for a target node is composed of
nodes on the source mesh. To build the stencil for a target node $\vec{x}_{0}$,
we first locate the triangle or quadrilateral $\tau$ on the source
mesh that contains $\vec{x}_{0}$, and then use the union of the $k$-ring
neighborhood of each node of $\tau$ on the source mesh. The specific
choice of $k$ depends on the degree of the polynomial and the specific
weighting scheme, as we will describe in Section~\ref{subsec:Optimizing-weights}.

We note some key difference between our WLS remap and MMLS in DTK
when applied to surfaces. First, our approach uses 2D Taylor series
within local $uv$ coordinates in a tangent space, whereas DTK uses
the global $xyz$ coordinate system. Note that for degree-$p$ polynomials,
the number of coefficients grows quadratically in 2D but cubically
in 3D. Hence, it is more efficient to use higher-degree (such as quartic)
fittings with CMF to achieve higher accuracy than DTK. Second, DTK
is purely mesh-less and it finds stencils by using $k$ nearest neighbors
(a.k.a. KNN), which requires logarithmic-time complexity on average
to determine the stencils. In contrast, our approach uses a mesh-based
algorithm to compute the stencils by taking advantage of locality
\citep{jiao1999mesh}, so the amortized cost to determine the stencil
of a node is constant.

\subsection{\label{subsec:Optimizing-weights}Optimizing weights for remap}

As shown in \citep{JZ08CCF} and \citep{RayWanJia12}, if the condition
number of the rescaled Vandermonde system $\tilde{\vec{A}}$ is bounded,
then the function values of a smooth function $f:\Gamma\rightarrow\mathbb{R}$
can be reconstructed to $\mathcal{O}(h^{p+1})$ on a discrete surface
$\Gamma$, where $h$ is proportional to the ``radius'' of the stencil.
Furthermore, if the stencils are (nearly) symmetric about the origin
of the local coordinate system, the reconstruction can superconverge
at $\mathcal{O}(h^{p+2})$ for even-degree $p$ due to error cancellation
\citep{JZ08CCF,Jiao2011RHO}, analogous to the error cancellation
in centered differences. When using CMF for remap, however, the near-symmetry
assumption is in general violated, because a target node may lie anywhere
within an element of a source mesh. Hence, it is impractical to achieve
$\mathcal{O}(h^{p+2})$ superconvergence. However, with proper choice
of stencils and weighting schemes, we can achieve $\mathcal{O}(h^{p+1+s})$
convergence for some $0<s\leq1$ for even $p$. To maximize this rate,
we observe two key criteria for the weights: \emph{asymptotic rates}
and \emph{smoothness} of the weights.

In \citep{JZ08CCF}, Jiao and Zha proposed to use weights that are
approximately inversely proportional to the square root of the residual
term associated with each node, i.e.,
\begin{equation}
\omega_{j}^{\text{ID}}=\gamma_{j}^{+}\left(\sqrt{r_{j}^{2}+\epsilon_{\text{ID}}}\right)^{-p/2},\label{eq:inverse-dist-weight}
\end{equation}
where $r_{j}=\left\Vert \vec{u}_{j}\right\Vert $, $p$ is the degree
of the WLS, $\epsilon$ is a safeguard against division by zero, and
\[
\gamma_{i}^{+}=\max\{0,\vec{m}_{j}^{T}\vec{m}_{0}\}
\]
is a safeguard against ``short circuiting'' for coarse meshes. We
refer to (\ref{eq:inverse-dist-weight}) as the \emph{inverse-distance-based}
weights. These weights perform well with nearly symmetric stencils.
However, they are not ideal for remap, because it has a large gradient
near $\vec{x}_{0}$, especially if $\epsilon_{\text{ID}}$ is small.
Therefore, if the nodes closest to $\vec{x}_{0}$ are highly asymmetric,
which tend to be the case in remap, then there is no error cancellation
for these nodes. To overcome this issue, it is desirable for the weights
to be smoother about the origin, or more precisely, the first, second,
and even higher-order derivatives of the weight should $0$ at the
origin. This analysis led us to the choice of a new weighting scheme,
or \emph{scaled Buhmann's function},
\begin{equation}
\omega_{j}^{\text{SB}}=\gamma_{j}^{+}\phi(r_{j}/\rho),\label{eq:Buhmann_weights}
\end{equation}
where $\gamma_{j}^{+}$ and $r_{j}$ are the same as those in (\ref{eq:inverse-dist-weight}),
$\rho$ is the \emph{cut-off radius}, and

\begin{equation}
\phi(r)=\begin{cases}
\frac{112}{45}r^{9/2}+\frac{16}{3}r^{7/2}-7r^{4}-\frac{14}{15}r^{2}+1/9 & \text{if }0\leq r\leq1\\
0 & \text{if }r>1
\end{cases}.\label{eq:Buhmann-function}
\end{equation}
The function $\phi(r)$ is a $C^{3}$ function due to Buhmann \citep{buhmann2001new},
who introduced it as a compactly supported RBF by generalizing the
polynomial radial functions due to Wu \citep{wu1995compactly} and
Wendland \citep{wendland1995piecewise}. Besides compactness and smoothness,
these radial functions are ``positive definite,'' which is a desirable
property for RBF \citep{buhmann2000radial,buhmann2003radial} but
is irrelevant in our context. In Section~\ref{subsec:Comparison-of-weighting-schemes},
we will compare our proposed scaled Buhmann weights versus other weighting
functions.

\begin{figure}
\subfloat[Weights based on inverse-distance.]{\begin{centering}
\includegraphics[width=0.48\columnwidth]{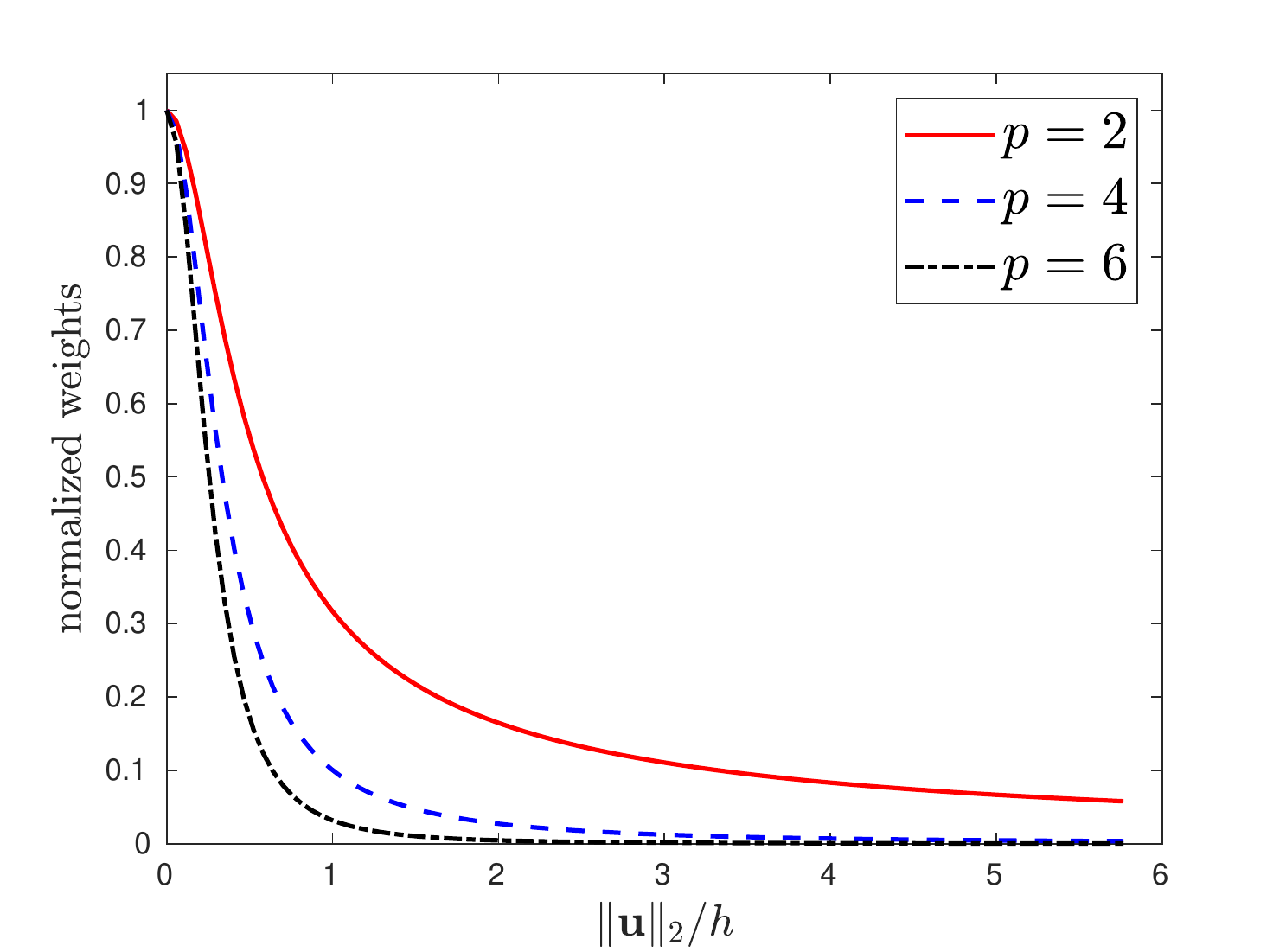}
\par\end{centering}
}\hfill\subfloat[Weights based on Buhmann's function]{\begin{centering}
\includegraphics[width=0.48\columnwidth]{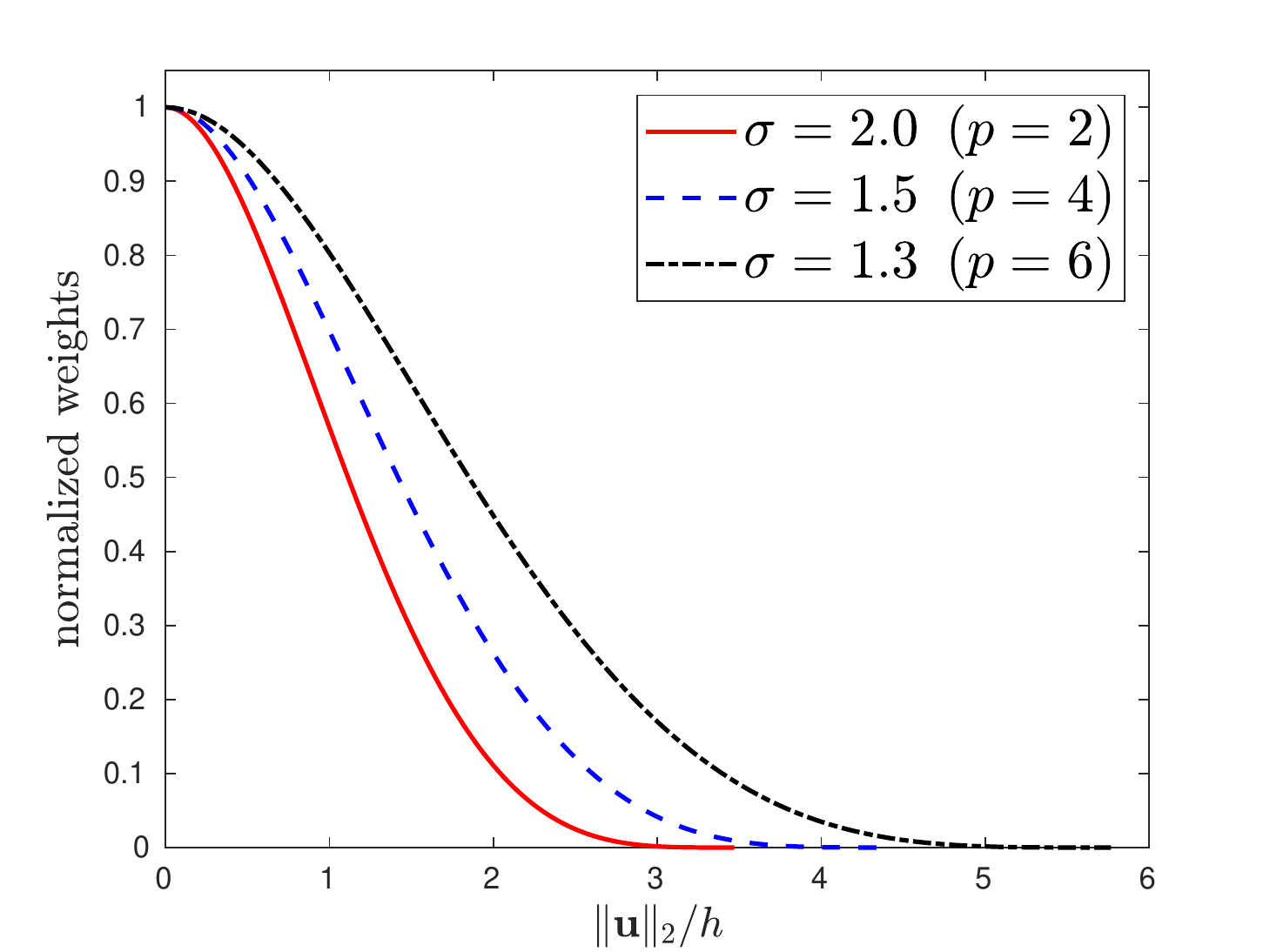}
\par\end{centering}
}

\caption{\label{fig:Example-shapes}Example shapes of inverse-distance-based
and scaled Bushman weights.}
\end{figure}

In (\ref{eq:Buhmann_weights}), the key question is how to choose
the cut-off radius $\rho$. For robustness, this radius must be large
enough to allow sufficient number of nonzero weights in the generalized
Vandermonde system. As in \citep{JZ08CCF}, we require the number
of points in the stencil to be at least 1.5 times of that of the number
of coefficients in the Taylor polynomial. Let $R=r_{k}$, where $k=\left\lceil 0.75(p+1)(p+2)\right\rceil $.
We define $\rho=\sigma R$ for some $\sigma>1$ and then choose $\sigma$
to minimize the remap errors via numerical optimization, as we will
detail in Section~\ref{subsec:Comparison-of-weighting-schemes}.
Through numerical optimization, we obtain $\sigma=2.0$, 1.6 and $1.4$
for degree-2, 4, and 6 WLS, respectively. Figure~\ref{fig:Example-shapes}
shows the example shapes of inverse-distance-based and scaled Buhmann
weights on a uniform 2D mesh for degree-2, 4, and 6 WLS. For inverse-distance
weights, we used $\epsilon_{\text{ID}}=0.01\bar{h}^{2}$, where $\bar{h}$
denotes the average edge length. It can be seen that due to its $C^{3}$
continuity, Buhmann's function is much smoother at $r=0$. In addition,
the different cut-off radii led to different asymptotic behaviors
of the scale Buhmann weights at the tails, similar to those of the
inverse-distance-based weights.

In terms of implementation, to select the candidate nodes in the stencils
for a target node $\vec{x}_{0}$, we start with the union of $\lfloor1.5p\rfloor/2$-rings
of the nodes of the source element that contains a given target node.
If there are less than $0.75(p+1)(p+2)$ points in the stencil, we
enlarge the ring size in 0.5 increments as defined in \citep{Jiao2011RHO}.
This adaptive computation of rings requires a proper mesh data structure
to achieve constant-time per node, for which we use an array-based
half-edge data structure \citep{DREJTAHF2014}. In addition, for efficiency,
we construct the transfer operator in the form of a sparse matrix,
of which each row corresponds to a node in the target mesh. After
obtaining this operators, the remap involves only a matrix-vector
multiplication.

\section{\label{sec:Detection-and-resolution}Detection and resolution of
$C^{0}$ and $C^{1}$ discontinuities}

Like other high-order methods, WLS with scaled Buhmann weights may
suffer from Gibbs phenomena near discontinuities. In this section,
we extend the WLS-ENO approach to remap. For accuracy and efficiency,
our approach first detects the discontinuous regions and then applies
WLS-ENO to the detected discontinuities. To this end, we introduce
a robust detector of discontinuities and a new weighting scheme for
WLS-ENO.

\subsection{Detection of discontinuities}

We first describe a robust technique to detect discontinuities in
the context of remap. Our approach is composed of four steps. First,
we compute element-based discontinuity indicators on the source mesh.
Second, we convert these indicators to node-based indicators on the
source mesh. Third, we obtain node-based discontinuity markers on
the source mesh using a novel dual-thresholding strategy. Finally,
we transfer the nodal markers to the target mesh. In the following,
we describe the four steps separately, including their derivations.

\subsubsection{\label{subsec:Element-based-indicators}Element-based indicators
on source mesh}

We first compute a \emph{discontinuity indicator }at each element
on the source mesh. While most other indicators are based on computing
the differences of neighboring values, our approach is different in
that it computes the difference between a second-order interpolation
and a higher-order reconstruction at the center of the element (i.e.,
triangle or quadrilateral). For the former, we use linear (or bilinear)
interpolation within the element, which is non-oscillatory because
it is a convex combination of the nodal values. For the latter, we
use quadratic WALF reconstruction, which is known to be more prone
to oscillations than CMF near discontinuities \citep{Jiao2011RHO}.

More precisely, given an element $e$, let $g_{e,1}$ denote the linear/bilinear
interpolation, i.e., the average of the nodal values. Let $g_{e,2}$
denote the quadratic WALF reconstruction at the element center, i.e.,
the average of the reconstructed values at the element center from
the quadratic CMF fittings using 1.5-ring stencil at the nodes of
the element. Then, the indicator value for element $e$ is given by
\begin{equation}
\alpha_{e}=g_{e,2}-g_{e,1}.\label{eq:element-based-indicator}
\end{equation}
We store the operator for computing $\alpha_{e}$ for all the elements
as a sparse matrix, of which each row contains the coefficients associated
with the nodes in its $k$-ring neighborhood. As for the transfer-operator
for smooth regions, we compute this sparse matrix \emph{a priori},
so that it can be reused in repeated transfer for efficiency.

The preceding definition of $\alpha_{e}$ is fairly straightforward,
but its numerical values have two important properties. First, its
sign is indicative, in that a positive and negative sign indicates
a local overshoot and undershoot, respectively. Second, its magnitude
has different asymptotic behavior at smooth regions and at discontinuities.
In particular, let $f(\vec{c})$ denote the exact value at the element
center. For smooth functions $g_{i}-f(\vec{c})=\mathcal{O}(h^{2})$
for $i=1$ and $2$, so $\alpha_{e}$ is $\mathcal{O}(h^{2})$. Assume
the function is at least $C^{1}$ within an element $e$. Near $C^{0}$
discontinuities, $g_{2}-f(\vec{c})=\mathcal{O}(1)$, so $\alpha_{e}=\mathcal{O}(1)$.
Near $C^{1}$ discontinuities, $g_{2}-f(\vec{c})=\mathcal{O}(h)$,
so $\alpha_{e}=\mathcal{O}(h)$. These different asymptotic behaviors
lead to clear gaps in the magnitude of the $\alpha$ values at $C^{0}$
discontinuities, $C^{1}$ discontinuities, and at smooth regions,
which will be useful for classification.

\subsubsection{Node-based indicators on source mesh}

After obtaining the element-based $\alpha$ values, we then use them
to define node-based indicators on the same mesh. Given a node $v$,
we define the node-based indicator as

\begin{equation}
\beta_{v}=\frac{\sum_{v\in e}\vert\alpha_{e}-\bar{\alpha}\vert}{\left|\sum_{v\in e}\alpha_{e}\right|+\epsilon_{\beta}\delta f_{g}h_{g}^{2}+\text{realmin}},\label{eq:node-based-indicator}
\end{equation}
where $\bar{\alpha}$ denotes the average of the $\alpha$ values
in the 1-ring neighborhood of $v$, i.e., $\bar{\alpha}=\sum_{v\in e}\alpha_{e}/k$,
where $k$ is the number of elements incident on $v$. The second
term in the denominator is a safeguard against division by a too small
value, which may happen if the function is locally linear. In particular,
$\delta f_{g}$ denotes the global range of function over the mesh,
$h_{g}$ denotes a global measure of average edge length in the $xyz$
coordinate system, and $\epsilon_{\beta}=10^{-3}$. The last term
$\text{realmin}$ in the denominator denotes the smallest positive
floating-point number, which is approximately 2.2251e-308 for double-precision
floating-point numbers and further protects against division by zero
if the input function is a constant. Note that $\alpha_{e}$ and $\delta f_{g}h_{g}^{2}$
have the same units for smooth regions, so $\beta_{v}$ is non-dimensional. 

From a practical point of view, given the element-based $\alpha$
values, it is efficient to compute the node-based $\beta$ values.
From a numerical point of view, the definition of $\beta_{v}$ requires
some justification. Let us assume $f$ is sufficiently nonlinear,
so that $\beta_{v}\approx\sum_{v\in e}\vert\alpha_{e}-\bar{\alpha}\vert/\left|\sum_{v\in e}\alpha_{e}\right|$.
First, note that this quantity is no longer dependent on $h$. Second,
the enumerator of $\beta_{v}$ captures the variance of $\alpha_{e}$,
and the denominator further amplifies the variance when $\alpha_{e}$
varies in signs, which tend to occur near $C^{0}$ discontinuities.
These properties are very useful in designing the thresholding strategy,
which we describe next.

\subsubsection{Node-based markers on source mesh via dual thresholding}

After obtaining the element-based indicators $\alpha_{e}$ and node-based
indicators $\beta_{v}$, we then use them to in a dual-thresholding
strategy for detecting discontinuities. In particular, we mark a node
$v$ as discontinuity if it has a large $\beta$ value and one of
its incident elements has a large $\alpha$ value, i.e., 
\[
\text{disc}_{v}=\begin{cases}
1 & \text{if }\ensuremath{\beta_{v}>\kappa}\text{ and }\exists e\ni v\text{ s.t. }\alpha_{e}>\tau\\
0 & \text{otherwise}
\end{cases},
\]
where $\kappa$ and $\tau$ are two thresholds. To detect both $C^{0}$
and $C^{1}$ discontinuities, we set $\kappa=0.3$ and 
\begin{equation}
\tau=\max\left\{ \underbrace{C_{\ell}\delta f_{\ell}h_{\ell}^{0.5}}_{\tau_{\ell}},\underbrace{C_{g}\delta f_{g}h_{g}^{1.5}}_{\tau_{g}}\right\} ,\label{eq:threshold2}
\end{equation}
where $\delta f_{\ell}$ denotes the local global range of $f$ within
the $k$-ring neighborhood for WALF reconstruction, $h_{\ell}$ is
the local average edge length in the local $uv$ coordinate system,
and $\delta f_{g}$ and $h_{g}$ are the same as those in (\ref{eq:node-based-indicator}).
$C_{\ell}$ and $C_{g}$ are two parameters, which we determine empirically.

\begin{figure}
\subfloat[$f_{3}$ with $C^{0}$ and $C^{1}$ discontinuities.]{\begin{centering}
\includegraphics[width=0.48\columnwidth]{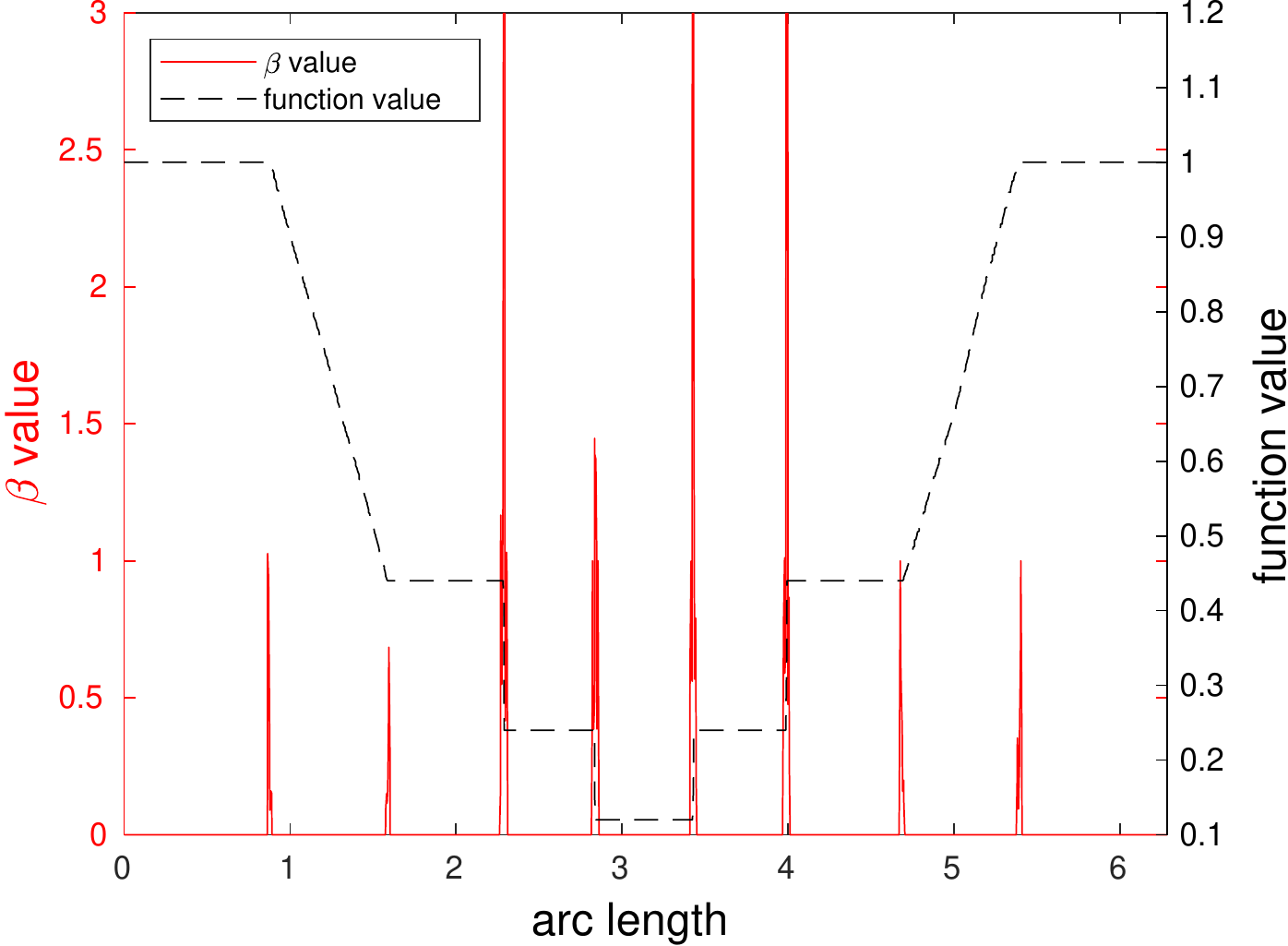}
\par\end{centering}
}\hfill\subfloat[$f_{4}$ with $C^{0}$, $C^{1}$, and $C^{2}$ discontinuities.]{\begin{centering}
\includegraphics[width=0.48\columnwidth]{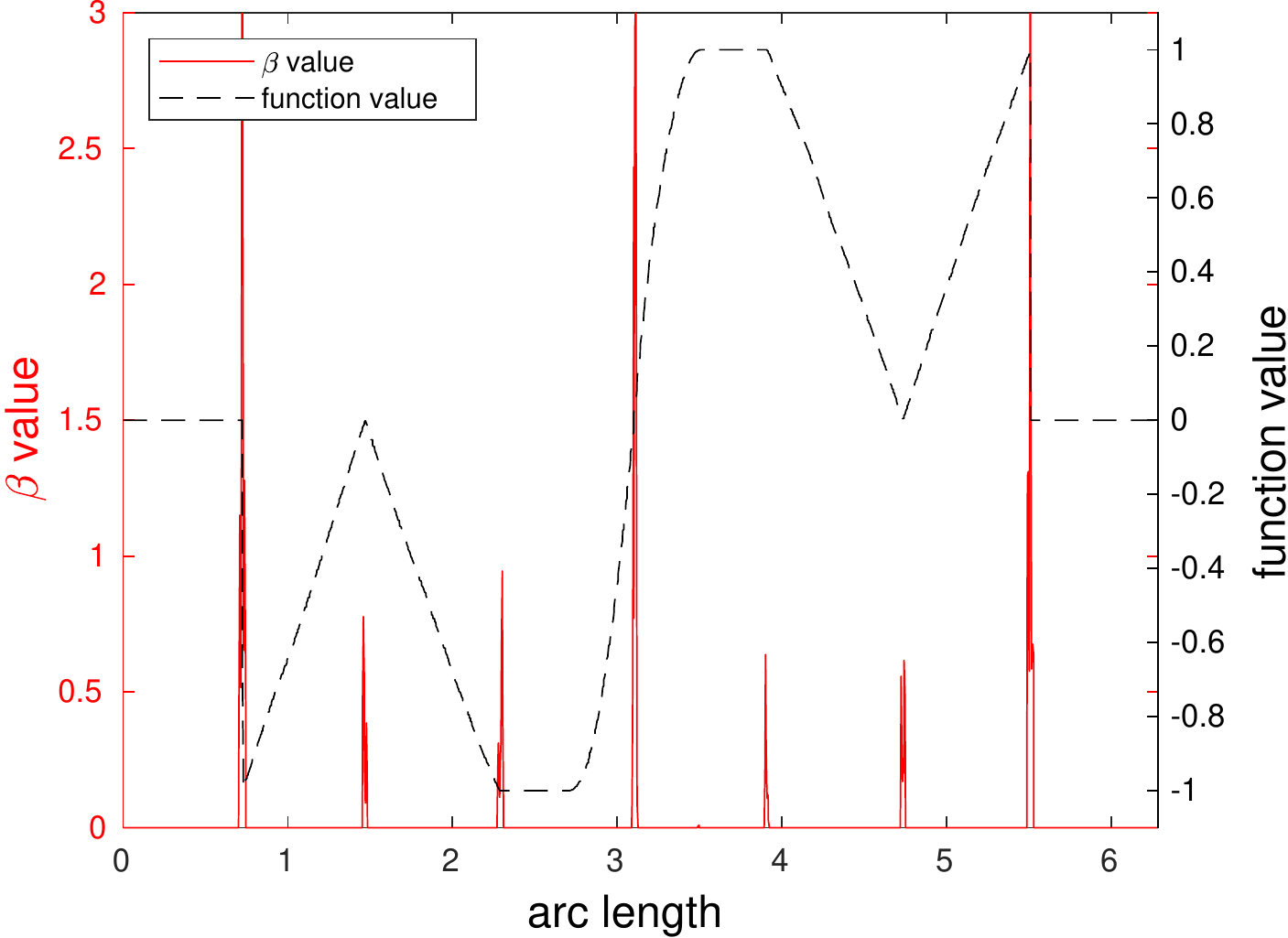}
\par\end{centering}
}

\caption{\label{fig:Example-node-based-indicator}Node-based $\beta$ indicator
values at discontinuities of example functions.}
\end{figure}

We derived the thresholds based on some asymptotic analysis and numerical
experimentation. First, let us focus on $\beta_{v}$. We observe that
$\beta_{v}\gtrsim3$ near $C^{0}$ discontinuities and $\beta_{v}\gtrsim0.5$
near $C^{1}$ discontinuities. This is evident in Figure~\ref{fig:Example-node-based-indicator},
which shows the $\beta$ values for two example functions (\ref{eq:f3})
and (\ref{eq:f4}) with $C^{0}$, $C^{1}$, and $C^{2}$ discontinuities
on a spherical triangulation shown in Figure~\ref{fig:Discontinuous-functions}.
These bounds can also be derived heuristically based on a 1D analysis,
which we omit here. For robustness, we chose $\kappa=0.3$, which
is small enough to identify $C^{1}$ discontinuities even on relatively
coarse meshes but large enough to avoid $C^{2}$ discontinuities on
relatively fine meshes. In addition, if the function is smooth and
convex (or concave), then no smooth points will be classified as discontinuities.
This is because the $\alpha$ values in these regions would have the
same signs, so $\beta\approx0$. However, at saddle points or inflection
points of a smooth function, the $\alpha$ values vary in signs, which
would lead to a nearly zero denominator in (\ref{eq:node-based-indicator})
and in turn a large $\beta$ value. This is partially the reason why
$\beta$ was large in the middle of Figure~\ref{fig:Example-node-based-indicator}(b).
Hence, using $\beta$ alone may lead to false positiveness, which
could activate the resolution of discontinuities at inflection points
and potentially reduce accuracy and efficiency.

To address the potential false positiveness, we introduced $\tau$
in (\ref{eq:threshold2}) as a second threshold. The local threshold
$\tau_{\ell}$ in (\ref{eq:threshold2}) is based on the Taylor series
analysis in Section~\ref{subsec:Element-based-indicators}: In smooth
regions, $\delta f_{\ell}=\mathcal{O}(h_{\ell})$, so $\delta f_{\ell}h_{\ell}^{0.5}=\mathcal{O}(h_{\ell}^{1.5})$.
This term separates the $\alpha$ values at smooth regions, which
are $\mathcal{O}(h_{\ell}^{2})$, from those at $C^{0}$ and $C^{1}$
discontinuities, which are $\mathcal{O}(1)$ and $\mathcal{O}(h_{\ell})$,
respectively. The use of $h_{\ell}$ (instead of $h_{g}$) makes $\tau_{\ell}$
to be independent of the scaling of the geometry and be insensitive
to the nonuniformity of the mesh. Similarly, the use of the local
range $\delta f_{\ell}$ makes $\tau_{l}$ independent of scaling
of $f$. However, if the function is locally nearly a constant, then
$\delta f_{\ell}\approx0$, so $\tau_{\ell}$ would fail to filter
out the false positiveness in such cases. In (\ref{eq:threshold2}),
$\tau_{g}$ is a safeguard for such cases, where $\delta f_{g}h_{g}^{1.5}=\mathcal{O}(h_{g}^{1.5})$
on quasiuniform meshes. In terms of the constant values $C_{\ell}$
and $C_{g}$, too large values can result in false negativeness, and
too small values can lead to some false positiveness. For robustness,
we chose $C_{\ell}=0.5$ and $C_{g}=0.05$, which are large enough
to avoid virtually all false positiveness on sufficiently fine meshes
without introducing false negativeness.

\subsubsection{Node-based markers on target mesh}

After obtaining the node-based markers, it is then straightforward
to map them to nodal markers on the target mesh. In particular, given
a target node $v$, if any of the source node in its stencil is marked
as discontinuities, we mark $v$ as a discontinuity on the target
mesh. This step further improves the robustness of the detection step,
in that some isolated false negativeness on the source mesh can be
corrected if any of their neighbor nodes is marked correctly. 

With a robust marker of discontinuities, one can apply various limiters
near discontinuities. For example, to preserve monotonicity near discontinuities,
we can bound the solution at a target node near discontinuities to
be bounded by the local extreme values of its containing element on
the source mesh. However, using such a limiter alone is insufficient
by itself. It is desirable to adapt the weighting scheme near the
detected discontinuities.

\subsection{\label{subsec:wls-eno}New weighting scheme for WLS-ENO}

After detecting the discontinuous regions, we resolve the potential
Gibbs phenomena from WLS with scaled Buhmann weights. Our basic idea
is motivated by WLS-ENO in \citep{liu2016wls}, which adapts the weights
to deprioritize the nodes that are on the other side of discontinuities.
 We propose a modification to the weighting scheme to take into account
the element-based $\alpha$ values. In particular, for the $j$th
node in the stencil of node $v$ on the target mesh, we define the
weight to be
\begin{equation}
\omega_{j}^{\text{ENO}}=\frac{\gamma_{j}^{+}\left(r_{j}^{2}+\epsilon_{\text{ID}}\right)^{-1/4}}{c_{0}\left|f_{j}-g_{0}\right|^{2}+c_{1}\delta f_{g}\max_{v_{j}\in e}\left|\alpha_{e}\right|+\epsilon_{\text{ENO}}\delta f_{g}^{2}\bar{h}^{2}},\label{eq:wls-eno-weights}
\end{equation}
where $f_{j}$ denotes the function value at the node in the source
mesh, $g_{0}$ denotes the linear interpolation at node $v$ from
the source mesh, and $\epsilon_{\text{ENO}}\bar{h}^{2}$ serves a
safeguard against division by a too small value. The scaling by the
different powers of $\delta f_{g}$ in the denominator consistent
units of all the terms for smooth functions, so that the weighting
scheme is invariant of shifting and scaling of the input function.
We use $\epsilon_{\text{ENO}}=10^{-3}$ and let $\bar{h}$ be the
average distance length in the local $uv$ coordinate system. In the
enumerator, $\gamma_{j}^{+}$, $r_{j}$, and $\epsilon^{\text{ID}}$
are the same as those in the inverse-distance based weight (\ref{eq:inverse-dist-weight}). 

The definition of $\omega_{j}^{\text{ENO}}$ was guided by an asymptotic
analysis, which we outline as follows. Without loss of generality,
let us first assume the enumerator is $1$. If $c_{0}=1$ and $c_{1}=0$,
then $\omega_{j}^{\text{ENO}}$ is essentially the same as that in
\citep{liu2016wls} for hyperbolic conservation laws, except that
we use linear interpolation to obtain $g_{0}$ in remap. With those
parameters, $\omega_{j}^{\text{ENO}}=\mathcal{O}(1)$, $\mathcal{O}(h^{-2})$
and $\mathcal{O}(h^{-2})$ at $C^{0}$ discontinuities, $C^{1}$ discontinuities,
and smooth regions, respectively, so it could not distinguish $C^{1}$
discontinuities from smooth regions. When $c_{1}\neq0$, however,
since $\max_{v_{j}\in e}\left|\alpha_{e}\right|=\mathcal{O}(h)$ at
$C^{1}$ discontinuities, $\omega_{j}^{\text{ENO}}=\mathcal{O}(1)$,
$\mathcal{O}(h^{-1})$ and $\mathcal{O}(h^{-2})$ at $C^{0}$ discontinuities,
$C^{1}$ discontinuities, and smooth regions, respectively. Note that
the denominator of (\ref{eq:wls-eno-weights}) does not take into
account distance or normals. As in (\ref{eq:inverse-dist-weight}),
$\gamma_{j}^{+}$ serves as a safeguard for sharp features or under-resolved
high-curvature regions in surface meshes. The inverse-distance part
of the enumerate introduces an $\mathcal{O}(\sqrt{h})$ term, which
would still allow $\omega_{j}^{\text{ENO}}$ to distinguish discontinuities
from smooth regions while reducing the influence of far-away points
and in turn reduce the potential interference of nearby discontinuities.

To obtain the parameters in $\omega_{j}^{\text{ENO}}$, we conducted
numerical experimentation and found that $c_{0}=1$ and $c_{1}=0.05$
worked well in practice. For robustness and efficiency, we use quadratic
WLS-ENO near discontinuities with $3$-rings, same as the ring size
for quartic WLS-ENO in smooth regions. With these parameters, our
proposed WLS-ENO scheme is relatively insensitive to false positiveness
of detected discontinuities. This is because even if all the nodes
are marked as discontinuities, the WLS-ENO weights (\ref{eq:wls-eno-weights})
with quadratic WLS would be used, and the method is still third-order
accurate for smooth functions, which is still higher than commonly
used low-order methods.

\section{\label{sec:Numerical-results}Numerical experimentation}

\begin{figure}
\begin{centering}
\includegraphics[width=0.45\columnwidth]{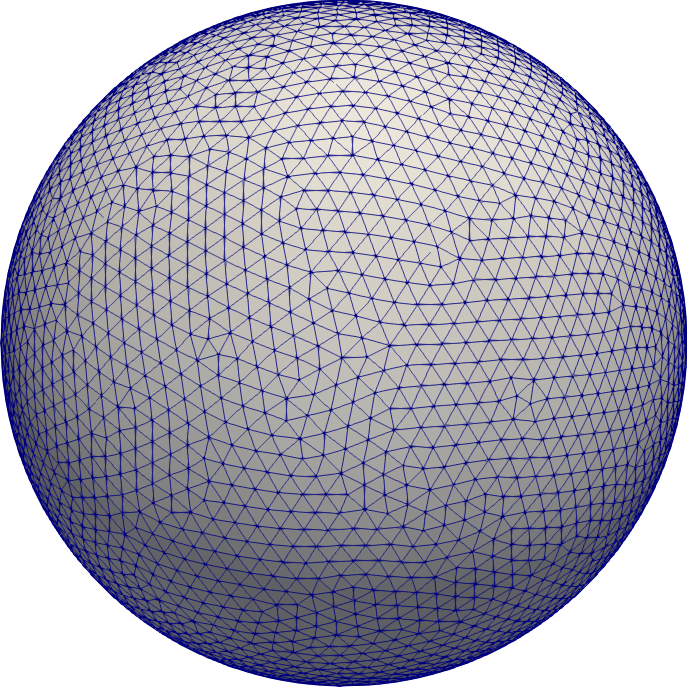}\hfill\includegraphics[width=0.45\columnwidth]{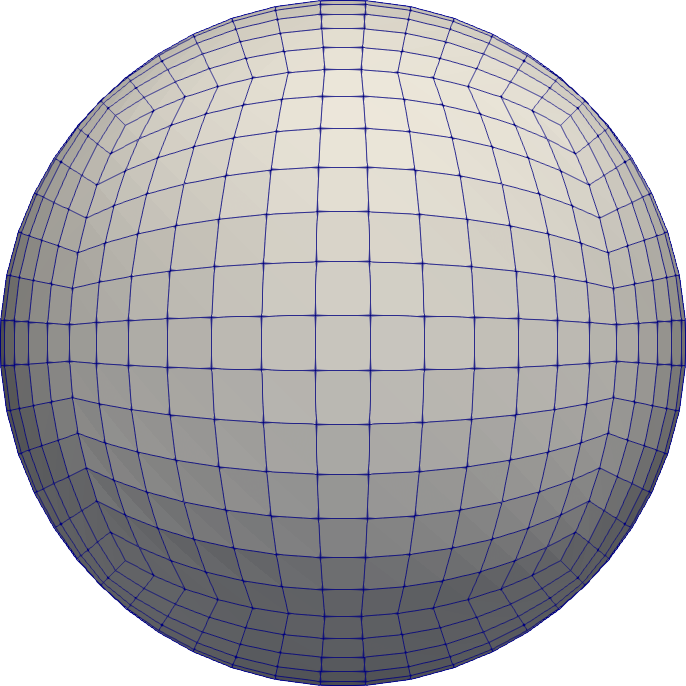}
\par\end{centering}
\caption{\label{fig:Example-triangular-and}Coarsest Delaunay and cubed-sphere
meshes used in our tests.}
\end{figure}

In this section, we report numerical experimentation with WLS-ENOR.
Since one of our main targeted application areas is climate and earth
modeling, we focus on on node-to-node transfer between Delaunay triangulations
and cubed-sphere meshes of the unit sphere. The former type is the
dual of the spherical centroidal Voronoi tessellations (SCVT) \citep{ju2011voronoi},
which are commonly used with finite volume methods in ocean modeling.
In this setting, the nodal values in the Delaunay mesh are equivalent
to cell-centered values in the Voronoi tessellation. The latter type
is often used in finite volume \citep{putman2007finite} or spectral
element methods \citep{terai2018atmospheric} in atmospheric models.
We assess WLS-ENOR for both smooth and discontinuous functions under
different mesh resolutions. To this end, we generated a series these
meshes of different resolutions using Ju's SCVT code \citep{ju2011voronoi}
and equidistant gnomonic projection of sphere \citep{sadourny1972conservative},
respectively.  Figure~\ref{fig:Example-triangular-and} shows the
coarsest meshes used in our tests, which we refer to as the level-1
meshes. The Delaunay mesh has 4,096 nodes and 8,188 triangles, and
the cubed-sphere mesh has 1,016 nodes and 1,014 quadrilaterals, respectively.
For each level of mesh refinement, the number of triangles or quadrilaterals
quadruples.

\subsection{\label{subsec:Transferring-smooth-functions}Transferring smooth
functions}

We first consider WLS-based transfer for smooth functions. We report
some representative results with two functions, namely, 
\begin{equation}
f_{1}(x,y,z)=\left(\sin(\pi x)+\cos(\pi y)\right)z\label{eq:trignometric function}
\end{equation}
and
\begin{equation}
f_{2}(x,y,z)=\left(11z^{2}-1\right)\left(x^{4}-6x^{2}y^{2}+y^{2}\right).\label{eq:spherical-harmonic}
\end{equation}
 The latter is the real part of an unnormalized version of the degree-$6$
spherical harmonic function $Y_{6}^{4}(\theta,\varphi)$ \citep{wiki:sphereical_harmonics},
i.e.,
\begin{equation}
f_{2}\left(\theta,\varphi\right)=\text{Re}\left(e^{4i\varphi}\sin^{4}\theta\left(11\cos^{2}\theta-1\right)\right),\label{eq:spherical-harmanic-alt}
\end{equation}
where $\theta\in[0,\pi${]} and $\varphi\in[0,2\pi]$ are the polar
(colatitudinal) and the azimuthal (longitudinal) angles, respectively.
For smooth functions, we consider the $\ell^{2}$-norm error. Let
$\vec{e}$ denote the vector composed of the pointwise error for each
node of a mesh, and let $N$ denote the number of nodes. The $\ell^{2}$-norm
error is measured as
\begin{equation}
\left\Vert \vec{e}\right\Vert _{\ell_{2}}=\frac{1}{\sqrt{N}}\left\Vert \vec{e}\right\Vert _{2}=\sqrt{\frac{1}{N}\sum_{i=1}^{N}e_{i}^{2}}.\label{eq:l2-pointwise}
\end{equation}
Given two meshes of different resolutions with $N_{1}$ and $N_{2}$
nodes, where $N_{1}<N_{2}$, and let $\vec{e}_{1}$ and $\vec{e}_{2}$
denote their respective error vectors. Assuming the meshes are uniformly,
which is the case in our tests, we evaluate the convergence rate in
$\ell_{2}$ norm between the two meshes as 
\begin{equation}
\text{convergence rate}=2\frac{\log(\left\Vert \vec{e}_{1}\right\Vert _{\ell_{2}}/\left\Vert \vec{e}_{2}\right\Vert _{\ell_{2}})}{\log(N_{2}/N_{1})},\label{eq:convergence-rate}
\end{equation}
which is equivalent to the standard definition based on edge length.
Although we recommend quartic WLS in smooth regions, we will report
results with quadratic, quartic, and sextic WLS for completeness and
for comparison.

\subsubsection{\label{subsec:Comparison-of-weighting-schemes}Numerical optimization
of cut-off radii}

\begin{figure}
\subfloat[Error profiles for $f_{1}$ in (\ref{eq:trignometric function}).]{\begin{centering}
\includegraphics[width=0.48\columnwidth]{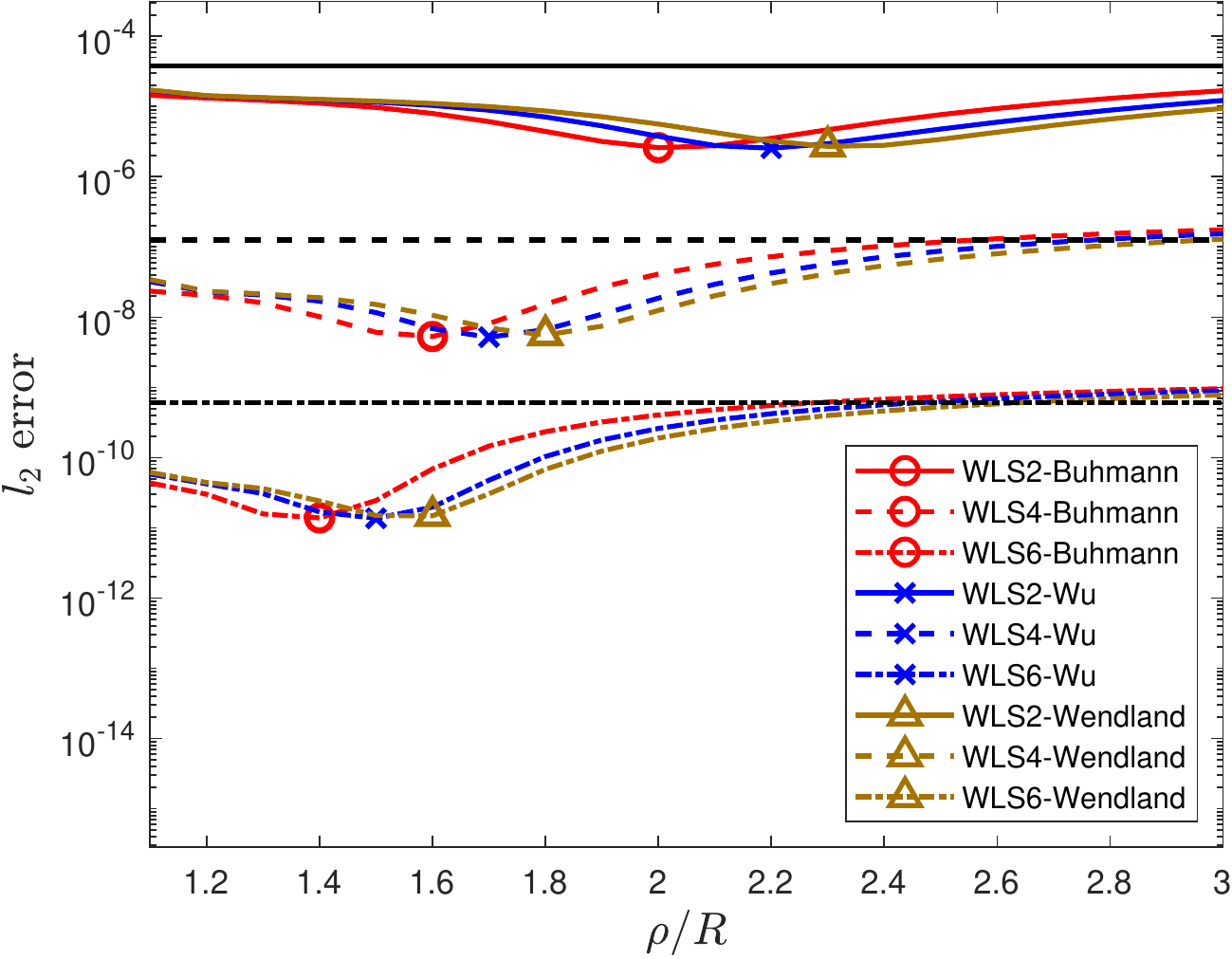}
\par\end{centering}
}\hfill\subfloat[Error profiles for $f_{2}$ in (\ref{eq:spherical-harmonic}).]{\begin{centering}
\includegraphics[width=0.48\columnwidth]{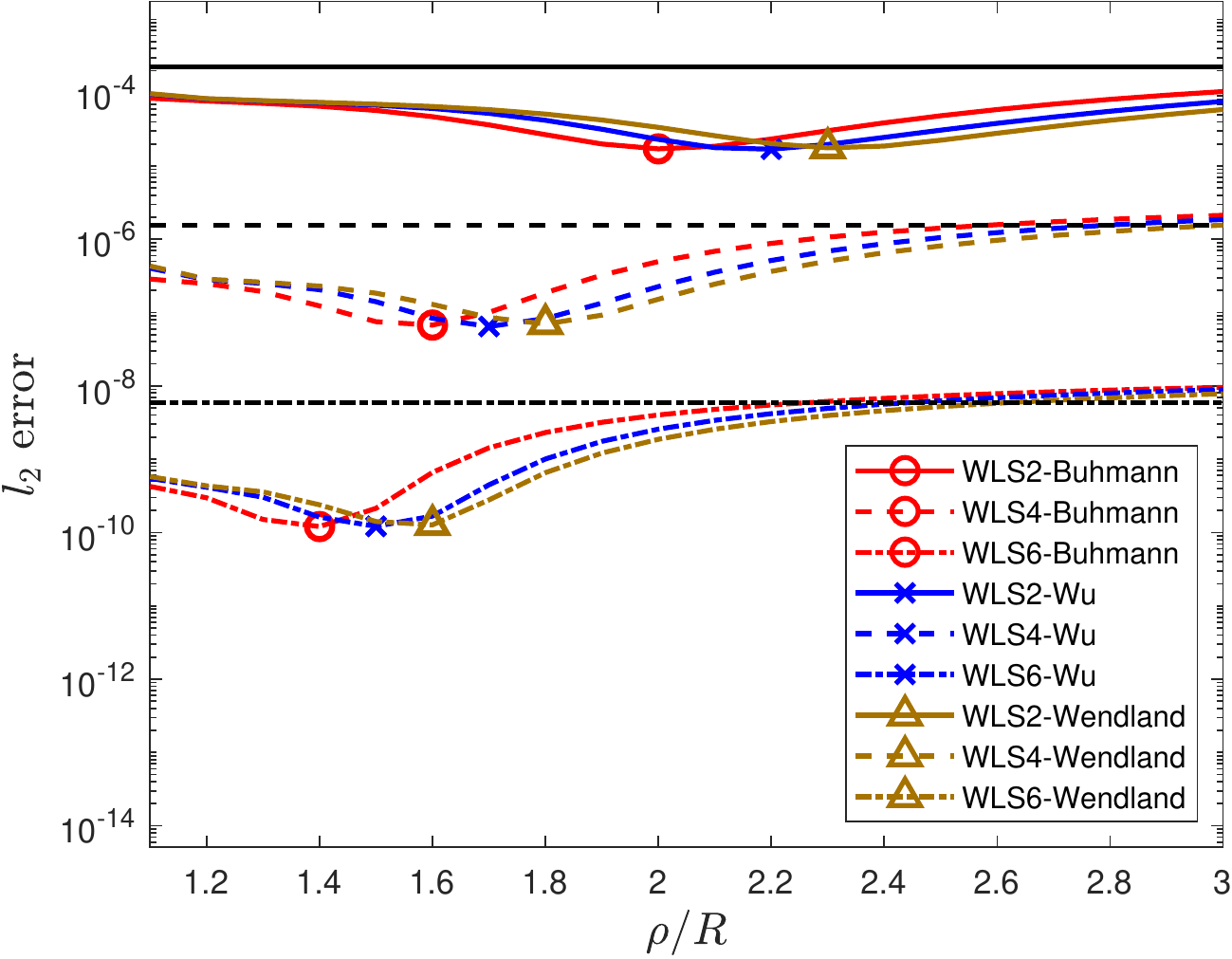}
\par\end{centering}
}\caption{\label{fig:Accuracy-versus-radius}$\ell^{2}$-norm errors in degree-2,
4, and 6 WLS with respect to the radius ratio $(\sigma=\rho/R)$ in
weights based on Buhmann's $C^{3}$, Wu's $C^{4}$, and Wendland's
$C^{4}$ radial functions. The marks indicate the optimal $\sigma$
values for each weighting strategy.}
\end{figure}

One of the key aspect of our WLS transfer is its weighting scheme,
namely the scaled Buhmann weights described in Section~\ref{subsec:Optimizing-weights}.
Here, we describe our procedure to optimize the cut-off radius $\rho=\sigma R$
in these weights. Our goal is to minimize the $\ell^{2}$-norm error
of the transferred solution with respect to $\sigma$. It is difficult,
if not impossible, to write down an equation for this objective function,
because $R$ depends on the combinatorial structure of the mesh. Fortunately,
a precise $\sigma$ is unnecessary, so we solve this optimization
problem approximately by using the two example $f_{1}$ and $f_{2}$
in (\ref{eq:trignometric function}) and (\ref{eq:spherical-harmonic})
and the Delaunay and cubed-sphere meshes depicted in Figure~\ref{fig:Example-triangular-and}.
We transferred the functions between the meshes for $\sigma$ between
1 and 3 with an increment of $0.1$. Figure~\ref{fig:Accuracy-versus-radius}
shows the $\ell^{2}$ norm errors of degree-2, 4, and 6 WLS in transferring
from the Delaunay to the cubed-sphere mesh. It is evident that the
error profile had a ``V'' shape with respect to $\sigma$ for both
functions, and the minima were at $\sigma=2.0$, $1.6$ and $1.4$
for degree-2, 4, and 6 WLS, respectively. As a reference, Figure~\ref{fig:Accuracy-versus-radius}
also shows the errors with the inverse-distance based weights (\ref{eq:inverse-dist-weight}).
The scaled Buhmann weights with optimal $\sigma$ improved the accuracy
by nearly two orders of magnitude compared to using (\ref{eq:inverse-dist-weight}).
It is also worth noting that the optimal $\sigma$ values remained
about the same when we transferred other smooth functions, transferred
from the cubed-sphere to the Delaunay meshes, or used different mesh
resolution. This shows that the numerical optimization is well posed
for smooth functions with a sufficiently smooth weight function.

We also applied the optimization procedure to other radial functions,
including those in \citep{wu1995compactly} and \citep{wendland1995piecewise}.
We observed that the well-posedness of our optimization procedure
required at least $C^{3}$ continuity. For comparison, Figure~\ref{fig:Accuracy-versus-radius}
also shows the error profiles with the following two weighting schemes:
\begin{equation}
w_{j}=\gamma_{j}^{+}\phi_{2,0}(r_{j}/\rho),\,\text{where}\,\phi_{2,0}(r)=\left(1-r\right)_{+}^{5}\left(r^{4}+5r^{3}+9r^{2}+5r+1\right),\label{eq:Wu4}
\end{equation}
and
\begin{equation}
w_{j}=\gamma_{j}^{+}\psi_{4,2}(r_{j}/\rho),\text{ where }\psi_{4,2}(r)=\left(1-r\right)_{+}^{6}\left(35r^{2}+18r+3\right).\label{eq:Wendland4}
\end{equation}
The radial functions $\phi_{2,0}$ and $\psi_{4,2}$ are $C^{4}$
functions due to Wu \citep{wu1995compactly} and Wendland \citep{wendland1995piecewise},
respectively. It is clear that the scaled Buhmann weights had the
smallest optimal cut-off radii, which result in smaller constant factors
in the leading error term and also fewer rows in the generalized Vandermonde
systems. The $C^{6}$ functions of Wu and Wendland required even larger
cut-off radii. Hence, the scaled Buhmann weights are the overall winners
in both accuracy and efficiency for WLS-based remapping. We note that
DTK contains a collection of weights, which were mostly based on Wu's
and Wendland's functions, and DTK requires the user to choose the
specific weights and the cut-off radius \citep{slattery2016mesh}.
In contrast, our approach is parameter-free from the user's perspective,
in that both the weighting functions and cut-off radii were pre-optimized.

\subsubsection{Accuracy and convergence}

\begin{figure}
\begin{raggedright}
\subfloat[Errors for $f_{1}$ in (\ref{eq:trignometric function}).]{\begin{centering}
\includegraphics[width=0.5\columnwidth]{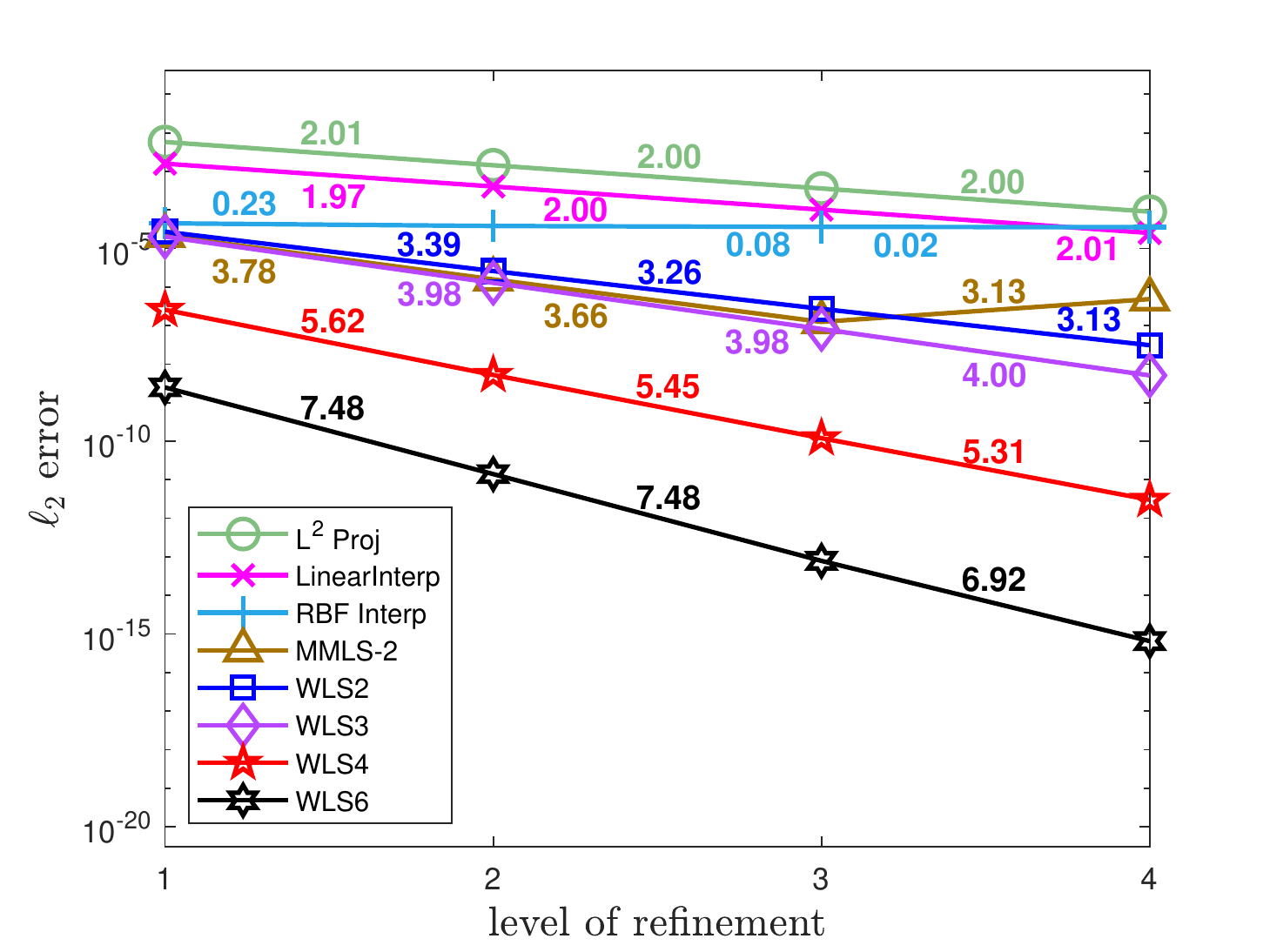}
\par\end{centering}
} \subfloat[Errors for $f_{2}$ in (\ref{eq:spherical-harmonic}).]{\begin{centering}
\includegraphics[width=0.5\columnwidth]{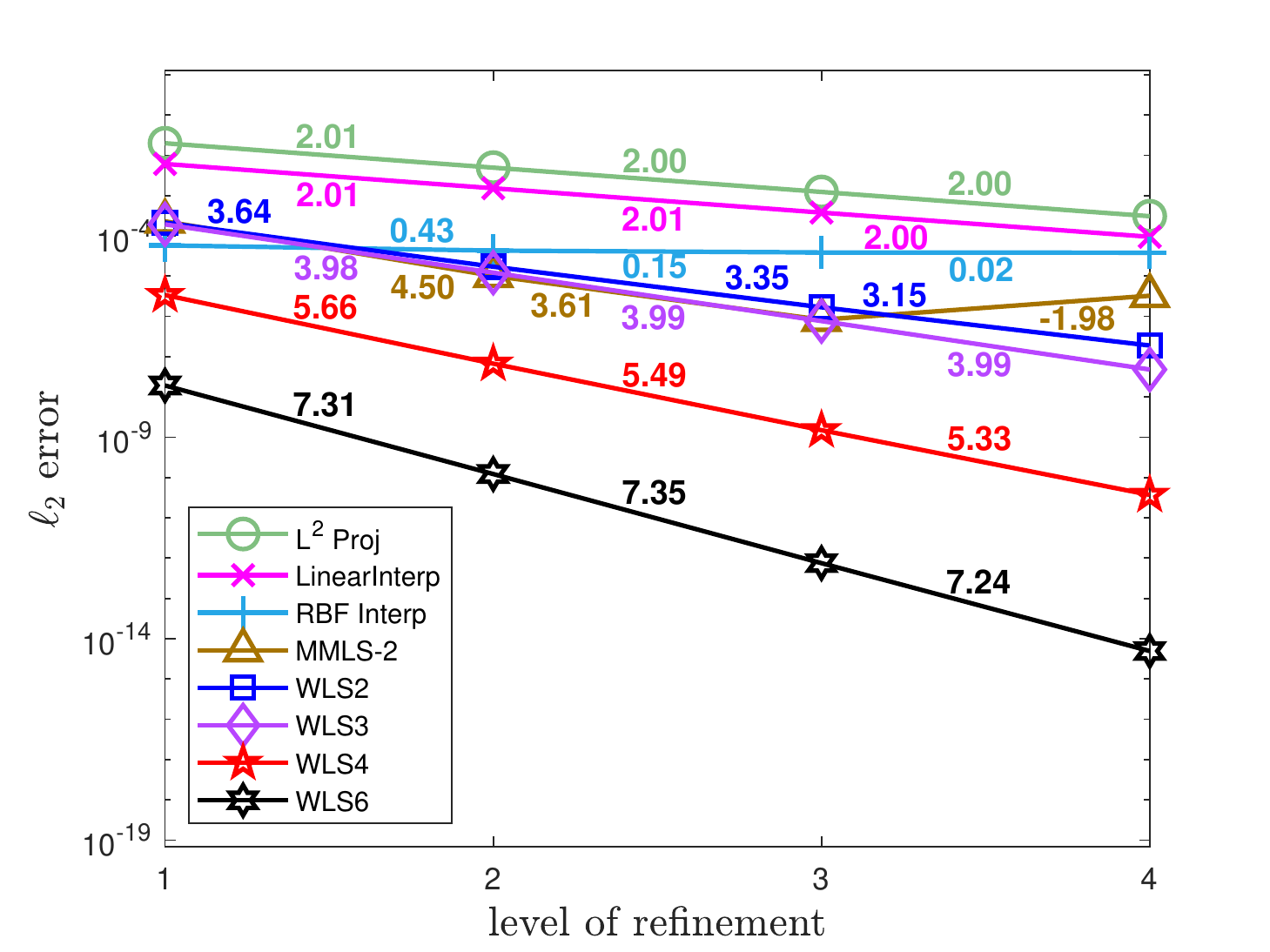}
\par\end{centering}
}
\par\end{raggedright}
\caption{\label{fig:Comparison-of-errors}Comparison of accuracy and convergence
of remap methods under mesh refinement.}
\end{figure}

We first assess the accuracy and convergence of WLS for smooth functions
in comparison with several other data-transfer methods, including
linear interpolation, common-refinement-based $L^{2}$ projection
\citep{jiao2004common}, MMLS in DTK \citep{slattery2016mesh}, and
RBF interpolation \citep{beckert2001multivariate}. For MMLS, we used
the same parameters as recommended in \citep{slattery2016mesh}, with
Wu's $C^{4}$ radial function as weights and a cut-off radius of about
five times the maximum edge length on the source mesh. For RBF interpolation,
we used the implementation in DTK as described in \citep{slattery2016mesh}.
However, we adhered to the original formulation of Beckert and Wendland
in \citep{beckert2001multivariate} with Wendland's $C^{2}$ radial
functions (instead of using Wu's $C^{4}$ functions as in \citep{slattery2016mesh}).
We set the cut-off radius for RBF to be five times the maximum edge
length on both source and target meshes, which was significantly larger
than the recommendation in \citep{beckert2001multivariate} but was
necessary for the solution to be stable and reasonably accurate. For
convergence study, we used four sets of triangular and quadrilateral
meshes with different resolutions. Figure~\ref{fig:Comparison-of-errors}
shows the $\ell^{2}$-norm errors of each of these methods together
with quadratic, quartic, and sextic WLS for functions $f_{1}$ and
$f_{2}$. For completeness, we also show the result with cubic WLS,
for which we used its optimal radius ratio $\sigma=1.2$. The number
associated with each line segment indicates the convergence rate for
its corresponding method under the corresponding mesh refinement.

We make several observations from the results. First, both $L^{2}$
projection and linear interpolation achieved second-order accuracy.
Note that on surface meshes, $L^{2}$ projection involves geometric
errors in its integration, so its accuracy was slightly worse that
that of linear interpolation, unlike the 2-D results in \citep{jiao2004common}.
The RBF interpolation of Beckert and Wendland produced smaller errors
than linear interpolation and $L^{2}$ projection, but its overall
convergence rate was close to zero. Note that MMLS in DTK uses quadratic
polynomials in the global $xyz$ coordinate system, which have the
same number of coefficients as cubic WLS in the local $uv$ coordinate.
As a result, MMLS performed slightly better than quadratic WLS on
the three coarser meshes but worse than cubic WLS. More importantly,
DTK lost convergence on the finest mesh, because the nearly planar
stencils on fine meshes lead to nearly singular Vandermonde systems.
Although DTK uses truncated SVD (TSVD) to resolve ill-conditioning
algebraically, the convergence was lost due to the random truncation
of low-degree terms by TSVD. In contrast, we use local 2D coordinate
systems within the local tangent space and use QRCP to solve the linear
system, and we truncate the highest-degree terms in the presence of
ill-conditioning. Hence, WLS is both stable and accurate, even when
using degree-6 polynomials, of which the solution achieved machine
precision on the finest mesh.

\subsubsection{Accuracy and conservation in repeated transfer}

\begin{figure}
\subfloat[Errors for $f_{1}$ in (\ref{eq:trignometric function}).]{\begin{centering}
\includegraphics[width=0.5\columnwidth]{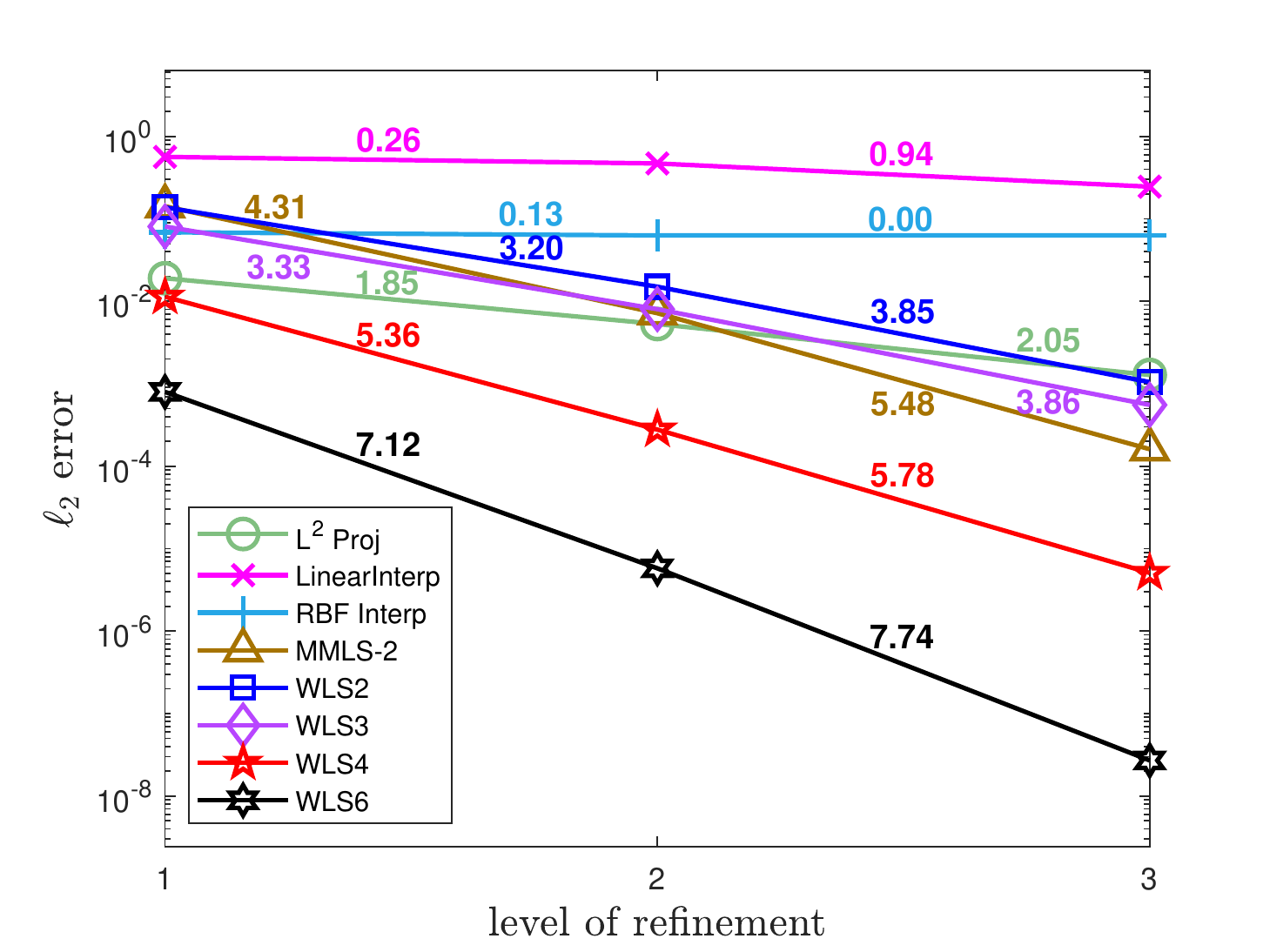}
\par\end{centering}
} \subfloat[Errors for $f_{2}$ in (\ref{eq:spherical-harmonic}).]{\begin{centering}
\includegraphics[width=0.5\columnwidth]{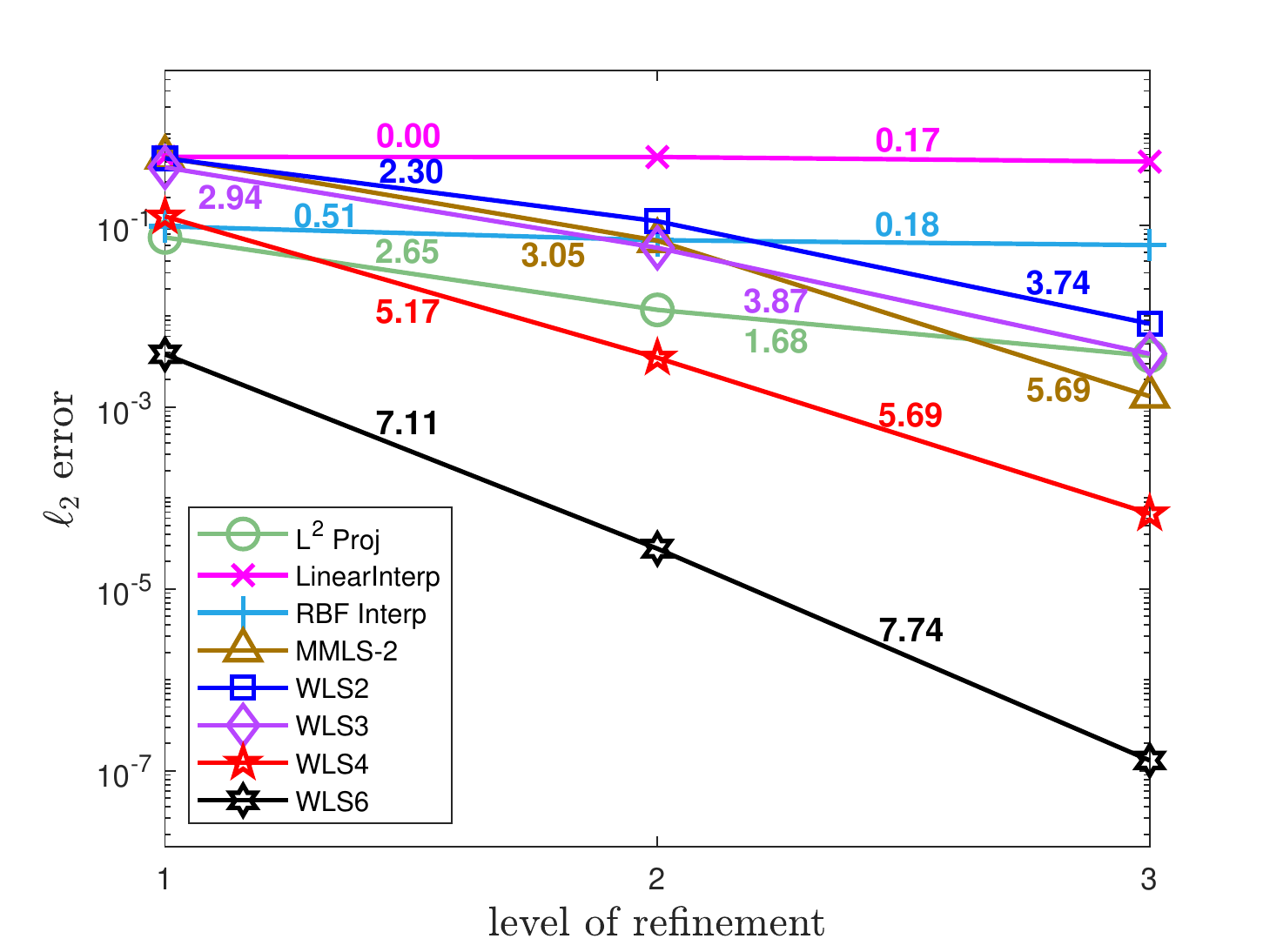}
\par\end{centering}
}\caption{\label{fig:repeted-transfer-smooth-accuracy}Accuracy and convergence
on Delaunay mesh after 1000 steps of repeated transfers.}
\end{figure}

The preceding section focused on accuracy and convergence in a single-step
transfer. In multiphysics applications, data must be exchanged between
meshes repeated. Hence, it is important to take into account the accuracy
and convergence in such a setting. Figure~\ref{fig:repeted-transfer-smooth-accuracy}
compared the accuracy when transferring smooth functions back and
forth for 1000 times between the Delaunay and cubed-sphere meshes.
Because MMLS was unstable and RBF was computationally too expensive
on the finest mesh, we only report results on the three coarser meshes.
Quartic and sextic WLS clearly stand out in their accuracy. Among
the other methods, linear interpolation became far less accurate in
repeated transfer. RBF interpolation performed better than linear
interpolation, but its convergence rate was close to zero. $L^{2}$
projection was remarkably accurate compared to all the other third
or lower-order methods, especially on coarse meshes, due to its conservation
property. Note that MMLS exhibited better accuracy than quadratic
WLS but similar accuracy as cubic WLS in repeated transfer. This is
because MMLS with quadratic polynomials in the global $xyz$ coordinate
systems has the same number of monomial basis functions as WLS with
cubic polynomials in the local $uv$ coordinate system.

\begin{figure}
\subfloat[Level-2 mesh.]{\includegraphics[width=0.48\columnwidth]{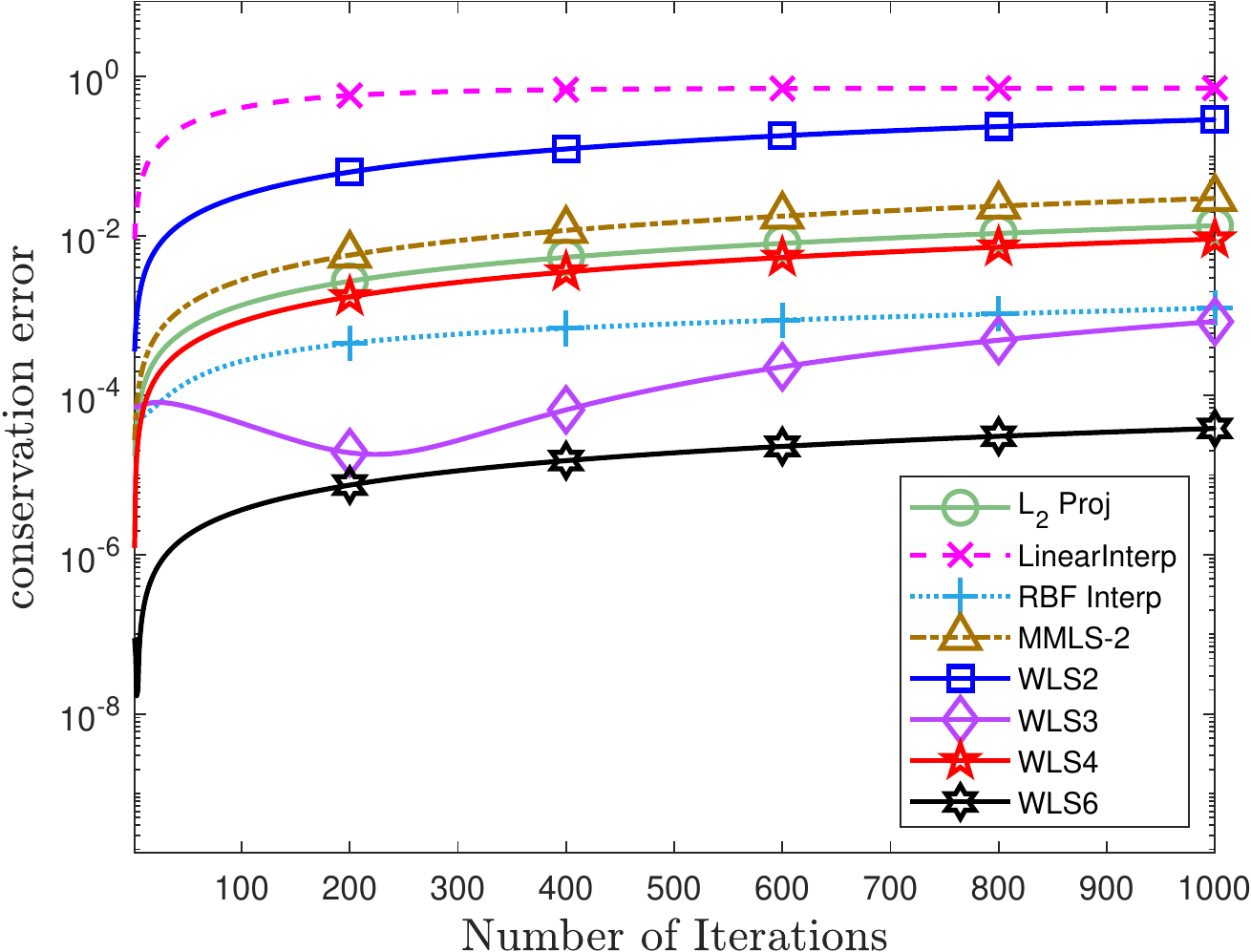}

}\hfill\subfloat[Level-3 mesh.]{\includegraphics[width=0.48\columnwidth]{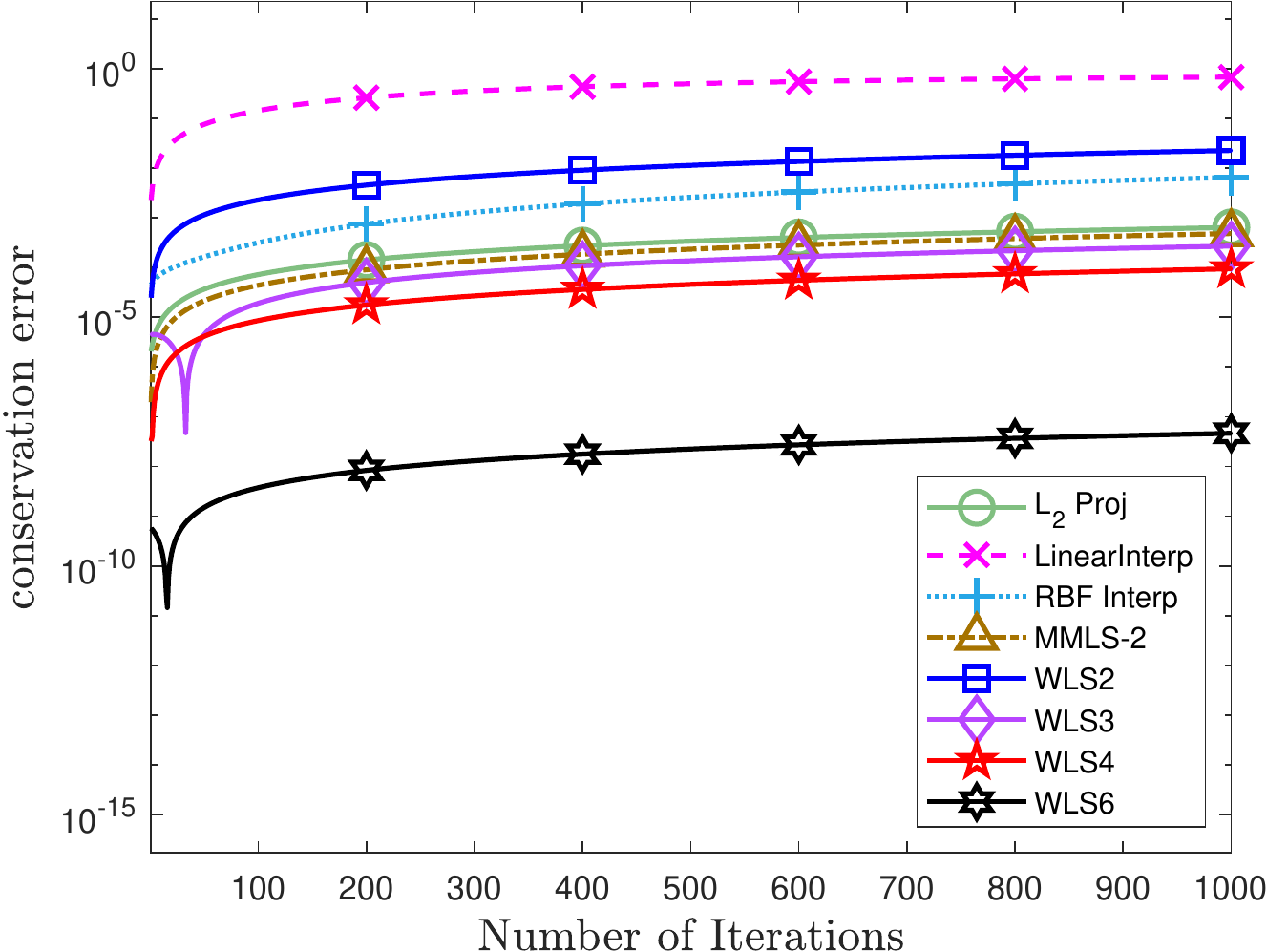}

}\caption{\label{fig:repeted-transfer-smooth-conservation}Evolution of conservation
errors in repeated transfer for $f_{2}$.}
\end{figure}

Besides accuracy, another important consideration for repeated transfer
is \emph{conservation} \citep{de2008comparison,jiao2004common}. Unfortunately,
the conservation error is not uniquely defined, especially for high-order
methods. In this work, we measure the conservation error as the difference
between the integral of the exact function and that of final reconstructed
function. We use numerical quadrature to compute the integrals the
exact spherical geometry, where we evaluated the values at the quadrature
points using the same-degree of WLS for quartic and sextic WLS and
using quadratic WLS reconstruction for all the other methods. Figure~\ref{fig:repeted-transfer-smooth-conservation}
compares the conservation errors on the level-2 and level-3 meshes
for $f_{2}$. Linear interpolation performed the worst due to its
low-order accuracy and non-conservation. $L^{2}$ projection was not
strictly conservative primarily because its numerical integration
is inexact on spheres, but it was nevertheless more conservative than
MMLS and WLS-2, although less conservative than WLS-3. RBF interpolation
was competitive with $L^{2}$ projection in terms of conservation,
despite its lower accuracy. Quartic and sextic WLS were more conservative
than $L^{2}$ projection on surfaces, although they do not enforce
conservation explicitly.

\subsection{\label{subsec:Discontinuous-functions}Transferring discontinuous
functions}

\begin{figure}
\begin{centering}
\includegraphics[width=0.48\columnwidth]{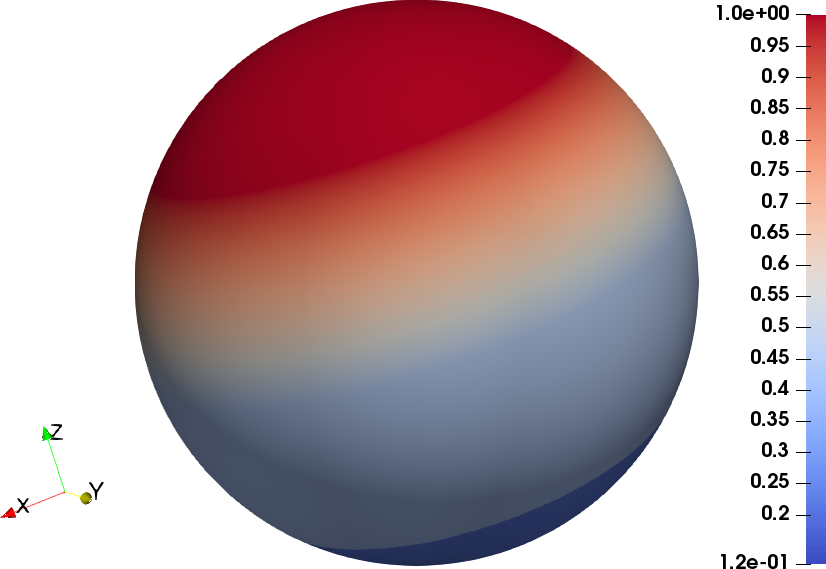} \includegraphics[width=0.48\columnwidth]{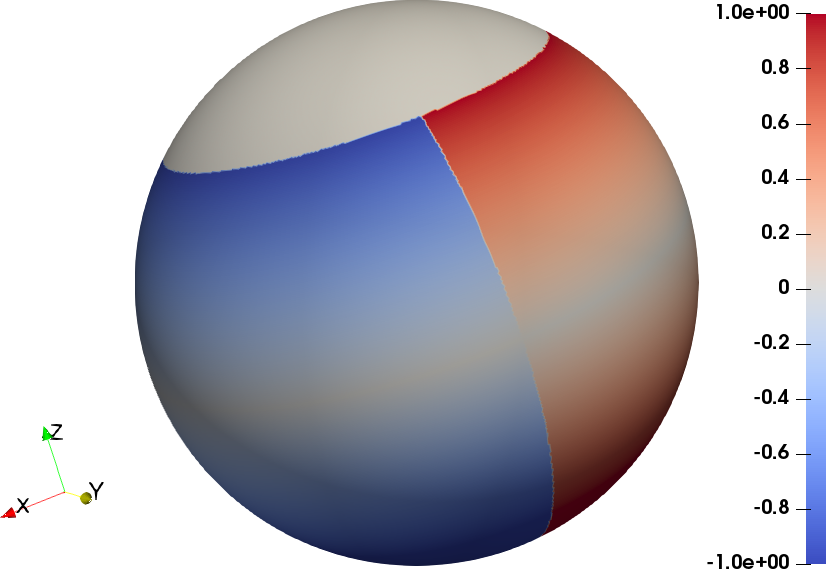}
\par\end{centering}
\caption{\label{fig:Discontinuous-functions}Discontinuous functions $f_{3}$
(left) and $f_{4}$ on level-4 Delaunay mesh.}
\end{figure}

\begin{figure}
\begin{centering}
\includegraphics[width=0.4\columnwidth]{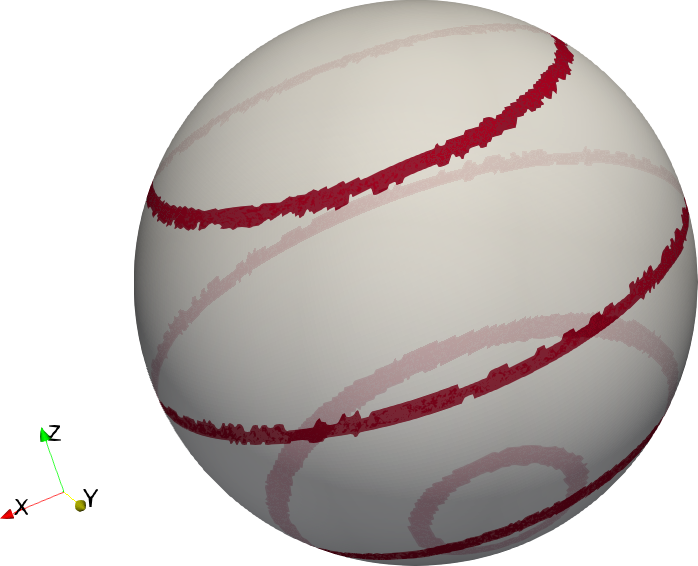}
\hspace{1.2cm}\includegraphics[width=0.4\columnwidth]{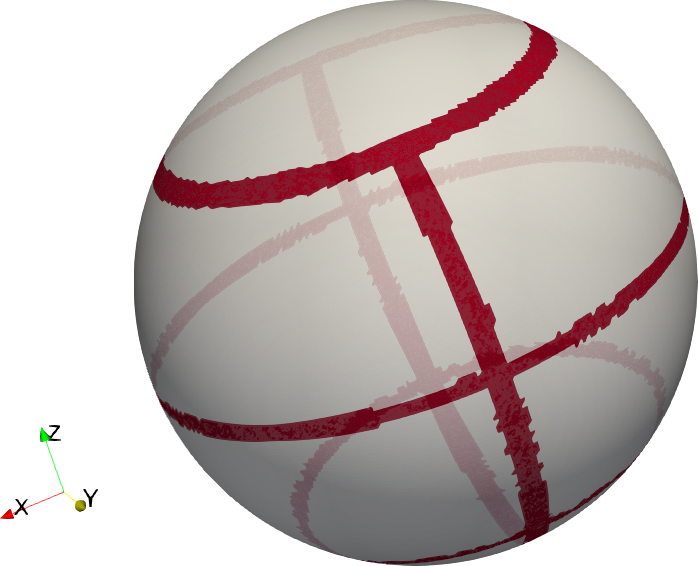}
\hspace{0.5cm}
\par\end{centering}
\caption{\label{fig:detected-discontinuities}Detected discontinuities for
$f_{3}$ (left) and $f_{4}$ on level-4 cubed-sphere mesh.}
\end{figure}

To assess the treatment of discontinuities, we use two piecewise smooth
functions, namely
\begin{equation}
f_{3}\left(\theta,\varphi\right)=\begin{cases}
\begin{array}{c}
1\\
1-0.8\left(\theta-0.87\right)\\
0.44\\
0.24\\
0.12
\end{array} & \begin{array}{c}
0\leq\theta<0.87\\
0.87\leq\theta<\pi/2\\
\pi/2\leq\theta<2.27\\
2.27\leq\theta<2.83\\
2.83\leq\theta\leq\pi
\end{array}\end{cases},\label{eq:f3}
\end{equation}
and
\begin{equation}
f_{4}\left(\theta,\varphi\right)=-1000+g(\varphi)\begin{cases}
0 & 0\leq\theta<\pi/4\\
-4\left(\theta/\pi-1/2\right) & \pi/4\le\theta<\pi/2\\
4\left(\theta/\pi-1/2\right) & \pi/2\le\theta<3\pi/4\\
1 & 3\pi/4\le\theta<7\pi/8\\
-64\theta^{2}/\pi^{2}+112\theta/\pi-48 & 7\pi/8\leq\theta\leq\pi
\end{cases},\label{eq:f4}
\end{equation}
with $g\left(\varphi\right)=\begin{cases}
-2000 & 0\leq\varphi<\pi\\
2000 & \pi\leq\varphi\leq2\pi
\end{cases}$. Function $f_{3}$ is an axial-symmetric function with both $C^{0}$
and $C^{1}$ discontinuities in the $\theta$ direction, which are
similar to an intermediate solution of two interacting blast waves
\citep{jiang1996efficient}. Hence, we refer to $f_{3}$ as \emph{interacting
waves}. Function $f_{4}$ has a richer structure of discontinuities,
including $C^{0}$, $C^{1}$, and $C^{2}$ discontinuities in the
$\theta$ direction, which intersect with the $C^{0}$ discontinuities
in the $\varphi$ direction. Hence, we refer to $f_{4}$ as \emph{crossing
waves}. The range of $f_{3}$ is between $0.12$ and 1, and that of
$f_{4}$ is approximately between $-3000$ and $1000$. The disparate
ranges allow us to demonstrate the independence of scaling and shifting
of the function values of our method. Figure~\ref{fig:Discontinuous-functions}
shows these functions on the level-4 Delaunay mesh. Figure~\ref{fig:detected-discontinuities}
shows the detected discontinuities on the level-4 cubed-sphere mesh,
computed from the indicators on the Delaunay mesh. All the $C^{0}$
and $C^{1}$ discontinuities were detected correctly on this mesh.
We note that they were also detected correctly on the coarser meshes.
These results demonstrate the low (nearly zero) false-negative rate
of our detection technique. In addition, our technique also has a
low false-positive rate, which is evident in Figure~\ref{fig:detected-discontinuities},
where the $C^{2}$ discontinuities in $f_{4}$ were not marked. We
also tested our discontinuity indicators for the smooth functions
$f_{1}$ and $f_{2}$ on all the meshes, and no false discontinuity
was marked on any of the meshes. Hence, our detection technique is
robust in terms of both false-negative and false-positive rates.

\subsubsection{Justification of WLS-ENO degree and weights}

\begin{figure}
\subfloat[After 500 steps on level-3 mesh.]{\begin{centering}
\includegraphics[width=0.48\columnwidth]{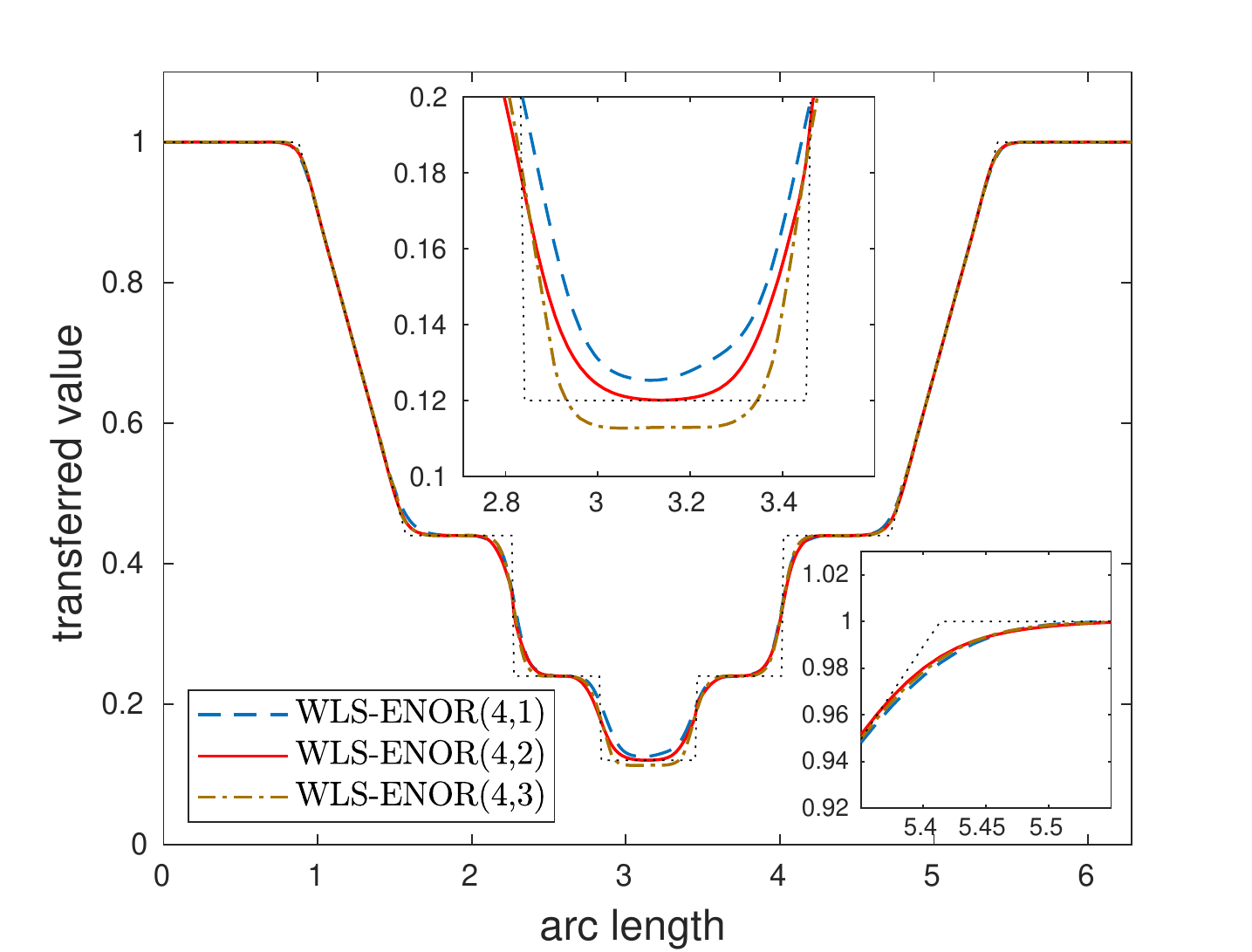}
\par\end{centering}
}\hfill\subfloat[After 1000 steps on level-3 mesh.]{\begin{centering}
\includegraphics[width=0.48\columnwidth]{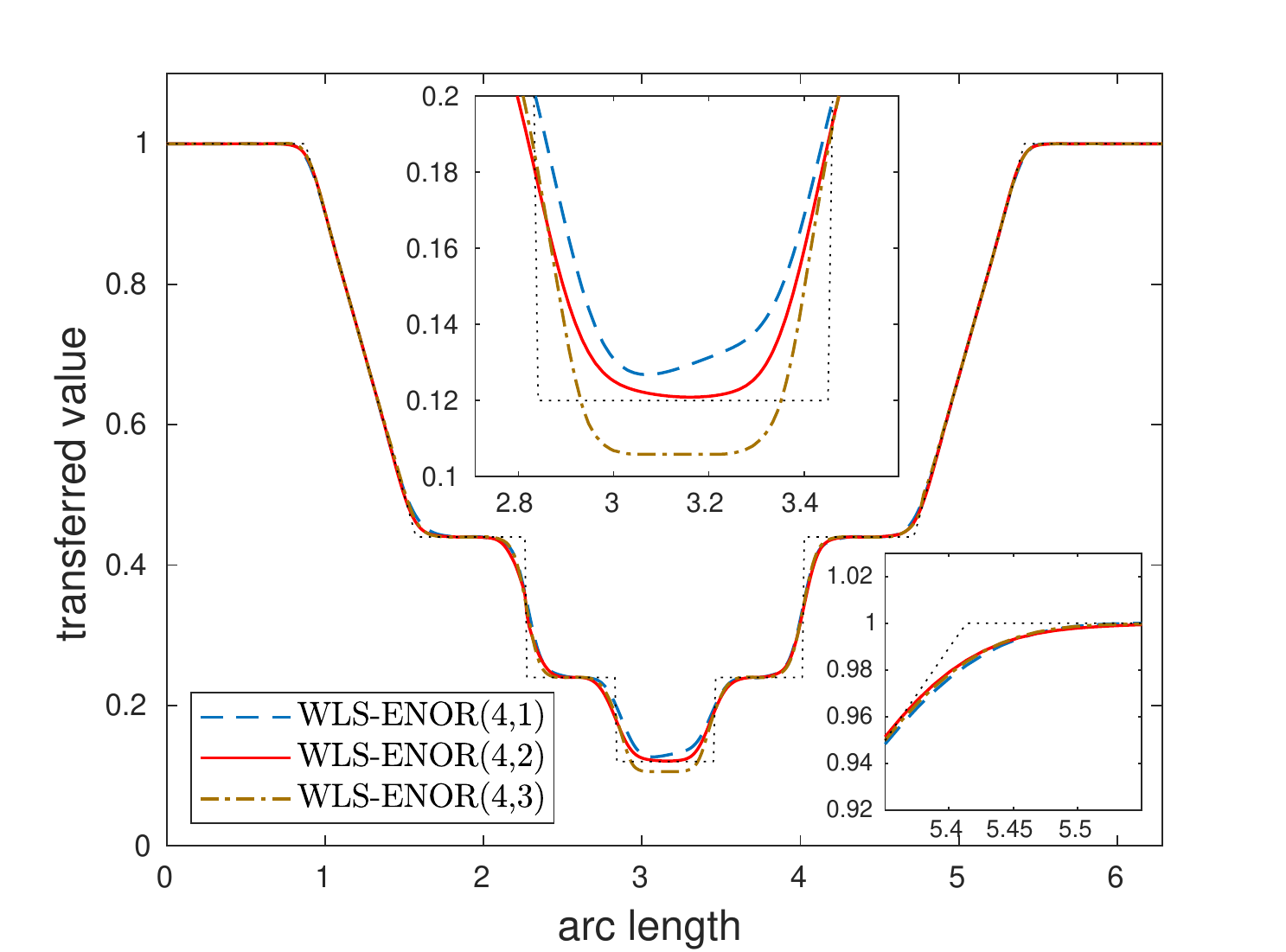}
\par\end{centering}
}

\caption{\label{fig:Comparison-of-WLS-ENOR}Comparison of WLS-ENOR(4,$q$)
for $q=1$, $2$, and $3$ in repeated transfer of $f_{3}$ between
level-3 Delaunay and cubed-sphere meshes.}
\end{figure}

In WLS-ENOR, we use quartic WLS in smooth regions and quadratic WLS-ENO
near discontinuities. We refer to this combination as \emph{WLS-ENOR(4,2)}.
As we have shown in Section~\ref{subsec:Transferring-smooth-functions},
quartic WLS delivers high accuracy and conservation for smooth functions.
We now justify the use of quadratic polynomials by comparing it with
WLS-ENOR(4,1) and WLS-ENOR(4,3), i.e., using quartic WLS in smooth
regions but linear or cubic WLS-ENO near discontinuities. For WLS-ENOR($p$,$q$),
we used $\lfloor1.5p\rfloor/2$-rings for smooth regions and $(q+0.5)$-rings
for discontinuities. We transferred $f_{3}$ back and forth between
the level-3 Delaunay and cubed-sphere meshes. Figure~\ref{fig:Comparison-of-WLS-ENOR}
shows the results after 500 and 1000 iterations along a great circle
passing through the poles, where the dotted solid line shows the ``exact''
solution. It is clear that WLS-ENOR(4,1) is overly diffusive at the
global minimum, although it worked well at other discontinuities.
On the other hand, WLS-ENOR(4,3) is less diffusive than WLS-ENOR(4,2)
near $C^{1}$ discontinuities, but it had a significant undershoot
at the global minimum where two $C^{0}$ discontinuities interact.
Note that this undershoot would vanish on the level-4 meshes, indicating
that WLS-ENOR(4,3) requires finer meshes than WLS-ENOR(4,2) to separate
$C^{0}$ discontinuities. Overall, WLS-ENOR(4,2) delivered a good
balance between robustness for coarser meshes and accuracy on finer
meshes, so we use it as the default. 

\begin{figure}
\subfloat[After 500 steps on level-4 mesh.]{\begin{centering}
\includegraphics[width=0.48\columnwidth]{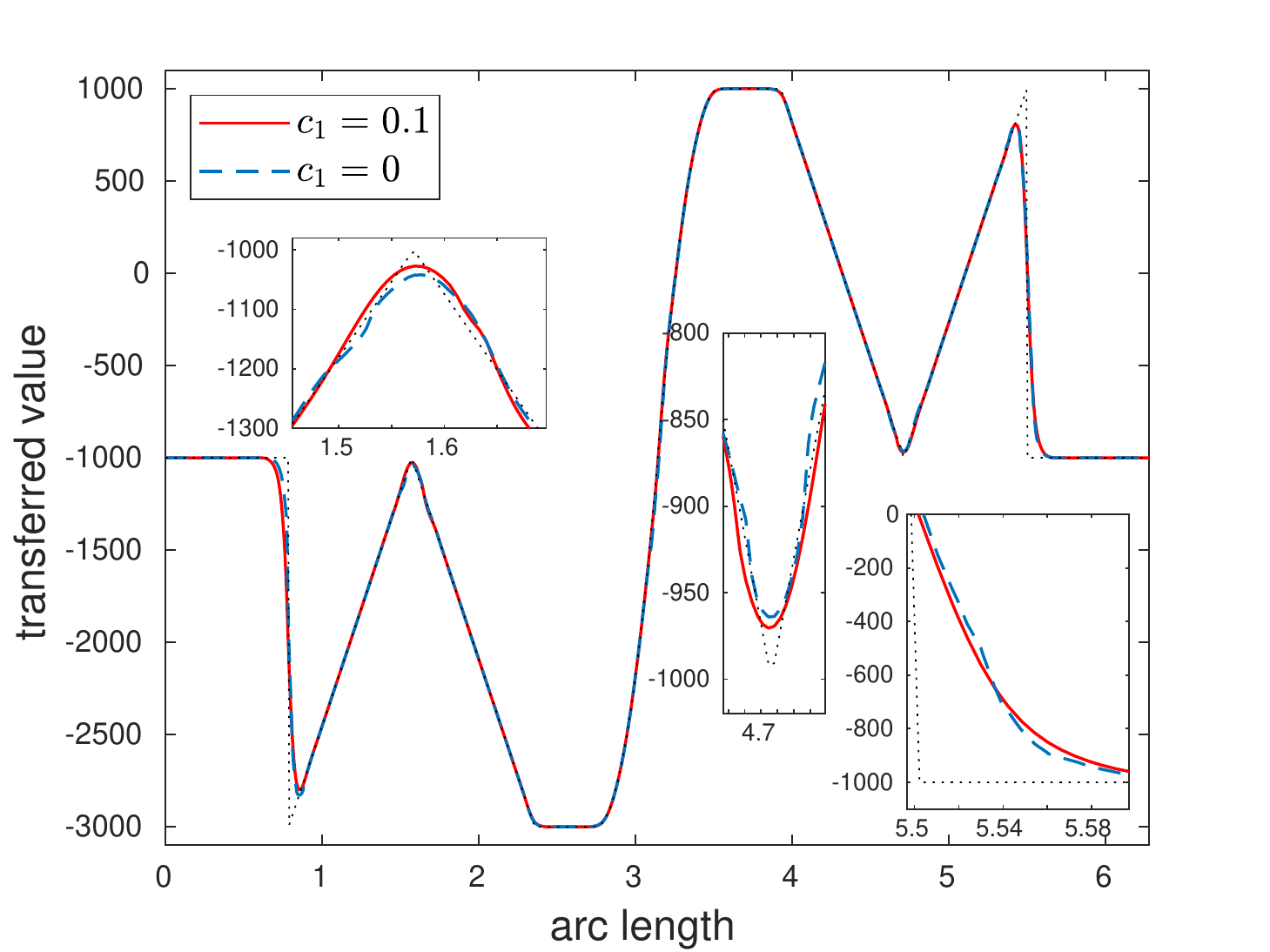}
\par\end{centering}
}\hfill\subfloat[After 1000 steps on level-4 mesh.]{\begin{centering}
\includegraphics[width=0.48\columnwidth]{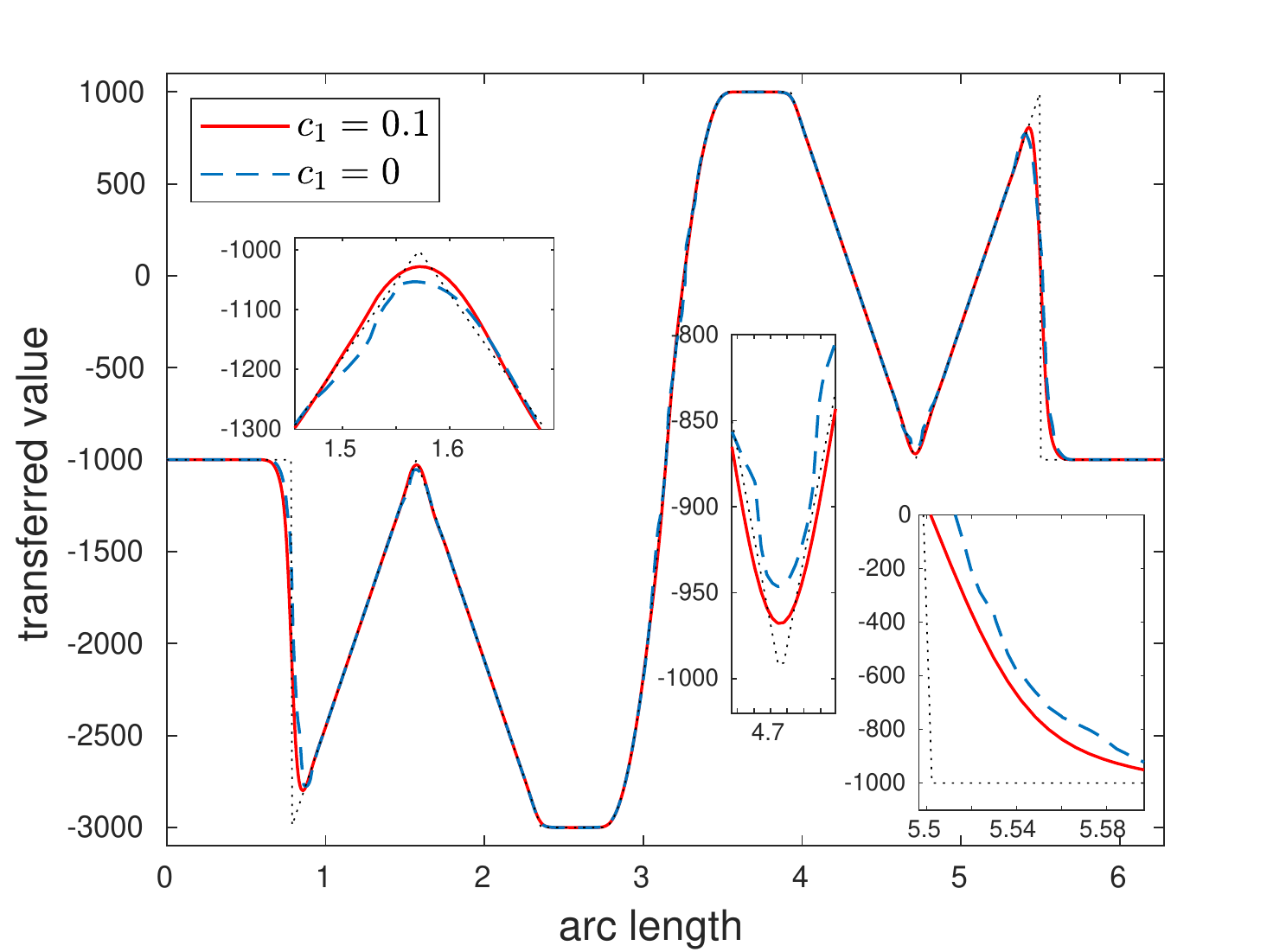}
\par\end{centering}
}

\caption{\label{fig:Comparison-of-WLS-ENO-weights}Comparison of WLS-ENOR with
and without explicit $C^{1}$ treatment in weighting scheme in repeated
transfer of $f_{4}$ between level-4 Delaunay and cubed-sphere meshes.}
\end{figure}

In WLS-ENOR, we introduced a new weighting scheme (\ref{eq:wls-eno-weights}),
which takes into account $C^{1}$ discontinuities and controls it
using the parameter $c_{1}$. To demonstrate its importance, Figure~\ref{fig:Comparison-of-WLS-ENO-weights}
compares WLS-ENOR with and without $c_{1}$ enabled after 500 and
1000 steps in repeated transfer of $f_{4}$ on the level-4 mesh. It
can be seen that when not considering the $C^{1}$ discontinuities
(i.e., $c_{1}=0$), the solution was less controlled near $C^{1}$
discontinuities. With our default nonzero $c_{1}$, the remap is more
robust, and the solution approximates the exact solution very well
even after 1000 steps of repeated transfer. 

\subsubsection{\label{subsec:Gibbs-phenomena}Gibbs phenomena in one-step transfer}

\begin{figure}
\subfloat[WLS-ENOR and linear interpolation.]{\begin{centering}
\includegraphics[width=0.48\columnwidth]{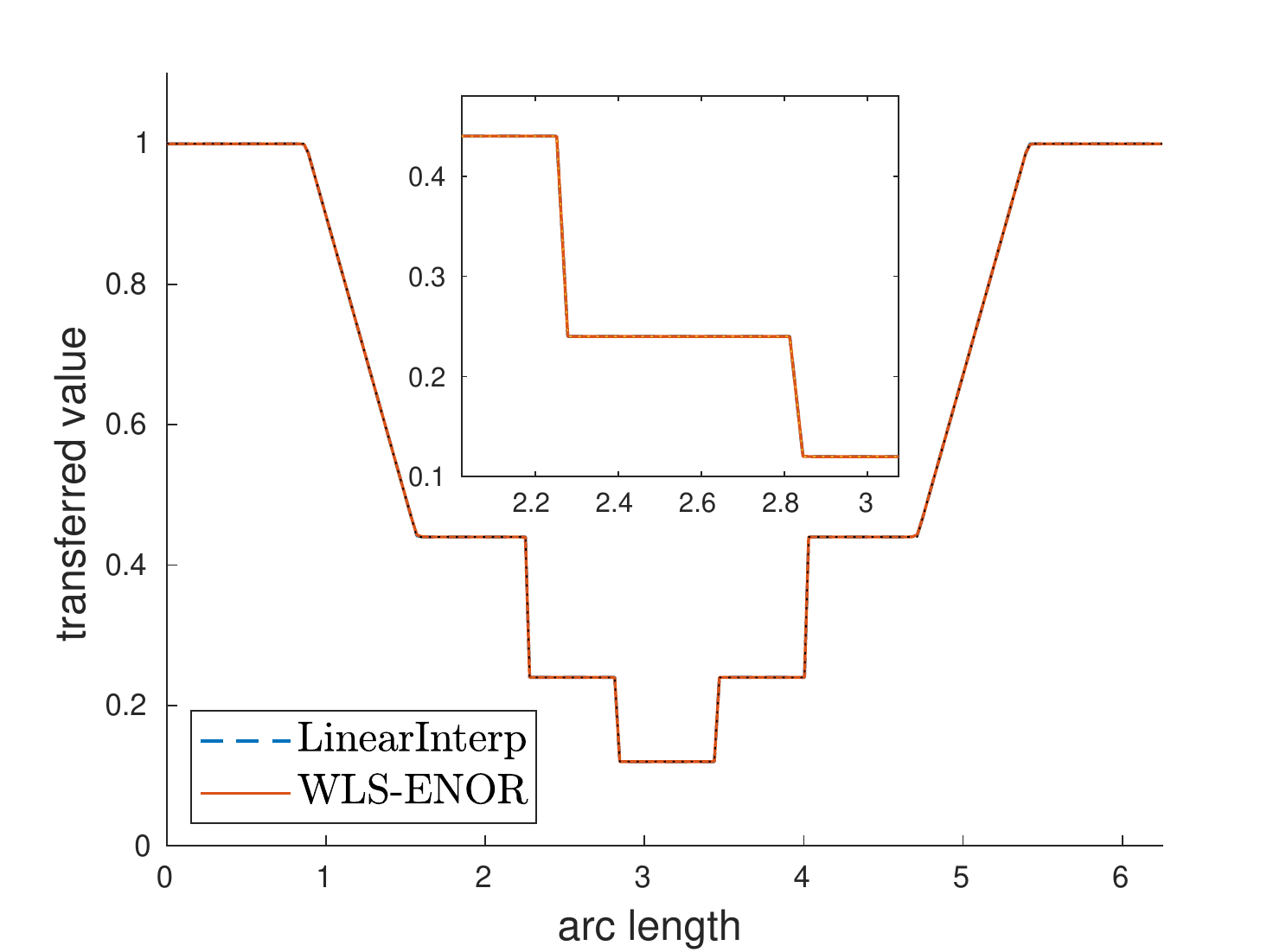}
\par\end{centering}
}\hfill\subfloat[Modified moving least squares (DTK) \citep{slattery2016mesh}.]{\begin{centering}
\includegraphics[width=0.48\columnwidth]{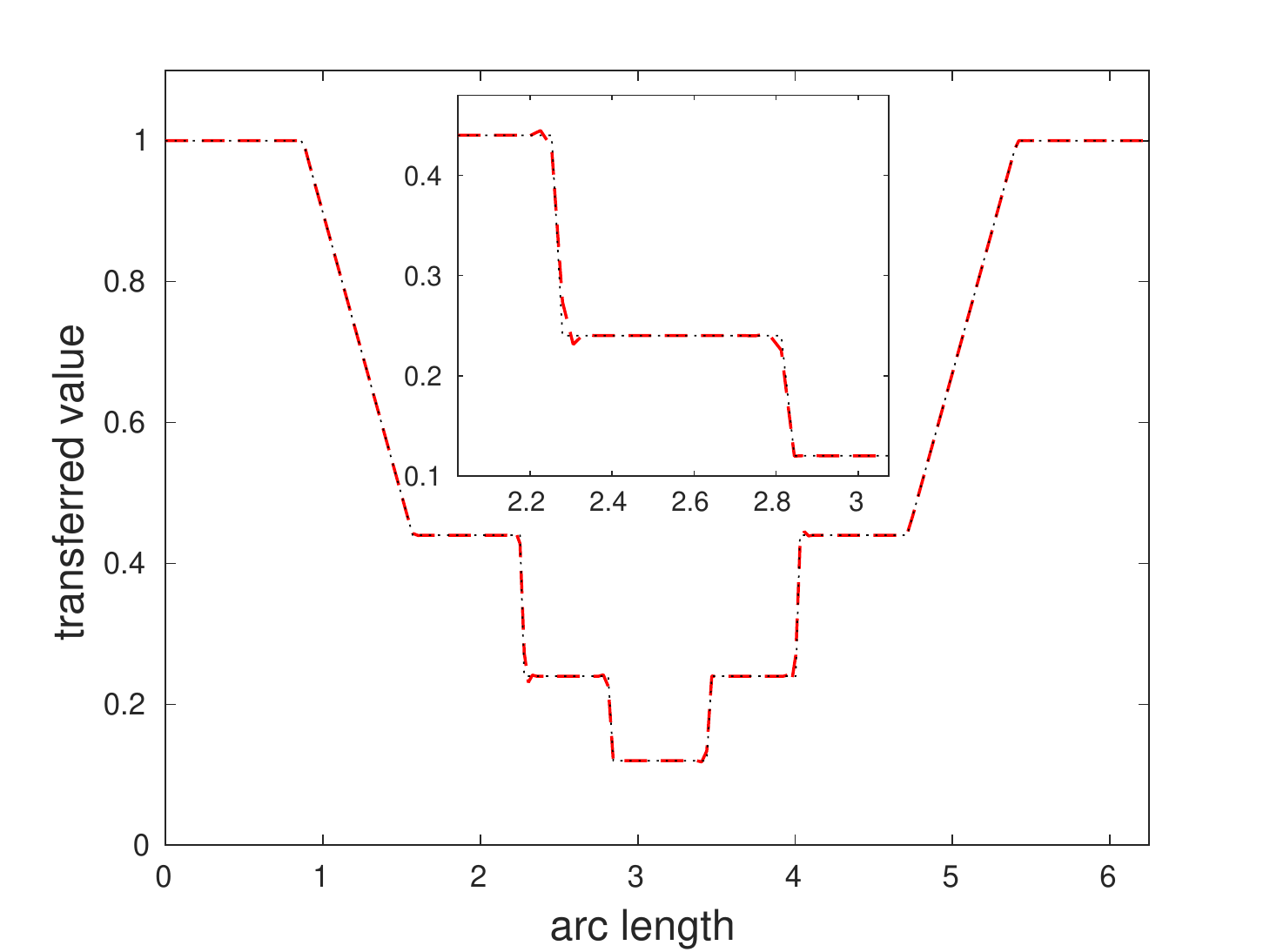}
\par\end{centering}
}

\subfloat[$L^{2}$ projection \citep{jiao2004common}.]{\begin{centering}
\includegraphics[width=0.48\columnwidth]{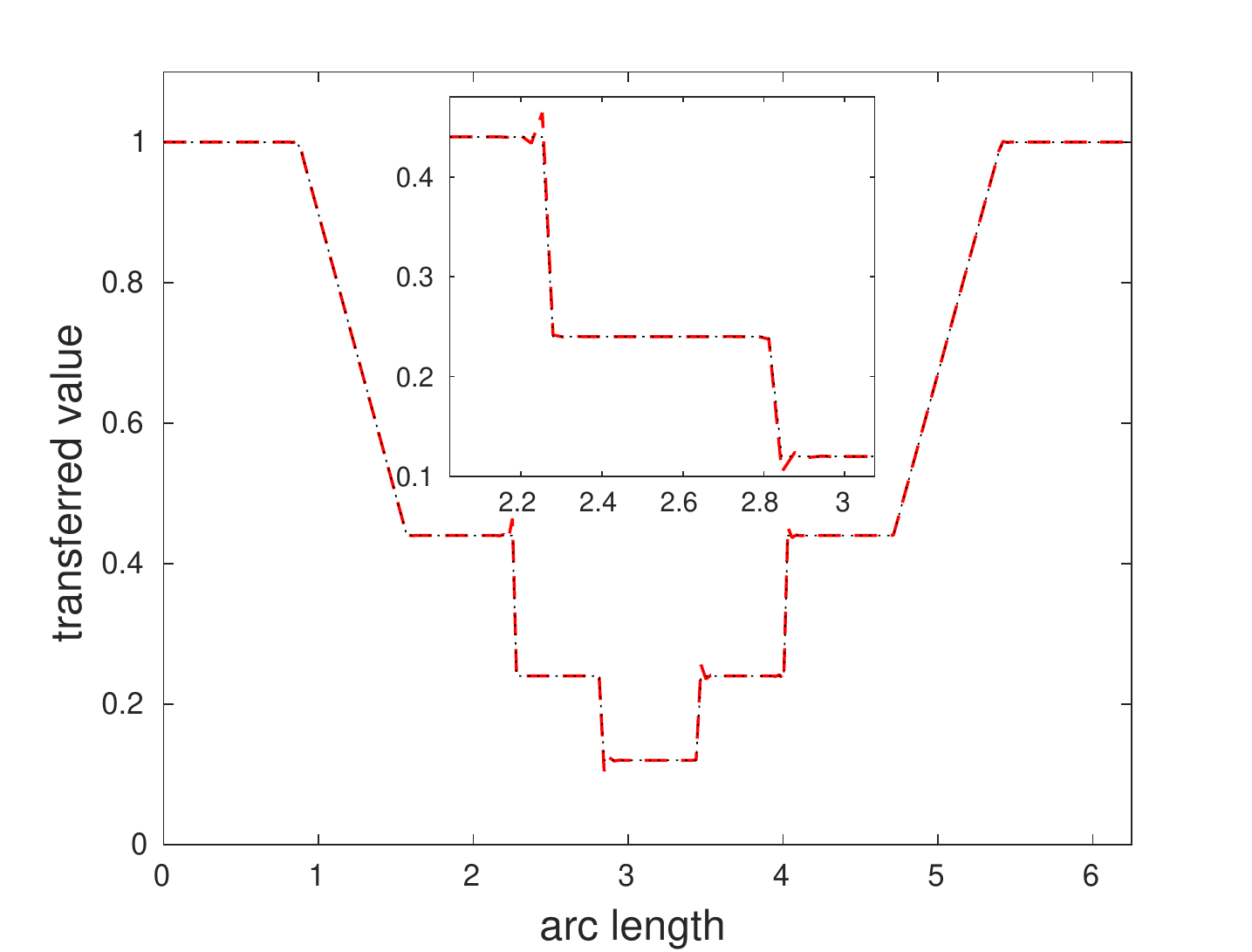}
\par\end{centering}
}\hfill\subfloat[RBF interpolation \citep{beckert2001multivariate}.]{\begin{centering}
\includegraphics[width=0.48\columnwidth]{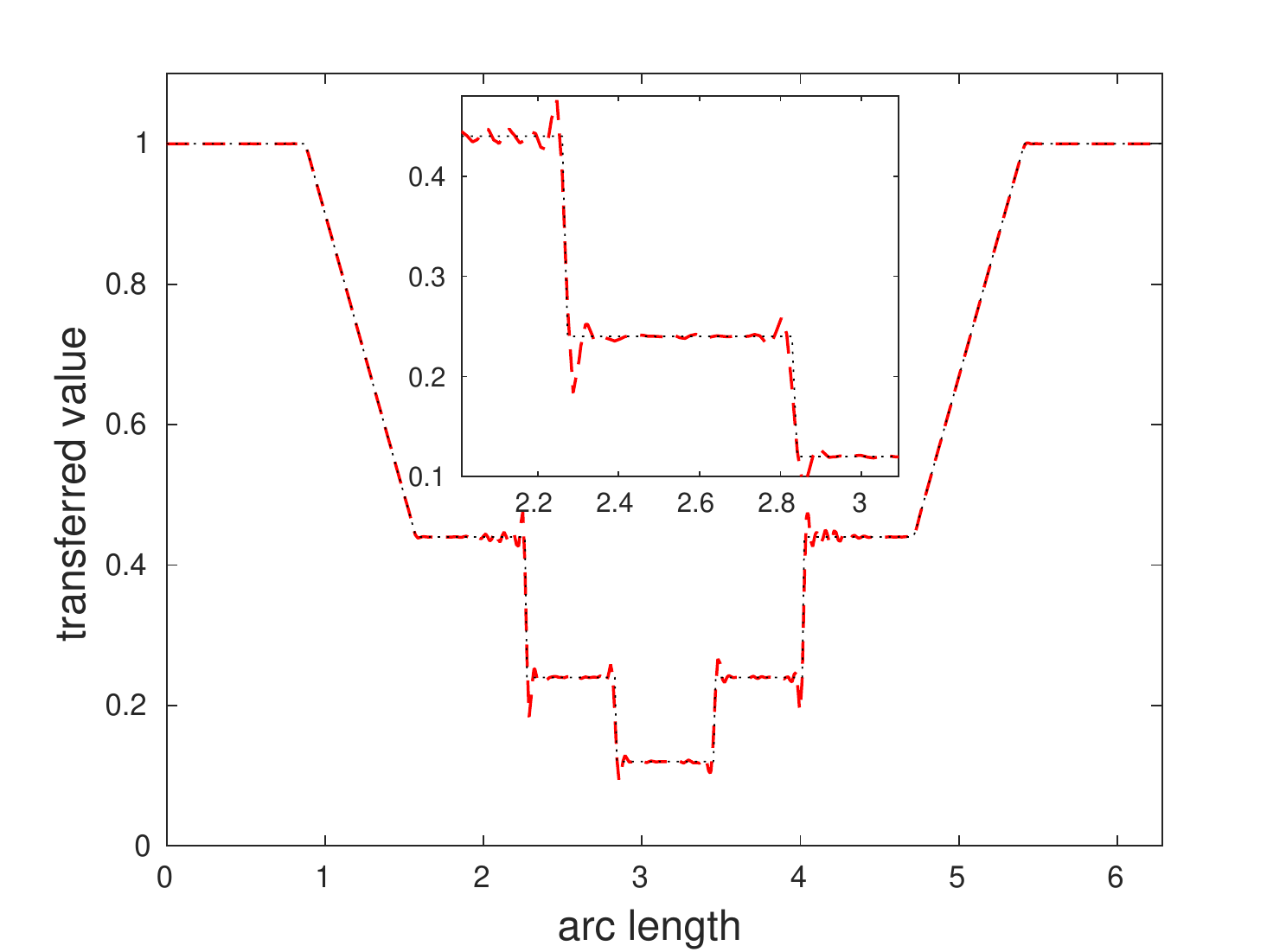}
\par\end{centering}
}

\caption{\label{fig:Comparison-discontinuous}Comparison of transferring $f_{3}$
from level-3 Delaunay to cubed-sphere mesh.}
\end{figure}

We now compare WLS-ENOR with some other methods in transferring discontinuous
function. In particular, we compare WLS-ENO(4,2) with the same methods
used in Section \ref{subsec:Transferring-smooth-functions} by transferring
$f_{3}$ from the level-3 Delaunay to the level-3 cubed-sphere mesh.
Figure~\ref{fig:Comparison-discontinuous} shows the result on the
cubed-sphere mesh along a great circle passing through the poles.
We plot WLS-ENO together with linear interpolation in Figure~\ref{fig:Comparison-discontinuous}(a),
because both are non-oscillatory. As can be seen in Figure~\ref{fig:Comparison-discontinuous}(b--d),
MMLS had visible overshoots and undershoots, although they were less
severe than those of $L^{2}$ projection and RBF interpolation.

\subsubsection{Gibbs phenomena and diffusion in repeated transfer}

\begin{figure}
\subfloat[After 500 steps on level-3 mesh.]{\begin{centering}
\includegraphics[width=0.48\columnwidth]{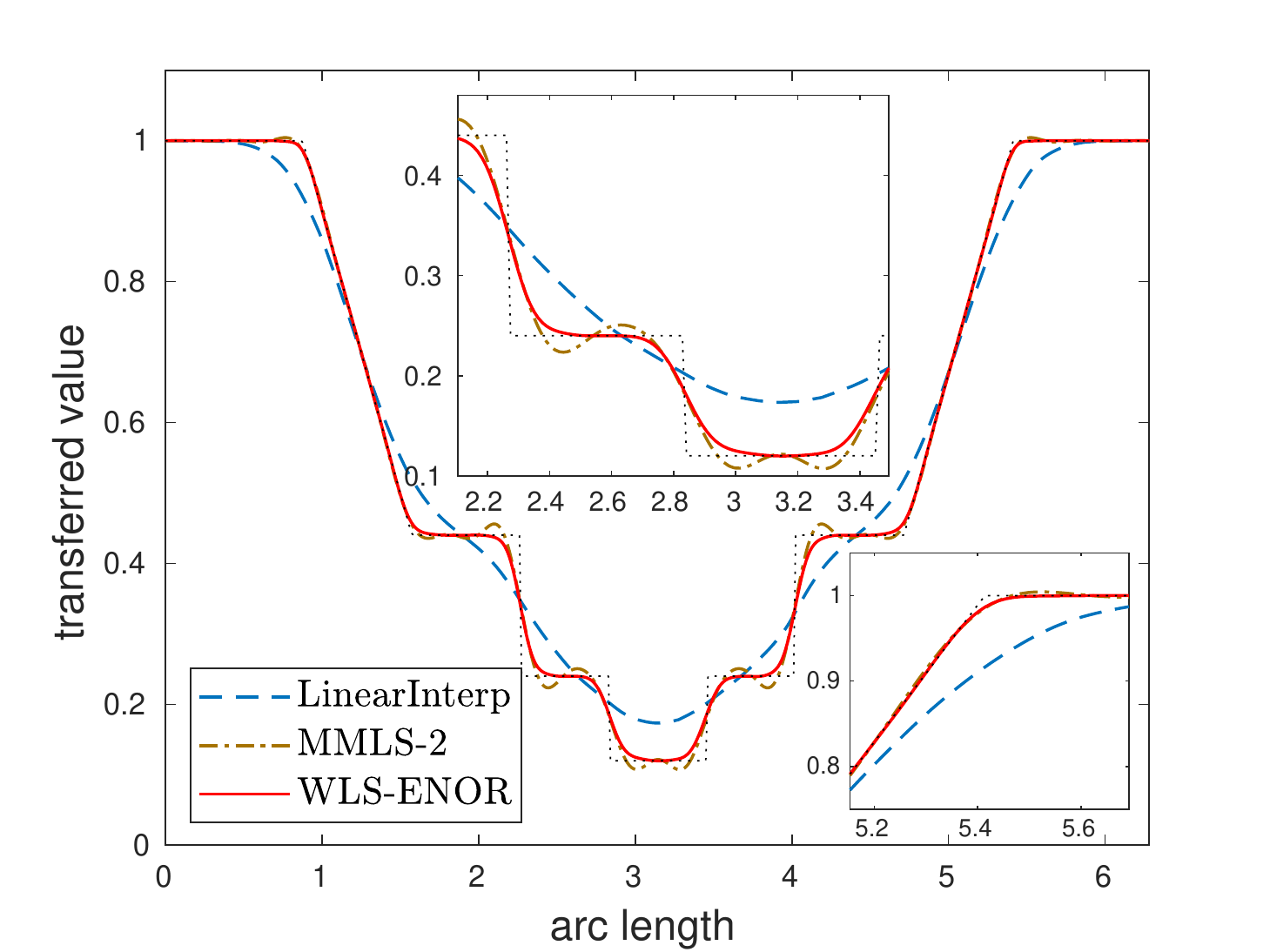}
\par\end{centering}
}\hfill\subfloat[After 1000 steps on level-3 mesh.]{\begin{centering}
\includegraphics[width=0.48\columnwidth]{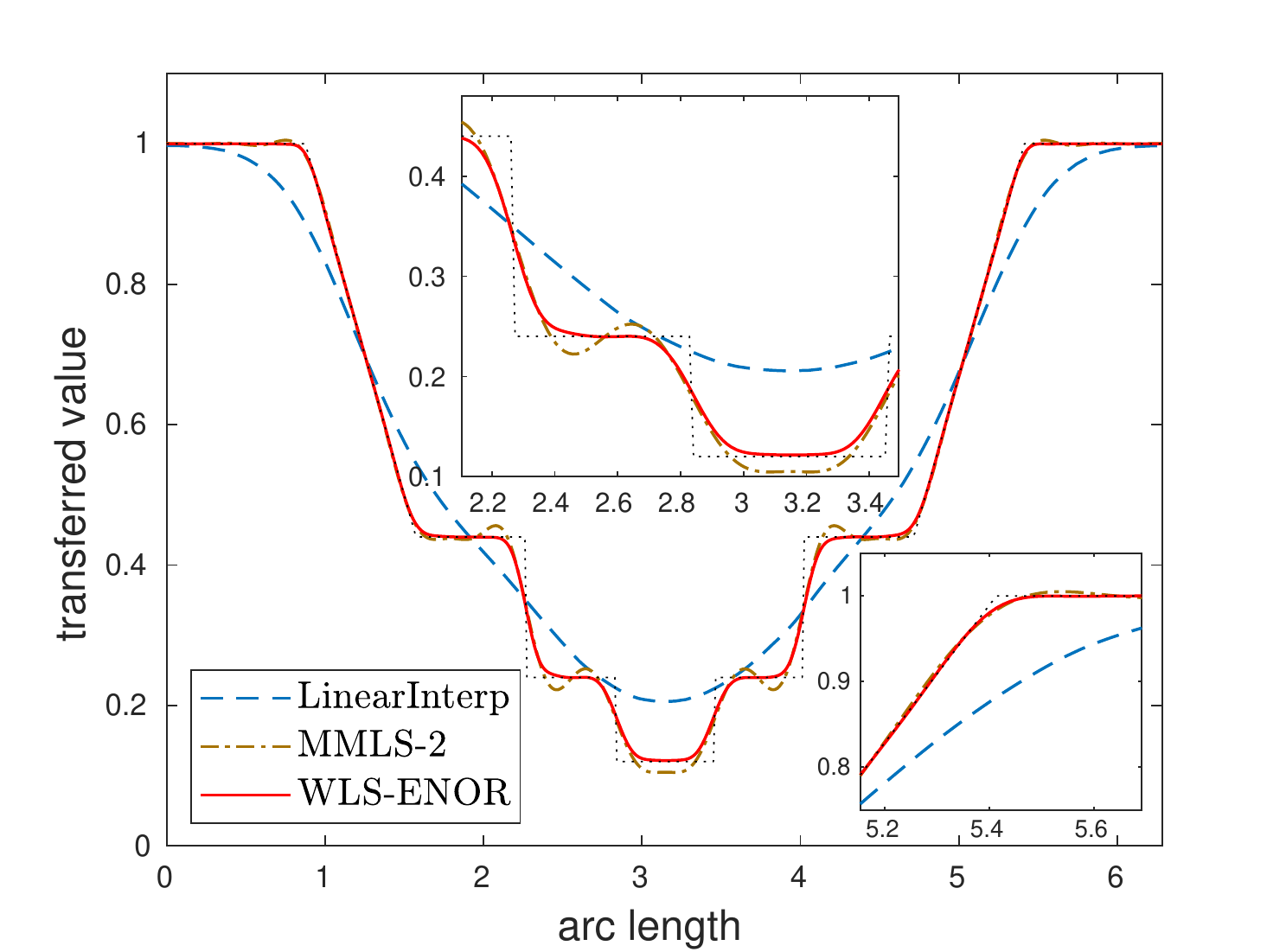}
\par\end{centering}
}

\caption{\label{fig:Comparison-of-diffusion}Comparison of diffusion of WLS-ENOR
versus linear interpolation and MMLS in repeated transfer of $f_{3}$
between level-3 Delaunay and cubed-sphere meshes.}
\end{figure}

We now compare repeated transfer of discontinuous functions with the
least oscillatory methods in the preceding section, namely, WLS-ENOR,
linear interpolation, and MMLS; we refer readers to \citep{slattery2016mesh}
for some comparison of MMLS, $L^{2}$ projection, and RBF interpolation
for repeated transfer of discontinuous functions. Figure~\ref{fig:Comparison-of-diffusion}
shows the results on the level-3 meshes after 500 and 1000 steps of
repeated transfer between the level-3 Delaunay and cubed-sphere meshes
for function $f_{3}$. WLS-ENOR performed the best in terms of both
preserving monotonicity and minimizing diffusion. In contrast, linear
interpolation was excessively diffusive, whereas MMLS was oscillatory
and had a severe undershoot at the global minimum and a noticeable
overshoot at the global maximum.

\begin{figure}
\subfloat[After 1000 iterations on level-2 mesh.]{\begin{centering}
\includegraphics[width=0.48\columnwidth]{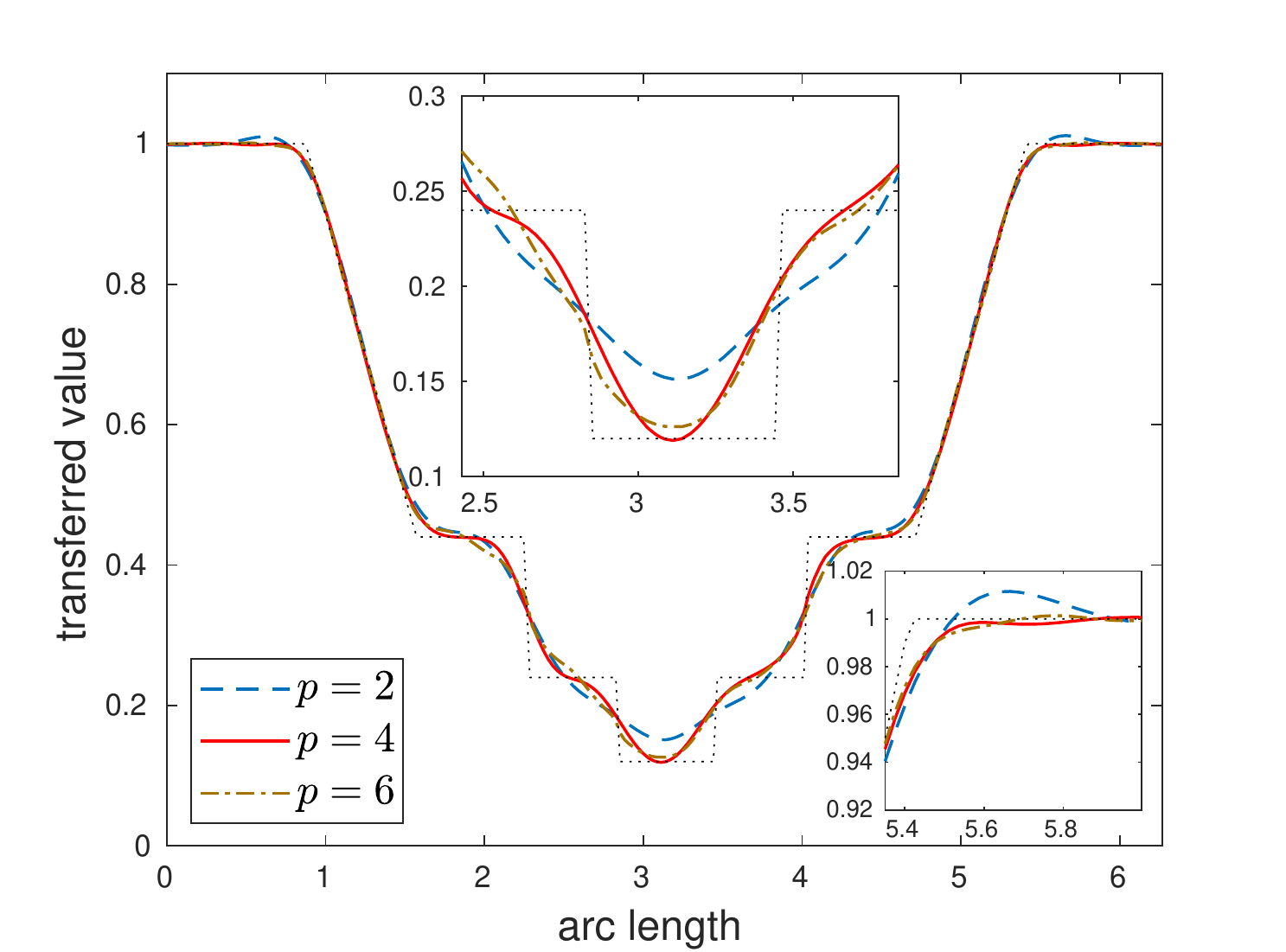}
\par\end{centering}
}\hfill\subfloat[After 1000 iterations on level-4 mesh.]{\begin{centering}
\includegraphics[width=0.48\columnwidth]{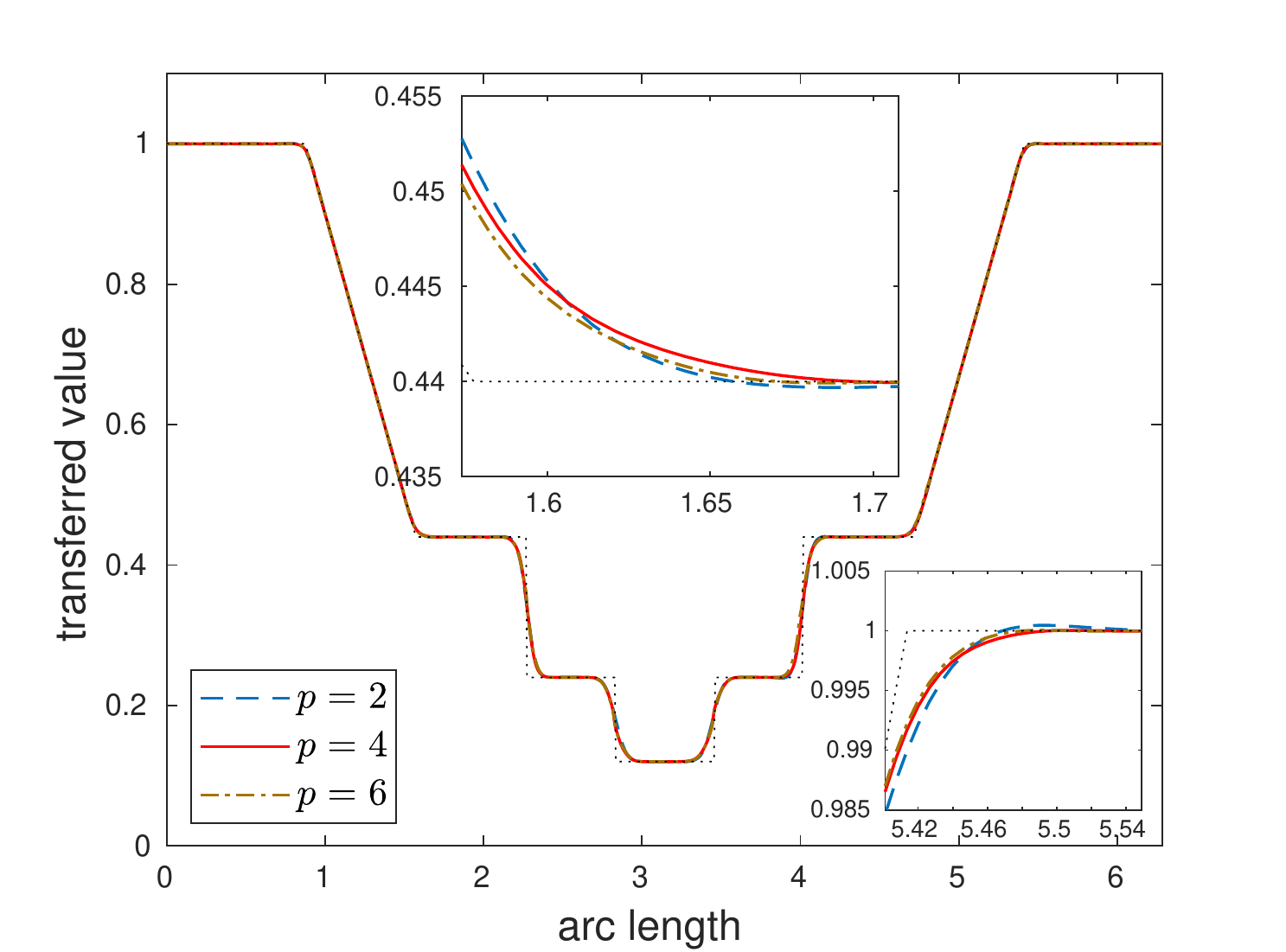}
\par\end{centering}
}

\caption{\label{fig:Comparison-of-wls-enor22}Comparison of WLS-ENOR($p$,2)
with $p=2$, 4, and 6 after 1000 repeated transfers for $f_{3}$ between
Delaunay and cubed-sphere meshes.}
\end{figure}

Finally, we report a comparison of WLS-ENOR using different degree
of polynomials in smooth regions. Figure~\ref{fig:Comparison-of-wls-enor22}
shows the results with quadratic, quartic, and sextic WLS (i.e., WLS-ENOR($p$,2)
with $p=2$, 4, and 6) after 1000 steps of repeated transfer of function
$f_{4}$ on the level-2 and level-4 meshes. On the level-2 mesh, quadratic
WLS was significantly more diffusive than the others near the global
minimum, where two $C^{0}$ discontinuities interact. It also had
some overshoots near $C^{1}$ discontinuities. This is because quadratic
WLS does not have sufficient accuracy on coarser meshes. Quartic and
sextic WLS produced similar results in smooth regions and near $C^{1}$
discontinuities. However, sextic WLS is more expensive, and it was
more diffusive than quartic WLS on the level-2 mesh, This is because
sextic WLS requires about 20 cells between $C^{0}$ discontinuities
to avoid interference from each other, but quartic WLS requires about
16 cells. These results show that WLS-ENOR(4,2) strikes a good balance
between accuracy, efficiency, and robustness.

\section{\label{sec:Conclusions}Conclusions}

In this paper, we proposed a new method, called \emph{WLS-ENO Remap
}or \emph{WLS-ENOR }in short, for transferring data between non-matching
meshes. The method utilizes quartic WLS reconstruction in smooth regions,
and it selectively applies quadratic WLS-ENO reconstruction near discontinuities.
For smooth functions, we proposed a new weighting scheme based on
Buhmann's $C^{3}$ radial function and optimized it to enable superconvergence
of WLS remap for smooth functions. For discontinuous functions, we
proposed a robust technique to detect $C^{0}$ and $C^{1}$ discontinuities
with two novel indicators and a new dual-thresholding strategy. We
also proposed a new weighting scheme for WLS-ENO to take into account
$C^{1}$ discontinuities for robustness. Overall, WLS-ENOR achieves
higher than fifth-order accuracy and highly conservative in smooth
regions, and is non-oscillatory and minimally diffusive near discontinuities.
It compares favorably with other commonly used node-to-node remap
techniques in accuracy, stability, conservation, and monotonicity-preservation
for almost all of our test cases.

In this work, we primarily focused on node-to-node transfer on spheres,
which is most relevant to climate and earth-system modeling. We plan
to extend the work to address cell-averaged data, which are commonly
used in finite volume methods. Future work will also include extending
our method to support more general surfaces, including those with
geometric discontinuities, such as ridges and corners, as well as
those with boundaries.

\section*{Acknowledgments}

This work was supported under the Scientific Discovery through Advanced
Computing (SciDAC) program in the US Department of Energy\textquoteright s
Office of Science, Office of Advanced Scientific Computing Research
through subcontract \#462974 with Los Alamos National Laboratory.

\bibliographystyle{abbrv}
\bibliography{wls-enor}

\begin{thebibliography}{100}

\bibitem{adam2016higher}
A.~Adam, D.~Pavlidis, J.~R. Percival, P.~Salinas, Z.~Xie, F.~Fang, C.~C. Pain,
  A.~Muggeridge, and M.~D. Jackson.
\newblock Higher-order conservative interpolation between control-volume
  meshes: Application to advection and multiphase flow problems with dynamic
  mesh adaptivity.
\newblock {\em J. Comput. Phys.}, 321:512--531, 2016.

\bibitem{archibald2005polynomial}
R.~Archibald, A.~Gelb, and J.~Yoon.
\newblock Polynomial fitting for edge detection in irregularly sampled signals
  and images.
\newblock {\em SIAM J. Numer. Ana.}, 43(1):259--279, 2005.

\bibitem{archibald2008determining}
R.~Archibald, A.~Gelb, and J.~Yoon.
\newblock Determining the locations and discontinuities in the derivatives of
  functions.
\newblock {\em Appl. Numer. Math.}, 58(5):577--592, 2008.

\bibitem{beckert2001multivariate}
A.~Beckert and H.~Wendland.
\newblock Multivariate interpolation for fluid-structure-interaction problems
  using radial basis functions.
\newblock {\em Aerosp. Sci. Technol.}, 5(2):125--134, 2001.

\bibitem{bochev2005constrained}
P.~Bochev and M.~Shashkov.
\newblock Constrained interpolation (remap) of divergence-free fields.
\newblock {\em Comput. Meth. Appl. Mech. Engrg.}, 194(2-5):511--530, 2005.

\bibitem{bozzini2014interpolation}
M.~Bozzini, L.~Lenarduzzi, M.~Rossini, and R.~Schaback.
\newblock Interpolation with variably scaled kernels.
\newblock {\em IMA J. Numer. Ana.}, 35(1):199--219, 2014.

\bibitem{bozzini2013detection}
M.~Bozzini and M.~Rossini.
\newblock The detection and recovery of discontinuity curves from scattered
  data.
\newblock {\em J. Comput. Appl. Math.}, 240:148--162, 2013.

\bibitem{buhmann2001new}
M.~Buhmann.
\newblock A new class of radial basis functions with compact support.
\newblock {\em Math. Comput.}, 70(233):307--318, 2001.

\bibitem{buhmann2000radial}
M.~D. Buhmann.
\newblock Radial basis functions.
\newblock {\em Acta Numer.}, 9:1--38, 2000.

\bibitem{buhmann2003radial}
M.~D. Buhmann.
\newblock {\em Radial Basis Functions: Theory and Implementations}, volume~12.
\newblock Cambridge University Press, 2003.

\bibitem{bungartz2016precice}
H.-J. Bungartz, F.~Lindner, B.~Gatzhammer, M.~Mehl, K.~Scheufele, A.~Shukaev,
  and B.~Uekermann.
\newblock {preCICE} -- a fully parallel library for multi-physics surface
  coupling.
\newblock {\em Comput. Fluids}, 141:250--258, 2016.

\bibitem{canny1987computational}
J.~Canny.
\newblock A computational approach to edge detection.
\newblock In {\em Readings in Computer Vision}, pages 184--203. Elsevier, 1987.

\bibitem{cates2007detecting}
D.~Cates and A.~Gelb.
\newblock Detecting derivative discontinuity locations in piecewise continuous
  functions from fourier spectral data.
\newblock {\em Numer. Algo.}, 46(1):59--84, 2007.

\bibitem{chesshire1994scheme}
G.~Chesshire and W.~D. Henshaw.
\newblock A scheme for conservative interpolation on overlapping grids.
\newblock {\em SIAM J. Sci. Comput.}, 15(4):819--845, 1994.

\bibitem{cockburn2000development}
B.~Cockburn, G.~E. Karniadakis, and C.-W. Shu.
\newblock {\em The Development of Discontinuous Galerkin Methods}.
\newblock Springer, 2000.

\bibitem{colella1984piecewise}
P.~Colella and P.~R. Woodward.
\newblock The piecewise parabolic method ({PPM}) for gas-dynamical simulations.
\newblock {\em J. Comput. Phys.}, 54(1):174--201, 1984.

\bibitem{craig2017development}
A.~Craig, S.~Valcke, and L.~Coquart.
\newblock Development and performance of a new version of the {OASIS} coupler,
  {OASIS3-MCT}\_3. 0.
\newblock {\em Geosci. Model Dev.}, 10(9):3297--3308, 2017.

\bibitem{de2008comparison}
A.~de~Boer, A.~H. van Zuijlen, and H.~Bijl.
\newblock Comparison of conservative and consistent approaches for the coupling
  of non-matching meshes.
\newblock {\em Comput. Methods Appl. Mech. Eng.}, 197(49-50):4284--4297, 2008.

\bibitem{de1990quasiinterpolants}
C.~de~Boor.
\newblock Quasiinterpolants and approximation power of multivariate splines.
\newblock In {\em Computation of curves and surfaces}, pages 313--345.
  Springer, 1990.

\bibitem{de2002vertical}
M.~C.~L. de~Silanes, M.~C. Parra, and J.~J. Torrens.
\newblock Vertical and oblique fault detection in explicit surfaces.
\newblock {\em J. Comput. Appl. Math.}, 140(1-2):559--585, 2002.

\bibitem{DREJTAHF2014}
V.~Dyedov, N.~Ray, D.~Einstein, X.~Jiao, and T.~Tautges.
\newblock {AHF}: Array-based half-facet data structure for mixed-dimensional
  and non-manifold meshes.
\newblock In J.~Sarrate and M.~Staten, editors, {\em Proceedings of the 22nd
  International Meshing Roundtable}, pages 445--464. Springer International
  Publishing, 2014.

\bibitem{farrell2011conservative}
P.~Farrell and J.~Maddison.
\newblock Conservative interpolation between volume meshes by local {G}alerkin
  projection.
\newblock {\em Comput. Meth. Appl. Mech. Engrg.}, 200(1-4):89--100, 2011.

\bibitem{farrell2009conservative}
P.~Farrell, M.~Piggott, C.~Pain, G.~Gorman, and C.~Wilson.
\newblock Conservative interpolation between unstructured meshes via supermesh
  construction.
\newblock {\em Comput. Meth. Appl. Mech. Engrg.}, 198(33-36):2632--2642, 2009.

\bibitem{fleishman2005robust}
S.~Fleishman, D.~Cohen-Or, and C.~T. Silva.
\newblock Robust moving least-squares fitting with sharp features.
\newblock {\em ACM Trans. Graph.}, 24(3):544--552, 2005.

\bibitem{flyer2007transport}
N.~Flyer and G.~B. Wright.
\newblock Transport schemes on a sphere using radial basis functions.
\newblock {\em J. Comput. Phys.}, 226(1):1059--1084, 2007.

\bibitem{fornberg2007gibbs}
B.~Fornberg and N.~Flyer.
\newblock The {Gibbs} phenomenon for radial basis functions.
\newblock {\em The Gibbs Phenomenon in Various Representations and
  Applications}, pages 201--224, 2007.

\bibitem{foster1991gibbs}
J.~Foster and F.~Richards.
\newblock The {G}ibbs phenomenon for piecewise-linear approximation.
\newblock {\em Amer. Math. Monthly}, 98(1):47--49, 1991.

\bibitem{gander2013algorithm}
M.~J. Gander and C.~Japhet.
\newblock Algorithm 932: {PANG}: software for nonmatching grid projections in
  {2D} and {3D} with linear complexity.
\newblock {\em ACM Trans. Math. Softw.}, 40(1):6, 2013.

\bibitem{gelb1999detection}
A.~Gelb and E.~Tadmor.
\newblock Detection of edges in spectral data.
\newblock {\em Appl. Comput. Harmon. Anal.}, 7(1):101--135, 1999.

\bibitem{gelb2000detection}
A.~Gelb and E.~Tadmor.
\newblock Detection of edges in spectral data ii. nonlinear enhancement.
\newblock {\em SIAM J. Numer. Anal.}, 38(4):1389--1408, 2000.

\bibitem{gelb2006adaptive}
A.~Gelb and E.~Tadmor.
\newblock Adaptive edge detectors for piecewise smooth data based on the minmod
  limiter.
\newblock {\em J. Sci. Comput.}, 28(2-3):279--306, 2006.

\bibitem{gibbs1899}
J.~W. Gibbs.
\newblock Letter to the editor.
\newblock {\em Nature}, page 200 and 606, 1899.

\bibitem{godunov1959difference}
S.~K. Godunov.
\newblock A difference method for numerical calculation of discontinuous
  solutions of the equations of hydrodynamics.
\newblock {\em Mat. Sb.}, 89(3):271--306, 1959.

\bibitem{Golub13MC}
G.~H. Golub and C.~F. {Van Loan}.
\newblock {\em Matrix Computations}.
\newblock Johns Hopkins, 4th edition, 2013.

\bibitem{gottlieb1997gibbs}
D.~Gottlieb and C.-W. Shu.
\newblock On the {G}ibbs phenomenon and its resolution.
\newblock {\em SIAM Rev.}, 39(4):644--668, 1997.

\bibitem{gottlieb1998total}
S.~Gottlieb and C.-W. Shu.
\newblock Total variation diminishing {R}unge-{K}utta schemes.
\newblock {\em Math. Comput.}, 67(221):73--85, 1998.

\bibitem{grandy1999conservative}
J.~Grandy.
\newblock Conservative remapping and region overlays by intersecting arbitrary
  polyhedra.
\newblock {\em J. Comput. Phys.}, 148(2):433--466, 1999.

\bibitem{harder1972interpolation}
R.~L. Harder and R.~N. Desmarais.
\newblock Interpolation using surface splines.
\newblock {\em J. Aircraft}, 9(2):189--191, 1972.

\bibitem{harten1989eno}
A.~Harten.
\newblock {ENO} schemes with subcell resolution.
\newblock {\em J. Comput. Phys.}, 83(1):148--184, 1989.

\bibitem{hewitt1979gibbs}
E.~Hewitt and R.~E. Hewitt.
\newblock The {G}ibbs-{W}ilbraham phenomenon: an episode in {F}ourier analysis.
\newblock {\em Arch. Hist. Exact Sci.}, 21(2):129--160, 1979.

\bibitem{Higham87SCN}
N.~J. Higham.
\newblock A survey of condition number estimation for triangular matrices.
\newblock {\em SIAM Rev.}, 29:575--596, 1987.

\bibitem{hill2004architecture}
C.~Hill, C.~DeLuca, M.~Suarez, A.~Da~Silva, et~al.
\newblock The architecture of the earth system modeling framework.
\newblock {\em Comput. Sci. Engrg.}, 6(1):18, 2004.

\bibitem{hu1999weighted}
C.~Hu and C.-W. Shu.
\newblock Weighted essentially non-oscillatory schemes on triangular meshes.
\newblock {\em J. Comput. Phys.}, 150(1):97--127, 1999.

\bibitem{humpherys2017foundations}
J.~Humpherys, T.~J. Jarvis, and E.~J. Evans.
\newblock {\em Foundations of Applied Mathematics, Volume I: Mathematical
  Analysis}.
\newblock SIAM, 2017.

\bibitem{hurrell2013community}
J.~W. Hurrell, M.~M. Holland, P.~R. Gent, S.~Ghan, J.~E. Kay, P.~J. Kushner,
  J.-F. Lamarque, W.~G. Large, D.~Lawrence, K.~Lindsay, et~al.
\newblock The community earth system model: a framework for collaborative
  research.
\newblock {\em B. Am. Meterol. Soc.}, 94(9):1339--1360, 2013.

\bibitem{jerri2013gibbs}
A.~J. Jerri.
\newblock {\em The {G}ibbs Phenomenon in {F}ourier Analysis, Splines and
  Wavelet Approximations}, volume 446.
\newblock Springer, 2013.

\bibitem{jiang1996efficient}
G.-S. Jiang and C.-W. Shu.
\newblock Efficient implementation of weighted {ENO} schemes.
\newblock {\em J. Comput. Phys.}, 126(1):202--228, 1996.

\bibitem{jiao1999mesh}
X.~Jiao, H.~Edelsbrunner, and M.~T. Heath.
\newblock Mesh association: formulation and algorithms.
\newblock In {\em International Meshing Roundtable}, pages 75--82. Citeseer,
  1999.

\bibitem{jiao2004common}
X.~Jiao and M.~T. Heath.
\newblock Common-refinement-based data transfer between non-matching meshes in
  multiphysics simulations.
\newblock {\em Int. J. Numer. Meth. Engrg.}, 61(14):2402--2427, 2004.

\bibitem{jiao2004overlayingI}
X.~Jiao and M.~T. Heath.
\newblock Overlaying surface meshes, part {I}: Algorithms.
\newblock {\em Int. J. Comput. Geom. Appl.}, 14(06):379--402, 2004.

\bibitem{jiao2004overlayingII}
X.~Jiao and M.~T. Heath.
\newblock Overlaying surface meshes, part {II}: Topology preservation and
  feature matching.
\newblock {\em Int. J. Comput. Geom. Appl.}, 14(06):403--419, 2004.

\bibitem{Jiao2011RHO}
X.~Jiao and D.~Wang.
\newblock Reconstructing high-order surfaces for meshing.
\newblock {\em Engrg. Comput.}, 28:361--373, 2012.

\bibitem{JZ08CCF}
X.~Jiao and H.~Zha.
\newblock Consistent computation of first- and second-order differential
  quantities for surface meshes.
\newblock In {\em ACM Solid and Physical Modeling Symposium}, pages 159--170.
  ACM, 2008.

\bibitem{joldes2015modified}
G.~R. Joldes, H.~A. Chowdhury, A.~Wittek, B.~Doyle, and K.~Miller.
\newblock Modified moving least squares with polynomial bases for scattered
  data approximation.
\newblock {\em Appl. Math. Comput.}, 266:893--902, 2015.

\bibitem{jones1999first}
P.~W. Jones.
\newblock First- and second-order conservative remapping schemes for grids in
  spherical coordinates.
\newblock {\em Mon. Weather Rev.}, 127(9):2204--2210, 1999.

\bibitem{joppich2006mpcci}
W.~Joppich and M.~K{\"u}rschner.
\newblock {MpCCI}--a tool for the simulation of coupled applications.
\newblock {\em Concurr. Comp.-Pract. E.}, 18(2):183--192, 2006.

\bibitem{ju2011voronoi}
L.~Ju, T.~Ringler, and M.~Gunzburger.
\newblock Voronoi tessellations and their application to climate and global
  modeling.
\newblock In {\em Numerical Techniques for Global Atmospheric Models}, pages
  313--342. Springer, 2011.

\bibitem{jung2007note}
J.-H. Jung.
\newblock A note on the {G}ibbs phenomenon with multiquadric radial basis
  functions.
\newblock {\em Appl. Numer. Math.}, 57(2):213--229, 2007.

\bibitem{jung2011iterative}
J.-H. Jung, S.~Gottlieb, and S.~O. Kim.
\newblock Iterative adaptive {RBF} methods for detection of edges in
  two-dimensional functions.
\newblock {\em Appl. Numer. Math.}, 61(1):77--91, 2011.

\bibitem{jung2004generalization}
J.-H. Jung and B.~D. Shizgal.
\newblock Generalization of the inverse polynomial reconstruction method in the
  resolution of the {G}ibbs phenomenon.
\newblock {\em J. Comput. Appl. Math.}, 172(1):131--151, 2004.

\bibitem{karniadakis2013spectral}
G.~Karniadakis and S.~Sherwin.
\newblock {\em Spectral/hp Element Methods for Computational Fluid Dynamics}.
\newblock Oxford University Press, 2013.

\bibitem{kelly1996gibbs}
S.~E. Kelly.
\newblock Gibbs phenomenon for wavelets.
\newblock {\em Appl. Comput. Harmon. Anal.}, 3(1):72--81, 1996.

\bibitem{keyes2013multiphysics}
D.~E. Keyes, L.~C. McInnes, C.~Woodward, W.~Gropp, E.~Myra, M.~Pernice,
  J.~Bell, J.~Brown, A.~Clo, J.~Connors, et~al.
\newblock Multiphysics simulations: challenges and opportunities.
\newblock {\em Int. J. High Perform. Comput. Appl.}, 27(1):4--83, 2013.

\bibitem{lancaster1981surfaces}
P.~Lancaster and K.~Salkauskas.
\newblock Surfaces generated by moving least squares methods.
\newblock {\em Math. Comput.}, 37(155):141--158, 1981.

\bibitem{Lancaster1981}
P.~Lancaster and K.~Salkauskas.
\newblock {Surfaces Generated by Moving Least Squares Methods}.
\newblock {\em Mathematics of Computation}, 37(155):141 -- 158, 1981.

\bibitem{larson2005model}
J.~Larson, R.~Jacob, and E.~Ong.
\newblock The model coupling toolkit: a new {F}ortran90 toolkit for building
  multiphysics parallel coupled models.
\newblock {\em Int. J. High Perform. Comput. Appl.}, 19(3):277--292, 2005.

\bibitem{lauritzen2008monotone}
P.~H. Lauritzen and R.~D. Nair.
\newblock Monotone and conservative cascade remapping between spherical grids
  (cars): Regular latitude--longitude and cubed-sphere grids.
\newblock {\em Mon. Weather Rev.}, 136(4):1416--1432, 2008.

\bibitem{lauritzen2018ncar}
P.~H. Lauritzen, R.~D. Nair, A.~Herrington, P.~Callaghan, S.~Goldhaber,
  J.~Dennis, J.~Bacmeister, B.~Eaton, C.~Zarzycki, M.~A. Taylor, et~al.
\newblock {NCAR} release of {CAM-SE} in {CESM2.0}: A reformulation of the
  spectral element dynamical core in dry-mass vertical coordinates with
  comprehensive treatment of condensates and energy.
\newblock {\em J. Adv. Model Earth Sy.}, 10(7):1537--1570, 2018.

\bibitem{leveque1992numerical}
R.~J. LeVeque.
\newblock {\em Numerical Methods for Conservation Laws}, volume 132.
\newblock Springer, 1992.

\bibitem{lipman2009approximating}
Y.~Lipman and D.~Levin.
\newblock Approximating piecewise-smooth functions.
\newblock {\em IMA J. Numer. Ana.}, 30(4):1159--1183, 2009.

\bibitem{liu2016wls}
H.~Liu and X.~Jiao.
\newblock {WLS-ENO}: Weighted-least-squares based essentially non-oscillatory
  schemes for finite volume methods on unstructured meshes.
\newblock {\em J. Comput. Phys.}, 314:749--773, 2016.

\bibitem{liu1994weighted}
X.-D. Liu, S.~Osher, and T.~Chan.
\newblock Weighted essentially non-oscillatory schemes.
\newblock {\em J. Comput. Phys.}, 115(1):200--212, 1994.

\bibitem{liuzhang2013robust}
Y.~Liu and Y.-T. Zhang.
\newblock A robust reconstruction for unstructured {WENO} schemes.
\newblock {\em J. Sci. Comput.}, 54(2,3):603--621, 2013.

\bibitem{luo2013reconstructed}
H.~Luo, Y.~Xia, S.~Spiegel, R.~Nourgaliev, and Z.~Jiang.
\newblock A reconstructed discontinuous {G}alerkin method based on a
  hierarchical {WENO} reconstruction for compressible flows on tetrahedral
  grids.
\newblock {\em J. Comput. Phys.}, 236:477--492, 2013.

\bibitem{mahadevan2015sigma}
V.~Mahadevan, I.~R. Grindeanu, N.~Ray, R.~Jain, and D.~Wu.
\newblock {SIGMA} release v1. 2-capabilities, enhancements and fixes.
\newblock Technical report, Argonne National Laboratory (ANL), Argonne, IL
  (United States), 2015.

\bibitem{margolin2003second}
L.~Margolin and M.~Shashkov.
\newblock Second-order sign-preserving conservative interpolation (remapping)
  on general grids.
\newblock {\em J. Comput. Phys.}, 184(1):266--298, 2003.

\bibitem{michelson1898letter}
A.~A. Michelson.
\newblock Letter to the editor.
\newblock {\em Nature}, 58:544--545, 1898.

\bibitem{petersen2012hypothesis}
A.~Petersen, A.~Gelb, and R.~Eubank.
\newblock Hypothesis testing for {F}ourier based edge detection methods.
\newblock {\em J. Sci. Comput.}, 51(3):608--630, 2012.

\bibitem{putman2007finite}
W.~M. Putman and S.-J. Lin.
\newblock Finite-volume transport on various cubed-sphere grids.
\newblock {\em J. Comput. Phys.}, 227(1):55--78, 2007.

\bibitem{qiu2005comparison}
J.~Qiu and C.-W. Shu.
\newblock A comparison of troubled-cell indicators for {R}unge--{K}utta
  discontinuous {G}alerkin methods using weighted essentially nonoscillatory
  limiters.
\newblock {\em SIAM J. Sci. Comput.}, 27(3):995--1013, 2005.

\bibitem{RayWanJia12}
N.~Ray, D.~Wang, X.~Jiao, and J.~Glimm.
\newblock High-order numerical integration over discrete surfaces.
\newblock {\em SIAM J. Numer. Ana.}, 50:3061--3083, 2012.

\bibitem{rendall2008unified}
T.~Rendall and C.~Allen.
\newblock Unified fluid--structure interpolation and mesh motion using radial
  basis functions.
\newblock {\em Int. J. Numer. Meth. Engrg.}, 74(10):1519--1559, 2008.

\bibitem{richards1991gibbs}
F.~Richards.
\newblock A {G}ibbs phenomenon for spline functions.
\newblock {\em J. Approx. Theory}, 66(3):334--351, 1991.

\bibitem{romani2019edge}
L.~Romani, M.~Rossini, and D.~Schenone.
\newblock Edge detection methods based on {RBF} interpolation.
\newblock {\em J. Comput. Appl. Math.}, 349:532--547, 2019.

\bibitem{rossini2018interpolating}
M.~Rossini.
\newblock Interpolating functions with gradient discontinuities via variably
  scaled kernels.
\newblock {\em Dolomites Research Notes on Approximation}, 11(2), 2018.

\bibitem{sadourny1972conservative}
R.~Sadourny.
\newblock Conservative finite-difference approximations of the primitive
  equations on quasi-uniform spherical grids.
\newblock {\em Mon. Weather Rev.}, 100(2):136--144, 1972.

\bibitem{saxena2009high}
R.~Saxena, A.~Gelb, and H.~Mittelmann.
\newblock A high order method for determining the edges in the gradient of a
  function.
\newblock {\em Comm. Comput. Phys.}, 5(2-4):694--711, 2009.

\bibitem{shi2002technique}
J.~Shi, C.~Hu, and C.-W. Shu.
\newblock A technique of treating negative weights in {WENO} schemes.
\newblock {\em J. Comput. Phys.}, 175(1):108--127, 2002.

\bibitem{shim1996gibbs}
H.-T. Shim and H.~Volkmer.
\newblock On the {G}ibbs phenomenon for wavelet expansions.
\newblock {\em J. Approx. Theory}, 84(1):74--95, 1996.

\bibitem{shizgal2003towards}
B.~D. Shizgal and J.-H. Jung.
\newblock Towards the resolution of the {G}ibbs phenomena.
\newblock {\em J. Comput. Appl. Math.}, 161(1):41--65, 2003.

\bibitem{shu1999high}
C.-W. Shu.
\newblock High order {ENO} and {WENO} schemes for computational fluid dynamics.
\newblock In {\em High-Order Methods for Computational Physics}, pages
  439--582. Springer, 1999.

\bibitem{shu2009high}
C.-W. Shu.
\newblock High order weighted essentially nonoscillatory schemes for convection
  dominated problems.
\newblock {\em SIAM Rev.}, 51(1):82--126, 2009.

\bibitem{slattery2013data}
S.~Slattery, P.~Wilson, and R.~Pawlowski.
\newblock The {Data Transfer Kit}: A geometric rendezvous-based tool for
  multiphysics data transfer.
\newblock In {\em International Conference on Mathematics \& Computational
  Methods Applied to Nuclear Science \& Engineering (M\&C 2013)}, pages 5--9,
  2013.

\bibitem{slattery2016mesh}
S.~R. Slattery.
\newblock Mesh-free data transfer algorithms for partitioned multiphysics
  problems: Conservation, accuracy, and parallelism.
\newblock {\em J. Comput. Phys.}, 307:164--188, 2016.

\bibitem{tadmor2007filters}
E.~Tadmor.
\newblock Filters, mollifiers and the computation of the {G}ibbs phenomenon.
\newblock {\em Acta Numer.}, 16:305--378, 2007.

\bibitem{tautges2004moab}
T.~J. Tautges, C.~Ernst, C.~Stimpson, R.~J. Meyers, and K.~Merkley.
\newblock {MOAB}: a mesh-oriented database.
\newblock Technical report, Sandia National Laboratories, 2004.

\bibitem{terai2018atmospheric}
C.~R. Terai, P.~M. Caldwell, S.~A. Klein, Q.~Tang, and M.~L. Branstetter.
\newblock The atmospheric hydrologic cycle in the {ACME} v0. 3 model.
\newblock {\em Clim. Dyn.}, 50(9-10):3251--3279, 2018.

\bibitem{ullrich2016arbitrary}
P.~A. Ullrich, D.~Devendran, and H.~Johansen.
\newblock Arbitrary-order conservative and consistent remapping and a theory of
  linear maps: Part ii.
\newblock {\em Mon. Weather Rev.}, 144(4):1529--1549, 2016.

\bibitem{ullrich2015arbitrary}
P.~A. Ullrich and M.~A. Taylor.
\newblock Arbitrary-order conservative and consistent remapping and a theory of
  linear maps: Part i.
\newblock {\em Mon. Weather Rev.}, 143(6):2419--2440, 2015.

\bibitem{van1979towards}
B.~Van~Leer.
\newblock Towards the ultimate conservative difference scheme. {V}. a
  second-order sequel to {G}odunov's method.
\newblock {\em J. Comput. Phys.}, 32(1):101--136, 1979.

\bibitem{viswanathan2008detection}
A.~Viswanathan, D.~Cochran, A.~Gelb, and D.~M. Cates.
\newblock Detection of signal discontinuities from noisy {F}ourier data.
\newblock In {\em 2008 42nd Asilomar Conference on Signals, Systems and
  Computers}, pages 1705--1708. IEEE, 2008.

\bibitem{wendland1995piecewise}
H.~Wendland.
\newblock Piecewise polynomial, positive definite and compactly supported
  radial functions of minimal degree.
\newblock {\em Adv. Comput. Math.}, 4(1):389--396, 1995.

\bibitem{wendland2001local}
H.~Wendland.
\newblock Local polynomial reproduction and moving least squares approximation.
\newblock {\em IMA J. Numer. Ana.}, 21(1):285--300, 2001.

\bibitem{wiki:sphereical_harmonics}
{Wikipedia contributors}.
\newblock Table of spherical harmonics --- {Wikipedia}{,} the free
  encyclopedia, 2019.
\newblock [Online; accessed 24-July-2019].

\bibitem{wilbraham1848certain}
H.~Wilbraham.
\newblock On a certain periodic function.
\newblock {\em The Cambridge and Dublin Mathematical Journal}, 3:198--201,
  1848.

\bibitem{wu1995compactly}
Z.~Wu.
\newblock Compactly supported positive definite radial functions.
\newblock {\em Adv. Comput. Math.}, 4(1):283, 1995.

\bibitem{xu2011point}
Z.~Xu, Y.~Liu, H.~Du, G.~Lin, and C.-W. Shu.
\newblock Point-wise hierarchical reconstruction for discontinuous {G}alerkin
  and finite volume methods for solving conservation laws.
\newblock {\em J. Comput. Phys.}, 230(17):6843--6865, 2011.

\bibitem{yang2009parameter}
M.~Yang and Z.-J. Wang.
\newblock A parameter-free generalized moment limiter for high-order methods on
  unstructured grids.
\newblock In {\em 47th AIAA Aerospace Sciences Meeting}, page 605, 2009.

\bibitem{zerroukat2005monotonic}
M.~Zerroukat, N.~Wood, and A.~Staniforth.
\newblock A monotonic and positive--definite filter for a semi-{L}agrangian
  inherently conserving and efficient {(SLICE)} scheme.
\newblock {\em Quart. J. Roy. Meteor. Soc.,}, 131(611):2923--2936, 2005.

\bibitem{zerroukat2006parabolic}
M.~Zerroukat, N.~Wood, and A.~Staniforth.
\newblock The parabolic spline method ({PSM}) for conservative transport
  problems.
\newblock {\em Int. J. Numer. Meth. Fluids}, 51(11):1297--1318, 2006.

\bibitem{zhang2009third}
Y.-T. Zhang and C.-W. Shu.
\newblock Third order {WENO} scheme on three dimensional tetrahedral meshes.
\newblock {\em Comm. Comput. Phys.}, 5(2-4):836--848, 2009.

\bibitem{zhang1997convergence}
Z.~Zhang and C.~F. Martin.
\newblock Convergence and {G}ibbs' phenomenon in cubic spline interpolation of
  discontinuous functions.
\newblock {\em J. Comput. Appl. Math.}, 87(2):359--371, 1997.

\bibitem{ziou1998edge}
D.~Ziou and S.~Tabbone.
\newblock Edge detection techniques--an overview.
\newblock {\em Pattern Recognition and Image Analysis}, 8:537--559, 1998.

\bibitem{zygmund1959trigonometric}
A.~Zygmund.
\newblock {\em Trigonometric Series}.
\newblock Cambridge Press, 1959.

\end{thebibliography}

\end{document}